\documentclass[11pt]{amsart}
\usepackage{amsmath,amsthm,amsfonts,amssymb}
\usepackage{verbatim}
\usepackage{latexsym}
\usepackage{nicefrac}
\usepackage{color}
\usepackage{ulem}
\usepackage{enumerate}
\usepackage{mdwlist}
\usepackage[british]{babel}
\usepackage[utf8]{inputenc}
\usepackage{soul}
\newcommand{\R}{{\mathbb{R}}}
\newcommand{\C}{{\mathbb{C}}}

\newcommand{\N}{\mathbb{ N}}
\newcommand{\Z}{\mathbb{ Z}}
\newcommand{\T}{\mathbb{ T}}
\newcommand{\leb}{\lambda\raisebox{0,06cm}{\!\!\!\text{\scriptsize$\setminus$}}\,}

\newcommand{\intpar}[4]{\int\limits_{#1}^{#2}#3~\mathrm{d}#4}
\newcommand{\vek}[3]{\begin{pmatrix} #1\\ #2\\ #3\end{pmatrix}}

\newtheorem{theorem}{Theorem}[subsection]
\newtheorem{lemma}[theorem]{Lemma}

\newtheorem{definition}[theorem]{Definition}
\newtheorem{corollary}[theorem]{Corollary}

\newtheorem{example}[theorem]{Example}

\newtheorem{remark}[theorem]{Remark}

\begin{document}
\title{Injectivity conditions for STFT phase retrieval on $\Z$, $\Z_d$ and $\R^d$}
\author{David Bartusel}
\email{david.bartusel@matha.rwth-aachen.de}
\address{Lehrstuhl A f\"ur Mathematik, RWTH Aachen University, D-52056 Aachen, Germany}

\maketitle

\begin{abstract}
We study the phase retrieval problem for the short-time Fourier transform on the groups $\Z$, $\Z_d$ and $\R^d$. As is well-known, phase retrieval is possible, once the windows's ambiguity function vanishes nowhere. However, there are only few results for windows which don't meet this condition.
The goal of this paper is to establish new and complete characterizations for phase retrieval with more general windows and compare them to existing results. For a fixed window, our uniqueness conditions usually only depend on the signal's support and are therefore easily comprehensible.
Additionally, we discuss sharpness of both new and existing results by looking at various examples along the way.
\end{abstract}



\section{Introduction}

We are interested in the phase retrieval problem for the short-time Fourier transform (STFT) on $\Z$, $\Z_d$ and $\R^d$.
More precisely, we try to investigate whether a signal $f\in L^2(G)$ is uniquely determined (up to a global phase factor) by the measurement $|V_gf|$, where $g\in L^2(G)$ is a known window and $G\in\{\Z,\Z_d,\R^d\}$. Here, the STFT is defined by
\begin{align*}
V_gf(k,z)&=\sum\limits_{j\in\Z} f_j\overline{g_{j-k}}z^{-j}\quad\quad(k\in\Z,z\in\C,|z|=1),\\
V_gf(k,l)&=\sum\limits_{j=0}^{d-1} f_j\overline{g_{j-k}}e^{\nicefrac{-2\pi ijl}{d}}\quad\quad(k,l\in\{0,\dots,d-1\})\quad\text{or}\\
V_gf(x,\omega)&=\intpar{\R^d}{}{f(t)\overline{g(t-x)}e^{-2\pi i\langle t,\omega\rangle}}{t}\quad\quad(x,\omega\in\R^d)
\end{align*}
for elements $f$ and $g$ of $\ell^2(\Z)$, $\C^d$ or $L^2(\R^d)$ respectively.\\
Applications of this specific phase retrieval problem include speech processing (cf. e.g. \cite{GL}) and ptychography (cf. e.g. \cite{dSM}). For a more general overview of phase retrieval and its applications, we refer the interested reader to the review paper \cite{GKR}.

For a fixed window $g$, the goal of this paper is to determine which signals $f$ can be recovered from the measurement (up to global phase). The common ground among all our results is that they are quite easy to understand since they usually only involve elementary connectivity conditions on the signal's support.\\
It is already well-known that disconectedness causes non-trivial ambiguities or at least instabilities for phase retrieval, both in the discrete \cite{BCEMS,CHLSS,LNQ} and the continuous setting \cite{ADGY,GRSpecCl}. However, uniqueness for STFT phase retrieval is still far from being completely understood.

The main tool in our approach is the observation that the given phase retrieval problem essentially reduces to the analysis of the so-called \textit{ambiguity functions} $V_ff$ and $V_gg$. More precisely, it is crucial to determine whether the signal $f$ can be recovered from the restriction of $V_ff$ to the support of $V_gg$.\\
This is one of the most prominent insights into STFT phase retrieval and has the well-known consequence that all signals can be uniquely recovered (up to global phase), once the window's ambiguity function vanishes nowhere (see e.g. \cite{BF,GKR,GRSpecCl}).\\
The so-called ambiguity function relation (see Theorem \ref{thm:VgfvsVff} below) has also been used in a more direct way in \cite{AlaiBL}, in order to analyze phase retrieval for windows whose ambiguity functions potentially have much smaller supports.
Here, the authors have proven uniqueness for Paley-Wiener signals for the case of windows that are non-vanishing only on one or two line segments. In contrast, we will generally consider larger signal classes, which (sometimes) comes at the cost of having to be slightly more restrictive on the window.\\
Finally, note that similar results can also be obtained using different proof methods, e.g. for the discrete setting in \cite[Theorem 2.4]{BF} by application of the so-called \textit{PhaseLift operator}.

\subsection{Structure of the paper}

The paper is organized as follows: In section \ref{sec:gen}, we will start out by considering the given problem in a more general setting. Namely, we will recall some basic properties of the short-time Fourier transform on locally compact abelian groups, including most importantly the aforementioned ambiguity function relation (Theorem \ref{thm:VgfvsVff}).

We will then return to the more ``classical'' settings ($\Z$, $\Z_d$ and $\R^d$) in sections \ref{sec:Z}-\ref{sec:Rd}. Phase retrieval on $\Z$ is by far the easiest of the three. Roughly speaking, this is due to the different structures of $\Z$ and its dual group.
At least for one-sided windows, we will be able to provide a full characterization of phase retrieval (Theorem \ref{thm:Zoneside}), which will carry over to general windows under mild assumptions on their decay (Theorem \ref{thm:Zgeneral}).\\
On $\Z_d$ and $\R^d$, the problem of phase retrieval presents itself far more challenging. In both cases, the dual group is isomorphic to the original group, which takes away the advantages of the previous setting.\\
For the finite setting, we will first examine whether a non-vanishing ambiguity function of the window is a necessary condition for (global) phase retrieval, which we will prove to be false in every dimension $d\geq 4$ (Corollary \ref{cor:Vgg}). Afterwards, we will restrict ourselves to short windows.
It turns out that the ambiguity function behaves quite nicely for the majority of these windows, which yields once more a full characterization for phase retrieval (Theorem \ref{thm:Cdgeneric}). For the rest of this window class, the ambiguity function may behave quite pathologically, but we will be able to extend the uniqueness result at least to a certain class of sparse signals (Theorem \ref{thm:CdHoleL}).\\
For the continuous setting, we will finally face a different challenge. While the ambiguity function can be controlled when making certain assumptions on the window, the recovery of the signal can't be performed as easily as within the discrete settings, since most of the information has to be understood in an ``almost everywhere'' sense.
Nevertheless, the results from the infinite-dimensional discrete setting will carry over to the continuous setting, as we will prove a complete characterization of phase retrieval for one-sided windows in dimension one as well as for windows of exponential decay (Theorem \ref{thm:Rdexpdec}), which is perhaps the most intriguing result of the paper.

\subsection{Notation}

We conclude this section by fixing some basic notation.

For $d\geq 2$, let $\Z_d:=\Z/d\Z$ be the cyclic group of order $d$. Throughout section \ref{sec:gen}, we will usually consider $d\in\N$ in order to treat the group $\R^d$ simultaneously. Note however, that we will always assume implicitly that $d\geq 2$ whenever considering $\Z_d$.\\
Any locally compact group with the name $G$ will be equipped with an additive structure. Recall that there exists a Haar measure $\mu_G$ on $G$, which is unique up to a constanct factor. For $p\geq 1$, we will -- with a slight abuse of notation -- always think of $f\in L^p(G)$ as a measurable, $p$-integrable function $f:G\to\C$.
The inner product on $L^2(G)$ is chosen to be linear in the first and conjugate linear in the second argument, i.e.
\begin{equation*}
\langle f,g\rangle:=\int_G f(t)\overline{g(t)}~\mathrm{d}\mu_G(t).
\end{equation*}
For any measurable function $f:G\to\C$ and every $z\in G$, we define the time-shifted function $T_zf:G\to\C,~(T_zf)(x)=f(x-z)$.\\
When $G$ is discrete and $z\in G$, we define $\delta_z:G\to\C$ to equal the characteristic function $\chi_{\{z\}}$.

We will frequently identify $\Z_d\cong\{0,\dots,d-1\}$ and will accordingly write $x\in L^2(\Z_d)\cong\C^d$ as $x=(x_0,\dots,x_{d-1})^t$. However, we will still think of $x\in L^2(\Z_d)$ as a $d$-periodic signal, such that the evaluation $x_j$ is well-defined for every $j\in\Z$.\\
Finally, we define the torus $\T:=\{z\in\C~|~|z|=1\}$, which becomes a (locally) compact abelian group when equipped with the subspace topology and the usual multiplication in $\C$.

\section{STFT phase retrieval on locally compact abelian groups}\label{sec:gen}

Let $G$ be a $\sigma$-compact, locally compact abelian group with Haar measure $\mu_G$. Given a family $\left(\varphi_j\right)_{j\in I}$ of linear functionals on $L^2(G)$, \textit{phase retrieval} tries to recover a signal $f$ from the phaseless measurements $\left(\left|\varphi_j(f)\right|\right)_{j\in I}$.
In other words, one would like to find a large subset $Y\subseteq L^2(G)$ such that the mapping
\begin{equation*}
Y\to\C^{I},\quad f\mapsto\left(\left|\varphi_j(f)\right|\right)_{j\in I}
\end{equation*}
is injective. We are interested in the special case where $\left(\varphi_j\right)_{j\in I}$ is the so-called \textit{short-time Fourier transform} with respect to a given window function $g\in L^2(G)$. Here, we will focus on $G\in\left\{\Z_d,\Z,\R^d\right\}$. Throughout this section however, we will recall some basic properties within the more general setting.

\subsection{The short-time Fourier transform}

For the definition of the Fourier transform on $G$ and its basic properties, we will follow \cite{RudFAG}.\\
Let $\widehat{G}$ denote the dual group of $G$, i.e. the set of all continuous characters $\chi:G\to\T$, with group operation
\begin{equation*}
(\chi_1\cdot\chi_2)(y)=\chi_1(y)\cdot\chi_2(y).
\end{equation*}
From now on, we will make the additional assumption that $G$ is also metrizable, which obviously holds true for $G\in\{\Z,\Z_d,\R^d\}$.
By \cite[Theorem 7.1.4]{Deit}, this ensures that $\widehat{G}$, endowed with the topology of uniform convergence on compact sets, is also a $\sigma$-compact, locally compact abelian group, which allows us in particular to use Fubini's theorem not only for $G$, but also for $\widehat{G}$.\\
We fix a Haar measure $\mu_{\hat{G}}$ on $\hat{G}$. Note that, by the Pontryagin duality theorem, we may always identify $\widehat{\widehat{G}}\cong G$ via
\begin{equation*}
y\mapsto\left(\nu_{y}:\widehat{G}\to\T,~\chi\mapsto\chi(y)\right).
\end{equation*}
Given $h\in L^1(G)$, we may now define its Fourier transform $\widehat{h}\in\mathcal{C}_0\left(\widehat{G}\right)$ by
\begin{equation*}
\widehat{h}(\chi):=\int_G h(y)\overline{\chi(y)}~\mathrm{d}\mu_G(y)\quad\text{for every } \chi\in\widehat{G}.
\end{equation*}
Sometimes, it will be more convenient to write $\mathcal{F}_G(h)$ instead of $\widehat{h}$.

We will begin by recalling the inversion formula and Plancherel's theorem.

\begin{theorem}\label{thm:planch}
The following statements hold true.
\begin{enumerate}[a)]
\item There exists a Haar measure $\mu^{\ast}_{\widehat{G}}$ on $\widehat{G}$ such that for every $h\in L^1(G)$ satisfying $\widehat{h}\in L^1\left(\widehat{G}\right)$, it holds
\begin{equation*}
h(y)=\left(\mathcal{F}_{\left(\widehat{G},\mu^{\ast}_{\widehat{G}}\right)}\left(\overline{\nu_y}\right)\right)\left(\widehat{h}\right)=\int_{\widehat{G}} \widehat{h}(\chi)\chi(y)~\mathrm{d}\mu^{\ast}_{\widehat{G}}(\chi)
\end{equation*}
for $\mu_G$-almost every $y\in G$.
\item For every $h\in L^1(G)\cap L^2(G)$, it holds $\widehat{h}\in L^2\left(\widehat{G}\right)$.
\item $\mathcal{F}_G$ extends to an isometry $\left(L^2(G),\mu_G\right)\to\left(L^2\left(\widehat{G}\right),\mu^{\ast}_{\widehat{G}}\right)$, satisfying Parseval's formula
\begin{equation*}
\left\langle\mathcal{F}_G h,\mathcal{F}_G \tilde{h}\right\rangle=\left\langle h,\tilde{h}\right\rangle
\end{equation*}
for all $h,\tilde{h}\in L^2(G)$.
\end{enumerate}
\end{theorem}

Note that the inversion formula (part a)) does not appear in this form in \cite{RudFAG}. However, it is proved in \cite[Theorem 3.5.8]{DE}.

Since two Haar measures on $\widehat{G}$ coincide up to a constant factor, we may slightly reformulate Theorem \ref{thm:planch}. This allows us to freely choose the Haar measure on $\widehat{G}$.

\begin{corollary}\label{cor:planch}
There exists a constant $c_G=c\left(G,\mu_G\right)>0$ satisfying the following conditions.
\begin{enumerate}[a)]
\item For every $h\in L^1(G)$ satisfying $\widehat{h}\in L^1\left(\widehat{G}\right)$, it holds
\begin{equation*}
h(y)=c_G^{-1}\cdot\left(\mathcal{F}_{\widehat{G}}(\overline{\nu_y})\right)\left(\widehat{h}\right)=c_G^{-1}\cdot\int_{\widehat{G}} \widehat{h}(\chi)\chi(y)~\mathrm{d}\mu_{\widehat{G}}(\chi)
\end{equation*}
for $\mu_G$-almost every $y\in G$.
\item $\mathcal{F}_G$ extends to a bounded linear operator $L^2(G)\to L^2\left(\widehat{G}\right)$, satisfying
\begin{equation*}
\left\langle\mathcal{F}_G h,\mathcal{F}_G \tilde{h}\right\rangle=c_G\cdot\left\langle h,\tilde{h}\right\rangle
\end{equation*}
for all $h,\tilde{h}\in L^2(G)$. In particular, it holds
\begin{equation*}
\left\|\widehat{h}\right\|_2=c_G^{\nicefrac{1}{2}}\cdot\|h\|_2
\end{equation*}
for every $h\in L^2(G)$.
\end{enumerate}
\end{corollary}

\begin{proof}
This follows from Theorem \ref{thm:planch} when choosing $c_G:=\frac{\mu_{\widehat{G}}}{\mu_{\widehat{G}}^{\ast}}$ which is a constant.
\end{proof}

When $\widehat{h}\notin L^1\left(\widehat{G}\right)$, the inversion formulas in Theorem \ref{thm:planch} and Corollary \ref{cor:planch} can't be used to obtain $h$. However, uniqueness can still be assured. We refer again to \cite[Corollary 3.3.4]{DE} for the proof.

\begin{lemma}\label{lm:fourinj}
The Fourier transform $\mathcal{F}:L^1(G)\to \mathcal{C}_0\left(\widehat{G}\right)$ is injective.
\end{lemma}

For continuous functions, we obtain a pointwise version of Lemma \ref{lm:fourinj}.

\begin{corollary}\label{cor:fourinj}
Let $h,\tilde{h}:G\to\C$ be continuous and integrable. Then, $h(z)=\tilde{h}(z)$ holds for every $z\in G$ if and only if $\widehat{h}(\chi)=\widehat{\tilde{h}}(\chi)$ holds for every $\chi\in\widehat{G}$.
\end{corollary}

\begin{proof}
Since the ``only if'' part is trivial, suppose that $\widehat{h}(\chi)=\widehat{\tilde{h}}(\chi)$ holds for every $\chi\in\widehat{G}$. By Lemma \ref{lm:fourinj}, it follows $h=\tilde{h}$ almost everywhere on $G$. By continuity of $h$ and $\tilde{h}$, we obtain $h(z)=\tilde{h}(z)$ for every $z\in G$.
\end{proof}

Since the STFT of a signal $f\in L^2(G)$ will be defined on $G\times\widehat{G}$, the following results concerning the product of two groups will be useful.

\begin{lemma}\label{lm:G1xG2}
Let $G_1$ and $G_2$ be metrizable, $\sigma$-compact, locally compact abelian groups with Haar measures $\mu_{G_1}$ and $\mu_{G_2}$ respectively, and choose $\mu_{G_1\times G_2}:=\mu_{G_1}\otimes\mu_{G_2}$ as Haar measure on $G_1\times G_2$. Then, the following statements hold true.
\begin{enumerate}[a)]
\item It holds $\widehat{G_1\times G_2}\cong \widehat{G_1}\times\widehat{G_2}$. More precisely, $\chi\in\widehat{G_1\times G_2}$ is equivalent to the existence of $\alpha\in\widehat{G_1}$ and $\beta\in\widehat{G_2}$ such that
\begin{equation*}
\chi(y,z)=\alpha(y)\beta(z)\quad\text{for every } (y,z)\in G_1\times G_2.
\end{equation*}
Consequently, $\mu_{\widehat{G_1\times G_2}}:=\mu_{\widehat{G_1}}\otimes\mu_{\widehat{G_2}}$ can be chosen as Haar measure on $\widehat{G_1\times G_2}$.
\item For every $h\in L^1(G_1\times G_2)$ and every $\beta\in G_2$, it holds $h_\beta\in L^1(G_1)$, where
\begin{equation*}
h_{\beta}:G_1\to\C,\quad y\mapsto \mathcal{F}_{G_2}\left(h(y,\cdot)\right)(\beta).
\end{equation*}
Moreover, for every $\alpha\in G_1$, it follows
\begin{equation*}
\widehat{h}(\alpha,\beta)=\mathcal{F}_{G_1}\left(h_{\beta}\right)(\alpha).
\end{equation*}
\item It holds
\begin{equation*}
c_{G_1\times G_2}=c_{G_1}\cdot c_{G_2}.
\end{equation*}
\end{enumerate}
\end{lemma}

\begin{proof}\mbox{}
\begin{enumerate}[a)]
\item For $\alpha\in\widehat{G_1}$ and $\beta\in\widehat{G_2}$, it is easy to show that $\chi(y,z):=\alpha(y)\beta(z)$ defines an element $\chi\in\widehat{G_1\times G_2}$.\\
Conversely, every $\chi\in\widehat{G_1\times G_2}$ allows the factorization $\chi(y,z)=\alpha(y)\beta(z)$, where\linebreak $\alpha(y):=\chi(y,0)$ and $\beta(z):=\chi(0,z)$ clearly define elements $\alpha\in\widehat{G_1}$ and $\beta\in\widehat{G_2}$.
\item and c) follow as an application of Fubini's and Tonelli's theorem.
\end{enumerate}
\end{proof}

We will now consider a window $g\in L^2(G)$ and define
\begin{equation*}
\left(\pi(z,\chi)g\right)(y):=\chi(y)\cdot T_zg(y)
\end{equation*}
for every $(z,\chi)\in G\times\widehat{G}$ and every $y\in G$.\\
Given a signal $f\in L^2(G)$, the short-time Fourier transform $V_gf:G\times\widehat{G}\to\C$ is defined by
\begin{equation*}
V_gf(z,\chi)=\left\langle f,\pi(z,\chi)g\right\rangle=\int_G f(y)\overline{\chi(y) g(y-z)}~\mathrm{d}\mu_G(y)=\mathcal{F}_G\left(f\cdot \overline{T_z g}\right)(\chi)
\end{equation*}
for every $(z,\chi)\in G\times\widehat{G}$. (Note that $\pi(\chi,z) g\in L^2(G)$ and $f\cdot \overline{T_z g}\in L^1(G)$.)

It turns out that the STFT is again square-integrable on $G\times\widehat{G}$. As a byproduct of the proof, we will obtain the highly useful orthogonality relations.

\begin{theorem}\label{thm:orthrel}
Let $f,g,\tilde{f},\tilde{g}\in L^2(G)$. Then, $V_gf,V_{\tilde{g}}\tilde{f}\in L^2\left(G\times\widehat{G}\right)$ and it holds
\begin{equation*}
\left\langle V_gf,V_{\tilde{g}}\tilde{f}\right\rangle=c_G\cdot\left\langle f,\tilde{f}\right\rangle\cdot\left\langle \tilde{g},g\right\rangle.
\end{equation*}
\end{theorem}

\begin{proof}
A proof for the case $G=\R^d$ can be found in \cite{Groech}. Because of Plancherel's theorem (Corollary \ref{cor:planch}b)), the same reasoning applies to the general case.
\end{proof}

To conclude this section, we will establish an inversion formula for the STFT.

\begin{theorem}\label{thm:STFTinv}
Let $f,g\in L^2(G)$ For every $\eta\in L^2(G)$ satisfying $\left\langle \eta,g\right\rangle\neq 0$, it holds
\begin{equation*}
f=\frac{1}{c_G\cdot\left\langle\eta,g\right\rangle}\cdot\int_G\int_{\widehat{G}} V_gf(z,\chi)\pi(z,\chi)\eta~\mathrm{d}\mu_{\widehat{G}}(\chi)\mathrm{d}\mu_G(z),
\end{equation*}
in the sense of a vector-valued integral, i.e.
\begin{equation*}
\left\langle f,h\right\rangle=\frac{1}{c_G\cdot\left\langle\eta,g\right\rangle}\int_G\int_{\widehat{G}} \left\langle V_gf(z,\chi)\pi(z,\chi)\eta,h\right\rangle~\mathrm{d}\mu_{\widehat{G}}(\chi)\mathrm{d}\mu_G(z)
\end{equation*}
for every $h\in L^2(G)$.
\end{theorem}

\begin{proof}
Once again, the proof in \cite{Groech} (for $G=\R^d$) also applies to the general case.
\end{proof}

\subsection{Phase retrieval and ambiguity functions}

From now on, let $g\in L^2(G)$ be a fixed window. Phase retrieval tries to recover a signal $f\in L^2(G)$ from the phaseless measurement $|V_g f|$. Since the STFT is linear in $f$, we obviously obtain $|V_g f|=|V_g (\gamma f)|$ for any $\gamma\in\T$. Thus, there is only a chance to recover $f$ up to a global phase factor.
With that in mind, we call $f$ \textit{phase retrievable} w.r.t. $g$, iff the statement
\begin{equation*}
\left|V_g f\right|=\left|V_g \tilde{f}\right|\quad\Rightarrow\quad \tilde{f}=\gamma f~\text{for some }\gamma\in\T
\end{equation*}
holds true for every $\tilde{f}\in L^2(G)$.\\
We say that the window $g$ \textit{does phase retrieval} iff every $f\in L^2(G)$ is phase retrievable w.r.t. $g$.\\

We make a first simple observation regarding shifted and reflected windows.

\begin{lemma}\label{lm:windowshift}
Let $f,g\in L^2(G)$ as well as $y\in G$. The following statements are equivalent.
\begin{enumerate}[(i)]
\item $f$ is phase retrievable w.r.t. $g$.
\item $f$ is phase retrievable w.r.t. $T_yg$.
\item $\mathcal{I}f$ is phase retrievable w.r.t. $\mathcal{I}g$, where $\mathcal{I}:G\to G,~(\mathcal{I}h)(x)=h(-x)$.
\end{enumerate}
\end{lemma}

\begin{proof}
For any $(z,\chi)\in G\times\widehat{G}$, it holds
\begin{equation*}
\pi(z,\chi)T_yg=\chi\cdot T_zT_yg=\chi\cdot T_{z+y}g=\pi(z+y,\chi)g,
\end{equation*}
which implies $\left|V_{T_yg}f(z,\chi)\right|=\left|V_gf(z+y,\chi)\right|$, and thus the equivalence of (i) and (ii).\\
Additionally, we obtain
\begin{equation*}
V_{\mathcal{I}g}\mathcal{I}f(z,\chi)=\int_G f(-y)\overline{\chi(y)g(z-y)}~\mathrm{d}\mu_G=\int_G f(y)\overline{\chi(-y)g(y-(-z))}~\mathrm{d}\mu_G(y)=V_gf(-z,\overline{\chi}),
\end{equation*}
which establishes the equivalence of (i) and (iii).
\end{proof}

Within this subsection, we will see that phase retrieval w.r.t $g$ only depends on the set
\begin{equation*}
\Omega(g):=\left\{(z,\chi)\in G\times\widehat{G}~\middle|~V_gg(z,\chi)\neq 0\right\},
\end{equation*}
where $V_gg$ is called the \textit{ambiguity function} of $g$. (Note that for $G=\R^d$, the term ``ambiguity function'' usually refers to $A_g(x,\chi)=\chi\left(\frac{x}{2}\right)\cdot V_gg(x,\chi)$ in the literature. However, since both functions obviously have the same zeros, this won't cause any trouble.)\\
More precisely, we will show the following lemma.

\begin{lemma}\label{lm:Vff}
A signal $f\in L^2(G)$ is phase retrievable w.r.t. $g$, iff the statement
\begin{equation*}
V_ff\vert_{\Omega(g)}=V_{\tilde{f}}\tilde{f}\vert_{\Omega(g)}\quad\Rightarrow\quad\tilde{f}=\gamma f\text{ for some } \gamma\in\T
\end{equation*}
holds true for every $\tilde{f}\in L^2(G)$.
\end{lemma}

This approach is rather common and leads to the perhaps most prominent uniqueness condition in STFT phase retrieval (see Theorem \ref{thm:Vggnonvan} below), which states that every signal is phase retrievable w.r.t. $g$, whenever it holds $\Omega(g)=G\times\widehat{G}$.
The theorem is well-known and has appeared in various settings, cf. e.g. \cite{BF,GKR,GRSpecCl}.\\
However, Lemma \ref{lm:Vff} can also be used more directly, by choosing the window in such a way that $V_gg$ doesn't vanish at least on some (known) subset of $G\times\widehat{G}$. In this case, the goal of phase retrieval reduces to finding a (possibly large) class of signals that are phase retrievable w.r.t. $g$.
For instance, this strategy has been used in \cite{AlaiBL} to provide phase retrieval results for Paley-Wiener functions, requiring only mild conditions on the window.\\

Generally speaking, the main goals of this paper can be summarized as follows:
\begin{enumerate}[1)]
\item Provide conditions on a window $g$ that allow to characterize the corresponding set $\Omega(g)$.
\item Based on the knowledge of $\Omega(g)$, provide conditions on a signal $f$ that allow to decide whether $f$ is phase retrievable w.r.t. $g$.
\end{enumerate}
Moreover, we intend to give these conditions mostly in terms of the supports of $f$ and $g$, in order for the ensuing results to be easily understandable.\\

The proof of Lemma \ref{lm:Vff} crucially involves the Fourier transform of $|V_gf|^2\in L^1\left(G\times\widehat{G}\right)$. By Lemma \ref{lm:G1xG2} and the Pontryagin duality theorem, we may identify
\begin{equation*}
\widehat{G\times\widehat{G}}\cong \widehat{G}\times\widehat{\widehat{G}}\cong\widehat{G}\times G.
\end{equation*}
More precisely, every character $\nu\in\widehat{G\times\widehat{G}}$ is of the form
\begin{equation*}
\nu_{\chi,z}:G\times\widehat{G}\to\T,\quad \nu_{\chi,z}(y,\alpha)=\chi(y)\alpha(z)
\end{equation*}
for some $\chi\in\widehat{G}$ and $z\in G$. For convenience, we will from now on always use the notation
\begin{equation*}
\widehat{h}(\chi,z):=\widehat{h}\left(\nu_{\chi,z}\right)=\int_{\widehat{G}}\int_G h(y,\alpha)\overline{\chi(y)\alpha(z)}~\mathrm{d}\mu_{\widehat{G}}(\alpha)\mathrm{d}\mu_G(\alpha),
\end{equation*}
where $(\chi,z)\in\widehat{G}\times G$.

Using the results of the previous section, we can now prove the following important identity.

\begin{theorem}\label{thm:VgfvsVff}
For every $f\in L^2(G)$ and every $(\chi,z)\in\widehat{G}\times G$, it holds
\begin{equation*}
\widehat{\left|V_gf\right|^2}(\chi,z)=c_G\cdot V_ff(-z,\chi)\cdot\overline{V_gg(-z,\chi)}.
\end{equation*}
\end{theorem}

\begin{proof}
There are (slightly varying) proofs of this theorem in different settings, that can be found e.g. in \cite{AlaiBL,GKR,GRSpecCl}. For convenience, we will give the proof for the general setting, even though it is essentially the same as for the finite setting in \cite{GKR}:\\
Fix $f\in L^2(G)$ and $(\chi,z)\in\widehat{G}\times G$. For every $(y,\alpha)\in G\times\widehat{G}$, we compute
\begin{align*}
V_gf(y,\alpha)\cdot\overline{\chi(y)\alpha(z)}	&=\int_G f(s)\overline{\alpha(s)g(s-y)}\cdot\overline{\chi(y)\alpha(z)}~\mathrm{d}\mu_G(s)\\
							&=\int_G f(s-z)\cdot\overline{\alpha(s-z)\alpha(z)}\cdot\chi(-y)\overline{g(s-z-y)}~\mathrm{d}\mu_G(s)\\
							&=\int_G \overline{\chi(s)}f(s-z)\cdot\overline{\alpha(s)}\cdot\chi(s-y)\overline{g(s-z-y)}~\mathrm{d}\mu_G(s)\\
							&=\left(V_{\left(\pi(z,\overline{\chi})g\right)}\left(\pi(z,\overline{\chi})f\right)\right)(y,\alpha).
\end{align*}
This may be used together with the orthogonality relations (Theorem \ref{thm:orthrel}) to obtain
\begin{align*}
\widehat{\left|V_gf\right|^2}(\chi,z)	&=\int_{\widehat{G}}\int_G V_gf(y,\alpha)\cdot\overline{V_gf(y,\alpha)}\cdot\overline{\chi(y)\alpha(z)}~\mathrm{d}\mu_G(y)\mathrm{d}\mu_{\widehat{G}}(\alpha)\\
						&=\int_{\widehat{G}}\int_G \left(V_{\left(\pi(z,\overline{\chi})g\right)}\left(\pi(z,\overline{\chi})f\right)\right)(y,\alpha)\cdot\overline{V_gf(y,\alpha)}~\mathrm{d}\mu_G(y)\mathrm{d}\mu_{\widehat{G}}(\alpha)\\
						&=\left\langle V_{\left(\pi(z,\overline{\chi})g\right)}\left(\pi(z,\overline{\chi})f\right),V_gf\right\rangle\\
						&=c_G\cdot\left\langle\pi(z,\overline{\chi})f,f\right\rangle\cdot\left\langle g,\pi(z,\overline{\chi})g\right\rangle.
\end{align*}
Finally, we compute
\begin{align*}
\left\langle \pi(z,\overline{\chi})f,f\right\rangle	&=\int_G \overline{\chi(y)}f(y-z)\overline{f(y)}~\mathrm{d}\mu_G(y)\\
								&=\int_G \overline{\chi(y+z)}f(y)\overline{f(y+z)}~\mathrm{d}\mu_G(y)\\
								&=\overline{\chi(z)}\cdot\int_G f(y)\overline{\chi(y)f(y+z)}~\mathrm{d}\mu_G(y)\\
								&=\overline{\chi(z)}\cdot V_ff(-z,\chi),
\end{align*}
as well as
\begin{align*}
\left\langle g,\pi(z,\overline{\chi})g\right\rangle	&=\int_G g(y)\chi(y)\overline{g(y-z)}~\mathrm{d}\mu_G(y)\\
								&=\int_G g(y+z)\chi(y+z)\overline{g(y)}~\mathrm{d}\mu_G(y)\\
								&=\chi(z)\cdot\overline{\int_G g(y)\overline{\chi(y)g(y+z)}~\mathrm{d}\mu_G(y)}\\
								&=\chi(z)\cdot\overline{V_gg(-z,\chi)}.
\end{align*}
The fact that $\chi(z)\cdot\overline{\chi(z)}=1$ entails the identity.
\end{proof}

With Theorem \ref{thm:VgfvsVff} at hand, it will now be easy to prove Lemma \ref{lm:Vff}.

\begin{proof}[Proof of Lemma \ref{lm:Vff}]
According to Corollary \ref{cor:fourinj} (applied to $G\times\widehat{G}$) and Theorem \ref{thm:VgfvsVff}, it holds
\begin{align*}
\left|V_gf\right|=\left|V_g\tilde{f}\right|\quad	&\Leftrightarrow\quad\widehat{\left|V_gf\right|^2}=\widehat{\left|V_g\tilde{f}\right|^2}\\
							&\Leftrightarrow\quad\forall (z,\chi)\in G\times\widehat{G}:~ \widehat{\left|V_gf\right|^2}(\chi,-z)=\widehat{\left|V_g\tilde{f}\right|^2}(\chi,-z)\\
							&\Leftrightarrow\quad\forall (z,\chi)\in \Omega(G):~ V_ff(z,\chi)=V_{\tilde{f}}\tilde{f} (k,l),
\end{align*}
which implies the statement.
\end{proof}

As mentioned before, this yields a simple sufficient condition for the window $g$ to do phase retrieval. Once again, the known proofs differ slightly between various settings. We will opt for a proof that directly uses the STFT inversion formula.

\begin{theorem}\label{thm:Vggnonvan}
Suppose that $g\in L^2(G)$ satisfies $\Omega(g)=G\times\widehat{G}$. Then $g$ does phase retrieval.
\end{theorem}

\begin{proof}
By Lemma \ref{lm:Vff}, it suffices to show that for all $f,\tilde{f}\in L^2(G)$, the condition $V_ff=V_{\tilde{f}}\tilde{f}$ implies the existence of $\gamma\in\T$ satisfying $\tilde{f}=\gamma f$. W.l.o.g. we assume $f,\tilde{f}\not\equiv 0$.\\
Choosing
\begin{equation*}
\eta:=\begin{cases} f, &\text{if } \left\langle f,\tilde{f}\right\rangle\neq 0,\\ f+\tilde{f},&\text{otherwise},\end{cases}
\end{equation*}
ensures $\eta\in L^2(G)$ as well as $\left\langle\eta,f\right\rangle\neq 0\neq\left\langle\eta,\tilde{f}\right\rangle$ and thus, the inversion formula in Theorem \ref{thm:STFTinv} yields
\begin{equation*}
f\cdot\left\langle\eta,f\right\rangle=\tilde{f}\left\langle\eta,\tilde{f}\right\rangle.
\end{equation*}
Hence, it holds $\tilde{f}=c f$, where $c:=\frac{\left\langle\eta,f\right\rangle}{\left\langle\eta,\tilde{f}\right\rangle}$. It remains to show that $c\in\T$, but this follows from the fact that
\begin{equation*}
\widehat{\left|f\right|^2}=V_ff(0,\cdot)=V_{\tilde{f}}\tilde{f}(0,\cdot)=\widehat{\left|\tilde{f}\right|^2},
\end{equation*}
which implies $|f|=\left|\tilde{f}\right|$ by Lemma \ref{lm:fourinj}.
\end{proof}

It turns out that we won't even have to check the condition $V_gg(z,\chi)\neq 0$ for every $(z,\chi)\in G\times\widehat{G}$. Being the Fourier transform of $g\cdot \overline{T_z g}$, we would expect $V_gg(z,\cdot)$ to carry the same information as $V_gg(-z,\cdot)$, and the same for $f$. This intuition is justified by the following statement.

\begin{lemma}\label{lm:sym}
For every $h\in L^2(G)$ and every $(z,\chi)\in G\times\widehat{G}$, it holds
\begin{equation*}
V_hh(-z,\overline{\chi})=\overline{\chi(z)\cdot V_hh(z,\chi)}.
\end{equation*}
In particular, $(z,\chi)\in\Omega(h)$ holds if and only if $(-z,\overline{\chi})\in\Omega(h)$.
\end{lemma}

\begin{proof}
We compute
\begin{align*}
V_hh(-z,\overline{\chi})&=\int_G h(y)\chi(y)\overline{h(y+z)}~\mathrm{d}\mu_G(y)\\&=\int_G h(y-z)\chi(y-z)\overline{h(y)}~\mathrm{d}\mu_G(y)\\&=\overline{\chi(z)}\cdot\overline{\int_G h(y)\overline{\chi(y)h(y-z)}~\mathrm{d}\mu_G(y)}\\&=\overline{\chi(z)\cdot V_hh(z,\chi)}.
\end{align*}
\end{proof}

We will conclude this section by giving a brief outline of our principal recovery strategy. The main idea is already contained in the proof of Theorem \ref{thm:Vggnonvan} and has been featured prominently (at least for discrete settings), i.e. in \cite{GKR,LNQ}.

Suppose that for some $z\in G$, it holds $V_gg(z,\chi)\neq 0$ for all $\chi\in\widehat{G}$. Since phase retrieval is equivalent to recovering $f$ up to a global phase from $V_ff\vert_{\Omega(g)}$ by Theorem \ref{thm:VgfvsVff}, we may then assume to have access to
\begin{equation*}
V_ff(z,\cdot)=\mathcal{F}_G\left(f\cdot \overline{T_z f}\right)
\end{equation*}
and thus to $f\cdot \overline{T_z f}$ by Lemma \ref{lm:fourinj}. Sometimes (most notably when $\widehat{G}\in\{\T,\R^d\}$), it might even be enough to assume $V_ff(z,\chi)\neq 0$ for only almost every $\chi\in\widehat{G}$.\\
If we have access to this information for $z=0$, we obtain $|f|$ almost everywhere. (In many settings, this assumption can easily be ensured.) Intuitively, phase retrieval may then be performed using the following two steps.
\begin{enumerate}[1)]
\item Fix a phase for $f$ at some $y\in G$ satisfying $f(y)\neq 0$.
\item Propagate the phase along $G$ using the information about $f\cdot \overline{T_z f}$ for various $z\neq 0$.
\end{enumerate}
Obviously, we will lose the phase information when coming across zeros of $f$ within step 2). At this point, it becomes clear that certain connectivity assumptions on $\operatorname{supp}(f)$ can be a very useful tool for phase retrieval. Another serious problem arises from the fact that we only have information on $f\cdot \overline{T_z f}$ \textit{almost} everywhere.
However, this second problem does not occur when $G$ is discrete, which is why we will start out by considering this case in the following two sections.

\section{STFT phase retrieval on $\Z$}\label{sec:Z}

In this section, we analyze phase retrieval on the group $G=\Z$. Since we may identify $\widehat{G}\cong\T$, with characters
\begin{equation*}
\chi_z:\Z\to\T,\quad j\mapsto z^j,
\end{equation*}
the Fourier transform becomes
\begin{equation*}
\widehat{h}(z)=\sum\limits_{j\in\Z} h_j z^{-j}\quad (z\in\T),
\end{equation*}
whereas the STFT is given by
\begin{equation*}
V_gf(k,z)=\sum\limits_{j\in\Z} f_j\overline{g_{j-k}}z^{-j}\quad(k\in\Z,z\in\T).
\end{equation*}
The Haar measures on $\Z$ and $\T$ are the counting measure and the spherical measure $\sigma$ respectively, where the latter is normed such that $\sigma(\T)=1$.

\subsection{Uniqueness conditions for phase retrieval on $\Z$}

The main benefit of considering STFT phase retrieval on $\Z$ can be described as follows: Since the time-frequency plane $G\times\widehat{G}$ is given by $\Z\times\T$, we basically try to recover a countable number of values from an uncountable number of measurements.
Recall that our main strategy is to recover $f\cdot \overline{T_kf}$ from $\{V_ff(k,z)~|~(k,z)\in\Omega(g)\}$ for any given $k\in\Z$ by applying injectivity of the Fourier transform. Once $N:=\{z\in\T~|~(k,z)\in\Omega(g)\}$ is of measure zero, we can find $V_ff(k,\cdot)$ by continuity and thus gain full access to $f\cdot\overline{T_kf}$.
The reconstruction can also be expressed in a more ``computational'' way: Since $\T$ is compact, it holds $\hat{h}\in L^1(\T)$ for every $h\in L^1(\Z)=\ell^1(\Z)$ and we can apply the Fourier inversion formula to $V_ff(k,\cdot)$, which doesn't require knowledge of $V_ff(k,z)$ for $z\in N$.

For convenience, we will compute the constant $c_G$.

\begin{lemma}\label{lm:cGZ}
When $G=\Z$, it holds $c_G=1$.
\end{lemma}

\begin{proof}
By Corollary \ref{cor:planch} b), it suffices to compute $\frac{\left\|\widehat{h}\right\|_2}{\|h\|_2}$ for a single $h\neq 0$. Choosing $h=\delta_0$ yields $\widehat{h}(z)=1$ for every $z\in\T$. Hence, $\left\|\widehat{h}\right\|_2=\sigma^1(\T)=1=\|h\|_2$, which entails $c_G=1$.
\end{proof}

As mentioned above, Fourier inversion (or simply injectivity) allows us to prove the following statement which slightly generalizes Theorem \ref{thm:Vggnonvan}.

\begin{corollary}\label{cor:aenonvanZ}\mbox{}
\begin{enumerate}[a)]
\item Let $f,\tilde{f}\in\ell^2(\Z)$ as well as $k\in\Z$. If $V_ff(k,z)=V_{\tilde{f}}\tilde{f}(k,z)$ holds for almost every $z\in\T$, it follows $f\cdot \overline{T_k f}=\tilde{f}\cdot \overline{T_k\tilde{f}}$.
\item If $g\in\ell^2(\Z)$ satisfies $V_gg(k,z)\neq 0$ for almost every $(k,z)\in\Z\times\T$, it follows that $g$ does phase retrieval.
\end{enumerate}
\end{corollary}

\begin{proof}\mbox{}
\begin{enumerate}[a)]
\item Since $V_ff(k,\cdot)$ and $V_{\tilde{f}}\tilde{f}(k,\cdot)$ are both continuous on $\T$ and are therefore elements of $L^1(\T)$, we obtain
\begin{equation*}
f\cdot \overline{T_k f}=\int_{\T} V_ff(k,z)\cdot z^j~\mathrm{d}\sigma^1(z)=\int_{\T} V_{\tilde{f}}\tilde{f}(k,z)\cdot z^j~\mathrm{d}\sigma^1(z)=\tilde{f}\cdot \overline{T_k\tilde{f}}
\end{equation*}
almost everywhere on $\Z$ by Corollary \ref{cor:planch}. Since there are no non-empty sets of measure zero in $\Z$, it follows $f\cdot \overline{T_k f}=\tilde{f}\cdot \overline{T_k\tilde{f}}$.
\item Since
\begin{equation*}
\mu_{\Z\times\T}\left(\left\{(k,z)\in \Z\times\T~\middle|~V_gg(k,z)=0\right\}\right)=\sum\limits_{k\in\Z} \sigma\left(\left\{z\in\T~\middle|~V_gg(k,z)=0\right\}\right),
\end{equation*}
the assumption implies that the condition
\begin{equation*}
V_gg(k,z)\neq 0\quad\text{for almost every } z\in\T
\end{equation*}
holds true for every $k\in\Z$. Thus, when $f,\tilde{f}\in\ell^2(\Z)$ are such that $V_ff\vert_{\Omega(g)}=V_{\tilde{f}}\tilde{f}\vert_{\Omega(g)}$, it follows $f\cdot \overline{T_k f}=\tilde{f}\cdot \overline{T_k \tilde{f}}$ for every $k\in\Z$ because of part a).
Taking this back to the Fourier domain yields $V_ff=V_{\tilde{f}}\tilde{f}$, which finally implies $\tilde{f}=\gamma f$ for some $\gamma\in\T$, just as in the proof of Theorem \ref{thm:Vggnonvan}. Therefore, $g$ does phase retrieval by Lemma \ref{lm:Vff}.
\end{enumerate}
\end{proof}

The assumption of Corollary \ref{cor:aenonvanZ} a) is easy to verify when considering windows of finite length $L+1$, i.e. $\operatorname{supp}(g)=\{0,\dots,L\}$. According to Lemma \ref{lm:windowshift}, this particular choice of support is purely conventional and can be shifted freely along $\Z$.\\
For this kind of window, we obtain
\begin{equation*}
V_gg(k,z)=\sum\limits_{j=k}^d g_j\overline{g_{j-k}}z^{-j}=z^{-k}\cdot\sum\limits_{j=0}^{d-k}g_{j+k}\overline{g_j}\overline{z}^{j}
\end{equation*}
Since the polynomial $\sum_{j=0}^{d-k} g_{j+k}\overline{g_j} z^j$ doesn't vanish identically for $0\leq k\leq d$, it can have at most $d-k$ zeros (on $\T$).

Therefore, Corollary \ref{cor:aenonvanZ} implies that the condition $V_ff\vert_{\Omega(g)}=V_{\tilde{f}}\tilde{f}\vert_{\Omega(g)}$ yields $f\cdot \overline{T_k f}=\tilde{f}\cdot \overline{T_k\tilde{f}}$ for every $0\leq k\leq L$.

As suggested at the end of section \ref{sec:gen}, we can now follow the ``initialization and propagation'' approach and as a consequence, phase retrieval purely reduces to a question of connectivity.

\begin{definition} Consider a subset $M\subseteq\Z$ and a signal $f\in\ell^2(\Z)$.
\begin{enumerate}[a)]
\item Two indices $j,k\in M$ are called $L$-connected and we write $j\sim_{\Z,L}k$, iff there exist indices $l_0,\dots,l_n\in M$ satisfying $l_0=j$, $l_n=k$ and $|l_{m+1}-l_m|\leq L$ for every $0\leq m<n$.
\item Clearly, $\sim_{\Z,L}$ defines an equivalence relation on $M$. The equivalence classes of $M$ are called $L$-connectivity components of $M$, and $M$ is called $L$-connected, iff there exists at most one connectivity component.
\item $f$ is called $L$-connected, iff $\operatorname{supp}(f)$ is $L$-connected. The $L$-connectivity components of $\operatorname{supp}(f)$ are called $L$-connectivity components of $f$.
\end{enumerate}
\end{definition}

A set (or a signal) is thus $L$-connected, iff it is connected up to ``holes'' of length at most $L-1$.
With the observations above, it is easy to prove the following theorem. Note that the necessary conditions have already been observed, e.g. in \cite{LNQ}.

\begin{theorem}\label{thm:ZlengthL+1}
Let $g\in\ell^2(\Z)$ be such that $\operatorname{supp}(g)=\{0,\dots,L\}$. Furthermore, consider a signal $f\in\ell^2(\Z)$ and let $\left(C_m\right)_{m\in\N}$ be the $L$-connectivity components of $f$. Then, the following statements hold true.
\begin{enumerate}[a)]
\item A signal $\tilde{f}\in\ell^2(\Z)$ satisfies $\left|V_gf\right|=\left|V_g\tilde{f}\right|$, if and only if $\operatorname{supp}\left(\tilde{f}\right)=\operatorname{supp}(f)$ and for every $m\in\N$, there exists $\gamma_m\in\T$ such that $\tilde{f}\vert_{C_m}=\gamma_mf\vert_{C_m}$.
\item $f$ is phase retrievable w.r.t. $g$ if and only if $f$ is $L$-connected.
\end{enumerate}
\end{theorem}

\begin{proof}\mbox{}
\begin{enumerate}[a)]
\item Let $\tilde{f}\in\C^d$ and suppose that $\operatorname{supp}\left(\tilde{f}\right)=\operatorname{supp}(f)$ as well as $\tilde{f}\vert_{C_m}=\gamma_m f\vert_{C_m}$ for every $m\in\N$. Furthermore, let $(k,z)\in\Omega(g)$, i.e. in particular $-L\leq k\leq L$. Consequently, for every $j\in\Z$, it holds either
\begin{equation*}
\tilde{f}\overline{\tilde{f}_{j-k}}=f_j\overline{f_{j-k}}=0
\end{equation*}
or $j$ and $j-k$ belong to the same connectivity component $C_m$ of $f$. In this case, it follows
\begin{equation*}
\tilde{f}_j\overline{\tilde{f}_{j-k}}=\gamma_mf_j\overline{\gamma_m}\overline{f_{j-k}}=f_j\overline{f_{j-k}}.
\end{equation*}
Therefore, we obtain $V_ff(k,z)=V_{\tilde{f}}\tilde{f}(k,z)$ for every $(k,z)\in\Omega(g)$, which implies\linebreak $\left|V_gf\right|=\left|V_g\tilde{f}\right|$ by Theorem \ref{thm:VgfvsVff} and Corollary \ref{cor:fourinj}.\\
Conversely, suppose that $\left|V_gf\right|=\left|V_g\tilde{f}\right|$. First, Theorem \ref{thm:VgfvsVff} implies $V_ff(k,z)=V_{\tilde{f}}\tilde{f}(k,z)$ for all $(k,z)\in\Omega(g)$ and therefore $f\cdot \overline{T_k f}=\tilde{f}\cdot \overline{T_k \tilde{f}}$ for every $0\leq k\leq L$ by Corollary \ref{cor:aenonvanZ}.
In particular, $|f|^2=\left|\tilde{f}\right|^2$ and thus $\operatorname{supp}\left(\tilde{f}\right)=\operatorname{supp}(f)$.
Now, let $m\in\N$, fix $j_0\in C_m$ and define $\gamma:=\frac{\tilde{f}_{j_0}}{f_{j_0}}$. Finally, let $j\in C_m$. Then, there exists a finite sequence $\left(j_s\right)_{1\leq s\leq r}$ in $C_m$ such that $j_r=j$ and that $\left|j_{s+1}-j_s\right|\leq L$ holds for every $0\leq s\leq r-1$. Therefore, it holds
\begin{equation*}
\tilde{f}_{j_{s+1}}\overline{\tilde{f}_{j_s}}=f_{j_{s+1}}\overline{f_{j_s}},
\end{equation*}
i.e.
\begin{equation*}
\frac{\tilde{f}_{j_{s+1}}}{f_{j_{s+1}}}=\frac{\overline{f_{j_s}}}{\overline{\tilde{f_{j_s}}}}=\frac{\tilde{f}_{j_s}}{f_{j_s}},
\end{equation*}
once $\frac{\tilde{f}_{j_s}}{f_{j_s}}\in\T$. By induction, we obtain $\tilde{f}_j=\gamma_m f_j$.
\item This follows immediately from a).
\end{enumerate}
\end{proof}

\begin{remark}\label{rem:ZlengthL+1}
Note that the proof of Theorem \ref{thm:ZlengthL+1} shows that it is always possibly to compute $|f|$ from $V_ff\vert_{\Omega(g)}$, and thus from $|V_gf|$, regardless of whether $f$ is indeed phase retrievable.
Since phase retrievability only depends on $\operatorname{supp}(f)$, this implies that one can decide whether a signal is in fact phase retrievable, purely based on the measurement.
\end{remark}

Theorem \ref{thm:ZlengthL+1} strictly generalizes the uniqueness conditions given in \cite[Conditions 1]{LNQ}, since it isn't restricted to one-sided signals.

Next, we will take a look at (right-)one-sided windows, i.e. the case $\operatorname{supp}(g)\subseteq\N_0$. Once more, the index $0$ as start of the support is entirely arbitrary and can be shifted around using Lemma \ref{lm:windowshift}. Here, we obtain
\begin{equation*}
V_gg(k,z)=z^{-k}\sum\limits_{j=0}^{\infty} g_{j+k}\overline{g_j}\overline{z}^{-j}
\end{equation*}
for every $k\geq 0$ and $z\in\T$. Thus, we need to replace the polynomials from above by the power series
\begin{equation*}
P_k(z):=\sum\limits_{j=0}^{\infty} g_{j+k}\overline{g_j}z^j.
\end{equation*}
Obviously, there might exist some $k\geq 0$, for which all coefficients of $P_k(z)$ vanish identically. To manage this problem, we define the set
\begin{equation*}
D_g:=\{l-j~|~j,l\in\operatorname{supp}(g)\}.
\end{equation*}
Obviously, it holds $k\in D_g~\Leftrightarrow~-k\in D_g$.\\
Now, $P_k$ vanishes identically on $K_1(0)$ if and only if $k\notin D_g$ holds.\\
If we could assume that $P_k$ converges on a slightly larger disk $K_r(0)$ for some $r>0$, it would define a holomorphic function on this larger disk, and therefore allow only finitely many zeros on the compact set $\T$. Unfortunately, this assumption doesn't hold true in general. (Consider for instance $g_j:=\frac{1}{j}$.)\\
However, we are still able to show $P_k(z)\neq 0$ for almost every $z\in\T$. In order to do so, we will need some theory on Hardy spaces, provided in \cite{Simon3}.

\begin{lemma}\label{lm:BSpowser}
Let $\left(c_j\right)_{j\in\N_0}\in\ell^1\left(\N_0\right)$ and define
\begin{equation*}
h(z):=\sum\limits_{j=0}^{\infty} c_jz^j
\end{equation*}
for every $z\in\overline{K_1(0)}$. If $c_j\neq 0$ holds for some $j\in\N_0$, it holds $h(z)\neq 0$ for almost every $z\in\T$.
\end{lemma}

\begin{proof}
We will use the notation from \cite{Simon3}.
First, $h$ doesn't vanish identically on $K_1(0)$ by assumption. Additionally, note that $\left(c_j\right)_{j\in\N}\in\ell^1\left(\N_0\right)\subseteq\ell^2\left(\N_0\right)$. \cite[Theorem 5.2.1]{Simon3} therefore implies $h\vert_{K_1(0)}\in H^2$ as well as
\begin{equation*}
\lim\limits_{r\uparrow 1} h(rz)=h(z)\quad\quad\text{for almost every } z\in\T.
\end{equation*}
On the other hand, by Fatou's theorem \cite[Theorem 5.2.6]{Simon3}, there exists a function $h^{\ast}:\T\to\C$ satisfying
\begin{equation*}
\lim\limits_{r\uparrow 1} h(rz)=h^{\ast}(z)\neq 0
\end{equation*}
for almost every $z\in\T$. Together, we obtain $h(z)\neq 0$ for almost every $z\in\T$.
\end{proof}

With Lemma \ref{lm:BSpowser} at hand, it will be easy to generalize Theorem \ref{thm:ZlengthL+1} to (right-)one-sided windows.
However, we will first need to introduce a proper notion of connectivity. Note that this is again easy to check and only depends on the supports of the signal and the window.

\begin{definition} Consider a subset $M\subseteq\Z$ and a signal $f\in\ell^2(\Z)$.
\begin{enumerate}[a)]
\item Two indices $j,k\in M$ are called $g$-connected and we write $j\sim_{\Z,g}k$, iff there exist indices $l_0,\dots,l_n\in M$ satisfying $l_0=j$, $l_n=k$ and $l_{m+1}-l_m\in D_g$ for every $0\leq m<n$.
\item Clearly, $\sim_{\Z,g}$ defines an equivalence relation on $M$. The equivalence classes of $M$ are called $g$-connectivity components of $M$, and $M$ is called $g$-connected, iff there exists at most one connectivity component.
\item $f$ is called $g$-connected, iff $\operatorname{supp}(f)$ is $g$-connected. The $g$-connectivity components of $\operatorname{supp}(f)$ are called $g$-connectivity components of $f$.
\end{enumerate}
\end{definition}

Note that $g$-connectivity is clearly equivalent to $L$-connectivity for windows $g$ satisfying\linebreak$\operatorname{supp}(g)=\{0,\dots,L\}$. With that in my mind, we are able to generalize Theorem \ref{thm:ZlengthL+1}.

\begin{theorem}\label{thm:Zoneside}
Let $g\in\ell^2(\Z)$ be such that $\operatorname{supp}(g)\subseteq\N_0$. Furthermore, consider a signal $f\in\ell^2(\Z)$ and let $\left(C_m\right)_{m\in\N}$ be the $g$-connectivity components of $f$. Then, the following statements hold true.
\begin{enumerate}[a)]
\item A signal $\tilde{f}\in\ell^2(\Z)$ satisfies $\left|V_gf\right|=\left|V_g\tilde{f}\right|$, if and only if $\operatorname{supp}\left(\tilde{f}\right)=\operatorname{supp}(f)$ and for every $m\in\N$, there exists $\gamma_m\in\T$ such that $\tilde{f}\vert_{C_m}=\gamma_mf\vert_{C_m}$.
\item $f$ is phase retrievable w.r.t. $g$ if and only if $f$ is $g$-connected.
\item $g$ does phase retrieval if and only if $D_g=\Z$.
\end{enumerate}
\end{theorem}

\begin{proof}
Clearly, it holds $\Omega(g)\subseteq D_g\times\T$. On the other hand, for every $k\in D_g$, the power series $P_k$ satisfies the assumption of Lemma \ref{lm:BSpowser}. Hence, it follows $P_k(z)\neq 0$ for almost every $z\in\T$.
Consequently, it holds $(k,z)\in\Omega(g)$ for almost every $z\in\T$. The rest of the proof of a) and b) is analogous to the proof of Theorem \ref{thm:ZlengthL+1}.\\
For c), note that every $f\in\ell^2(\Z)$ is $g$-connected, once $D_g=\Z$. Conversely, if there exists $k\in\Z\setminus D_g$, it clearly follows $-k\notin D_g$ and thus, $f:=\delta_0+\delta_k$ is not $g$-connected.
\end{proof}

Note again that the counterexamples used for the proof of part c) as well as analogues thereof are quite common and appear e.g. in \cite{BCEMS,LNQ}. Continuous analogues can be found in \cite{ADGY,GRSpecCl}.\\

Once again, just as in Remark \ref{rem:ZlengthL+1}, it is always possible to determine phase retrievability, just based on the measurement $|V_gf|$.\\

Theorem \ref{thm:Zoneside} provides us with a large class of windows that do phase retrieval. In particular, there exist real-valued windows doing phase retrieval. (Let e.g. $g_j:=\frac{1}{j}$ for $j\in\N$.) This is in contrast to the finite setting, where we will see later on that real-valued windows doing phase retrieval do not exist for even dimensions.\\

Since reflection doesn't affect connectivity, Lemma \ref{lm:windowshift} immediately yields a version for left-one-sided windows, i.e. $\operatorname{supp}(g)\subseteq -\N_0$.

\begin{corollary}\label{cor:Zoneside}
Let $g\in\ell^2(\Z)$ be such that $\operatorname{supp}(g)\subseteq-\N_0$. Furthermore, consider a signal $f\in\ell^2(\Z)$ and let $\left(C_m\right)_{m\in\N}$ be the $g$-connectivity components of $f$. Then, the following statements hold true.
\begin{enumerate}[a)]
\item A signal $\tilde{f}\in\ell^2(\Z)$ satisfies $\left|V_gf\right|=\left|V_g\tilde{f}\right|$, if and only if $\operatorname{supp}\left(\tilde{f}\right)=\operatorname{supp}(f)$ and for every $m\in\N$, there exists $\gamma_m\in\T$ such that $\tilde{f}\vert_{C_m}=\gamma_mf\vert_{C_m}$.
\item $f$ is phase retrievable w.r.t. $g$ if and only if $f$ is $g$-connected.
\item $g$ does phase retrieval if and only if $D_g=\Z$.
\end{enumerate}
\end{corollary}

We give a straight-forward example of some sparse windows doing phase retrieval. The construction uses a so-called \textit{difference set} for the window's support, which isn't a new idea in the field, but more commonly used in the finite setting (see e.g. \cite[Proposition 2.2]{BF}).

\begin{example}\label{ex:sparsewindowZ}
Define sequences $\left(a_n\right)_{n\in\N}$ and $\left(b_n\right)_{n\in\N}$ recursively by $a_0:=b_0:=0$ and
\begin{equation*}
b_{n+1}:=\min\left(\N_0\setminus\{a_m-a_l~|~m,l\leq2n\}\right)
\end{equation*}
as well as
\begin{equation*}
a_{2n+1}:=2a_{2n}+2b_{n+1}\quad\quad\text{and}\quad\quad a_{2n+2}:=2a_{2n}+3b_{n+1}.
\end{equation*}
Furthermore, let $\left(c_n\right)_{n\in\N_0}\in\ell^2(\N_0)$ be such that $c_n\neq 0$ holds for every $n\in\N_0$, and let\linebreak $g:=\sum_{n\in\N_0}c_n\cdot\delta_{a_n}\in\ell^2(\Z)$.
Then, the following statements hold true.
\begin{enumerate}[a)]
\item For every $k\in\N$ there exists exactly one pair $(m,l)\in\N_0^2$ satisfying $a_m-a_l=k$.
\item $g$ does phase retrieval.
\item For every $n^{\ast}\in\N_0$, the window $\tilde{g}:=g-c_{n^{\ast}}\cdot\delta_{a_{n^{\ast}}}$ doesn't do phase retrieval.
\end{enumerate}
\end{example}

\begin{proof}\mbox{}
\begin{enumerate}[a)]
\item First, we will show that $b_n<b_{n+1}$ holds for every $n\in\N_0$. This is obviously true for $n=0$. When $n\in\N$, the definition of $b_n$ implies that for every $k<b_n$ there exists a pair $(m,l)\in\N_0^2$ satisfying $m,l\leq 2n-2$ and $a_m-a_l=k$. Since $a_{2n}-a_{2n-1}=b_n$, it follows $b_{n+1}>b_n$.
Moreover, it is clear that $a_n<a_{n+1}$ holds for every $n\in\N_0$ as well.\\
Now, fix $k\in\N$. Since $(b_n)_{n\in\N_0}$ is strictly increasing, there exists $n\in\N$ such that $k<b_{n+1}$. As before, this implies the existence of $m,l\leq 2n$ satisfying $a_m-a_l=k$.\\
To show uniqueness, choose $(m,l)$ in such a way that $m$ (and therefore also $l$) is minimal. In particular, it holds $k\leq a_m<a_{m+1}$. Suppose that there is a pair $(m^{\prime},l^{\prime})\in\N_0^2$ such that $a_{m^{\prime}}-a_{l^{\prime}}=k$, but $m^{\prime}>m$.\\
If $m^{\prime}$ is odd, i.e. $m^{\prime}=2n+1$ for some $n\in\N_0$, it follows $2n\geq m$. This implies
\begin{equation*}
a_{m^{\prime}}-a_{l^{\prime}}\geq a_{2n+1}-a_{2n}=a_{2n}+2b_{n+1}\geq a_m+2>k,
\end{equation*}
which is a contradiction. Thus, $m^{\prime}$ is even, say $m^{\prime}=2n+2$. When $l^{\prime}=2n+1$, it holds $k=a_{2n+2}-a_{2n+1}=b_{n+1}$. By definition of $b_{n+1}$, the fact that $a_m-a_l=k$ implies $m\geq 2n+1$ and therefore $m=2n+1$ (since $m<m^{\prime}$). We obtain
\begin{equation*}
a_m-a_l\geq a_{2n+1}-a_{2n}=a_{2n}+2b_{n+1}>b_{n+1}=k,
\end{equation*}
which is a contradiction. However, supposing $l^{\prime}\leq 2n$ yields
\begin{equation*}
k=a_{m^{\prime}}-a_{l^{\prime}}\geq a_{2n+2}-a_{2n}=a_{2n}+3b_{n+1}>a_{2n}\geq k,
\end{equation*}
since the case $m>2n$ can be ruled out as before.\\
Altogether, we obtain uniqueness of $(m,l)$.
\item By a), it follows $D_g=\Z$ and Theorem \ref{thm:Vggnonvan} implies that $g$ does phase retrieval.
\item Because of the uniqueness statement in a), there exists $k\in\N$ such that $a_m-a_l\neq k$ for all $m,l\in\N_0\setminus\{n^{\ast}\}$. This yields $k\notin D_g$ and therefore, $g$ can't do phase retrieval by Theorem \ref{thm:Vggnonvan}.
\end{enumerate}
\end{proof}

Finally, we consider the case of a general window $g\in\ell^2(\Z)$. Here, we will have to replace the power series $P_k$ by Laurent series
\begin{equation*}
L_k(z):=\sum\limits_{j\in\Z} g_j\overline{g_{j-k}}z^j.
\end{equation*}
In this case, the theory of Hardy spaces is no longer applicable. However, we may use the initial idea to obtain a theorem that works under the assumption of slightly better convergence.

\begin{theorem}\label{thm:Zgeneral}
Let $g\in\ell^2(\Z)$ and assume that the Laurent series
\begin{equation*}
L_0(z)=\sum\limits_{j\in\Z} \left|g_j\right|^2 z^j
\end{equation*}
converges on the annulus $K_{\nicefrac{1}{r},r}(0)$ for some $r>1$. Furthermore, consider a signal $f\in\ell^2(\Z)$ and let $\left(C_m\right)_{m\in\N}$ be the $g$-connectivity components of $f$. Then, the following statements hold true.
\begin{enumerate}[a)]
\item A signal $\tilde{f}\in\ell^2(\Z)$ satisfies $\left|V_gf\right|=\left|V_g\tilde{f}\right|$, if and only if $\operatorname{supp}\left(\tilde{f}\right)=\operatorname{supp}(f)$ and for every $m\in\N$, there exists $\gamma_m\in\T$ such that $\tilde{f}\vert_{C_m}=\gamma_mf\vert_{C_m}$.
\item $f$ is phase retrievable w.r.t. $g$ if and only if $f$ is $g$-connected.
\item $g$ does phase retrieval if and only if $D_g=\Z$.
\end{enumerate}
\end{theorem}

\begin{proof}
Let $z\in K_{\nicefrac{1}{r},r}(0)$ and write $z=z_1z_2$, where $z_1:=|z|^{\nicefrac{1}{2}}$ and $z_2:=\frac{z}{z_1}$. By Hölder's inequality, it follows that the Laurent series $L_k(z)$ converges for every $k\in\Z$.
Thus, for every $k\in\Z$, $L_k$ defines a holomorphic function on $K_{\nicefrac{1}{r},r}(0)$ which vanishes identically if and only if $k\notin D_g$. Hence, when $k\in D_g$, it may have only finitely many zeros on the compact set $\T$. The rest follows analogously as in the proof of Theorem \ref{thm:Zoneside}.
\end{proof}

\begin{remark}\label{rm:Zgeneral}
Clearly, the assumption of Theorem \ref{thm:Zgeneral} is equivalent to the convergence of
\begin{equation*}
\sum\limits_{j\in\Z} \left|g_j\right|^2e^{\sigma\cdot |j|}
\end{equation*}
for some $\sigma>0$. Therefore, we have established a complete characterization for phase retrieval with exponentially decaying windows.
\end{remark}

\subsection{Necessary conditions for global phase retrieval}

To conclude this section, we briefly want to discuss whether the condition
\begin{equation*}
V_gg(k,z)\neq 0\quad\quad\text{for almost every } (k,z)\in\Z\times\T
\end{equation*}
is also necessary for $g$ to do phase retrieval. Note that we can say at least that $D_g=\Z$ is a necessary condition for $g$ to do phase retrieval, since $f:=\delta_0+\delta_k$ is clearly not phase retrievable w.r.t. $g$, when $k\notin D_g$. (It holds $V_ff(k^{\prime},z)=0$ whenever $k^{\prime}\notin\{k,-k\}$ and the same applies to $\tilde{f}:=\delta_0-\delta_k$.)
Combining this with the fact, that for one-sided windows, the proof of Theorem \ref{thm:Zoneside} shows that it holds $V_gg(k,z)\neq 0$ for almost every $z\in\T$ whenever $k\in D_g$, yields the following statement.

\begin{corollary}\label{cor:prZnec}
Let $g\in\ell^2(\Z)$ be one-sided (i.e. $\operatorname{supp}(g)\subseteq\left(\N_0\right)$ or $\operatorname{supp}(g)\subseteq\left(-\N_0\right)$). Then $g$ does phase retrieval if and only if $V_gg(k,z)\neq 0$ holds for almost every $(k,z)\in\Z\times\T$.
\end{corollary}

For general windows, this statement doesn't remain true anymore. In order to see this, we first state the following reconstruction result.

\begin{theorem}\label{thm:recspecZ}
Let $f,\tilde{f}\in\ell^2(\Z)$ and suppose that there exists an index $k^{\ast}>0$ such that the following conditions hold true.
\begin{enumerate}[(i)]
\item It holds $V_ff(0,z)=V_{\tilde{f}}\tilde{f}(0,z)$ for almost every $z\in\T$.
\item For every $k>k^{\ast}$, it holds $V_ff(k,z)=V_{\tilde{f}}\tilde{f}(k,z)$ for almost every $z\in\T$.
\item It holds $|Z(k)|\geq k^{\ast}+1$ for every $0<k\leq k^{\ast}$, where
\begin{equation*}
Z(k):=\left\{z\in\T~\middle|~V_ff(k,z)=V_{\tilde{f}}\tilde{f}(k,z)\right\}.
\end{equation*}
\end{enumerate}
Then, there exists $\gamma\in\T$ satisfying $\tilde{f}=\gamma f$.
\end{theorem}

\begin{proof}
From (i), we obtain $|f|=\left|\tilde{f}\right|$ by Corollary \ref{cor:aenonvanZ}. In the same way, we obtain\linebreak $f\cdot \overline{T_k f}=\tilde{f}\cdot \overline{T_k\tilde{f}}$ for every $k>k^{\ast}$.\\
First, assume that there are indices $j_0,j_0+N\in\operatorname{supp}(f)=\operatorname{supp}\left(\tilde{f}\right)$, where $N\geq 2k^{\ast}+1$. In this case, fix $\gamma:=\frac{\tilde{f}_{j_0}}{f_{j_0}}$. This implies $\tilde{f}_j=\gamma f_j$ for every $j\geq j_0+k^{\ast}+1$ and in particular $\tilde{f}_{j_0+N}=\gamma f_{j_0+N}$.
This can now be used to obtain $\tilde{f}_j=\gamma f_j$ for every $j\leq j_0+N-k^{\ast}-1$. Since $j_0+N-k^{\ast}-1\geq j_0+k^{\ast}$ holds by assumption, it follows $\tilde{f}=\gamma f$.\\
For the remainder of the proof, we may therefore assume that
\begin{equation*}
\operatorname{supp}(f)=\operatorname{supp}\left(\tilde{f}\right)\subseteq\{j_0,\dots,j_0+N\}
\end{equation*}
holds for some $j_0\in\Z$ and $N\in\N$, where $N\leq 2k^{\ast}$ and $f_{j_0}\neq 0\neq f_{j_0+N}$. As above, we let $\gamma:=\frac{\tilde{f}_{j_0}}{f_{j_0}}$.\\
Notice that $|Z(k)|\geq k^{\ast}+1$ holds for \textit{every} $k\in\N$ by assumptions (ii) and (iii). Therefore, we obtain $\tilde{f}_{j_0+N}=\gamma f_{j_0+N}$ because of
\begin{equation*}
\tilde{f}_{j_0+N}\overline{\tilde{f}_{j_0}}z^{-(j_0+N)}=V_{\tilde{f}}\tilde{f}(N,z)=V_ff(N,z)=f_{j_0+N}\overline{f_{j_0}}z^{-(j_0+N)}
\end{equation*}
for every $z\in Z(N)\neq\emptyset$. We will now prove the statement by induction. Assume that there is some $l\geq 0$ satisfying $\tilde{f}_j=\gamma f_j$ for every $j\in\{j_0,\dots,j_0+l\}\cup\{j_0+N-l,\dots,j_0+N\}$ (which holds true for $l=0$ by what has just been shown).\\
Since this is already equivalent to $\tilde{f}=\gamma f$ whenever $l\geq\frac{N}{2}$, we will assume $l<\frac{N}{2}\leq k^{\ast}$. Our goal is to show that this implies $\tilde{f}_{j_0+l+1}=\gamma f_{j_0+l+1}$ as well as $\tilde{f}_{j_0+N-l-1}=\gamma f_{j_0+N-l-1}$.\\
For every $z\in\T$, it holds
\begin{equation*}
V_ff(N-(l+1),z)=\sum\limits_{j\in\Z} f_j\overline{f_{j-N+l+1}} z^{-j}=\sum\limits_{j=j_0+N-l-1}^{j_0+N} f_j\overline{f_{j-N+l+1}} z^{-j}
\end{equation*}
and the same for $\tilde{f}$. When $j_0+N-l\leq j\leq j_0+N-1$, the induction hypothesis implies both $\tilde{f}_j=\gamma f_j$ and $\tilde{f}_{j-N+l+1}=\gamma f_{j-N+l+1}$, which yields
\begin{equation*}
\sum\limits_{j=j_0+N-l}^{j_0+N-1}\tilde{f}_j\overline{\tilde{f}_{j-N+l+1}}z^{-j}=\sum\limits_{j=j_0+N-l}^{j_0+N-1}f_j\overline{f_{j-N+l+1}}z^{-j}.
\end{equation*}
Since $V_ff(N-(l+1),z)=V_{\tilde{f}}\tilde{f}(N-(l+1),z)$ holds for every $z\in Z(N-(l+1))$, it follows $W(z,l)=\widetilde{W}(z,l)$, where we write
\begin{align*}
W(z,l)&:=f_{j_0+N}\overline{f_{j_0+l+1}}z^{-(j_0+N)}+f_{j_0+N-l-1}\overline{f_{j_0}}z^{-(j_0+N-l-1)}\quad\text{as well as}\\\widetilde{W}(z,l)&:=\tilde{f}_{j_0+N}\overline{\tilde{f}_{j_0+l+1}}z^{-(j_0+N)}+\tilde{f}_{j_0+N-l-1}\overline{\tilde{f}_{j_0}}z^{-(j_0+N-l-1)}.
\end{align*}
Combining $l+1\leq k^{\ast}$ with the fact that $|Z(N-(l+1))|\geq k^{\ast}+1$, yields the existence of $z_0,z_1\in Z(N-(l+1))$ satisfying $z_0^{l+1}\neq z_1^{l+1}$. This allows us to compute
\begin{equation*}
z_0^{j_0+N}\cdot W(z_0,l)-z_1^{j_0+N}\cdot W(z_1,l)=f_{j_0+N-l-1}\overline{f_{j_0}}\cdot\underbrace{\left(z_0^{l+1}-z_1^{l+1}\right)}_{\neq 0}.
\end{equation*}
Since the same calculation applies to $\tilde{f}$, we finally obtain
\begin{equation*}
\tilde{f}_{j_0+N-l-1}\overline{\tilde{f}_{j_0}}=f_{j_0+N-l-1}\overline{f_{j_0}},
\end{equation*}
which entails $\tilde{f}_{j_0+N-l-1}=\gamma f_{j_0+N-l-1}$ because of $\tilde{f}_{j_0}=\gamma f_{j_0}\neq 0$.\\
Plugging this into $W(z_0,l)=\widetilde{W}(z_0,l)$ yields
\begin{equation*}
\tilde{f}_{j_0+N}\overline{\tilde{f}_{j_0+l+1}}=f_{j_0+N}\overline{f_{j_0+l+1}}
\end{equation*}
and therefore $\tilde{f}_{j_0+l+1}=\gamma f_{j_0+l+1}$ because of $\tilde{f}_{j_0+N}=\gamma f_{j_0+N}\neq 0$.\\
By induction, we obtain $\tilde{f}=\gamma f$.
\end{proof}

We may reformulate (and slightly generalize) Theorem \ref{thm:recspecZ} in terms of phase retrieval.

\begin{corollary}\label{cor:prspecZ}
Let $g\in\ell^2(\Z)$ be such that $D_g=\Z$ as well as $V_gg(0,z)\neq 0$ for almost every $z\in\T$. Furthermore, assume that there exists $k^{\ast}\in\N$ such that for every $k>k^{\ast}$, it holds $V_gg(k,z)\neq 0$ for almost every $z\in\T$. Then $g$ does phase retrieval.
\end{corollary}

\begin{proof}
Let $f,\tilde{f}\in\ell^2(\Z)$ be such that $V_ff(k,z)=V_{\tilde{f}}\tilde{f}(k,z)$ holds for every $(k,z)\in\Omega(g)$. In particular, whenever $k>k^{\ast}$ or $k=0$, it holds $V_ff(k,z)=V_{\tilde{f}}\tilde{f}(k,z)$ for almost every $z\in\T$.\\
Since $D_g=\Z$, it holds $g\cdot \overline{T_k g}\not\equiv 0$ for every $k\in\Z$. Together with the injectivity of the Fourier transform, there exists at least one (and by continuity of $V_gg$ therefore infinitely many) $z\in\T$ such that $V_gg(k,z)\neq 0$. Consequently, there are infinitely many $z\in\T$ satisfying $V_ff(k,z)=V_{\tilde{f}}\tilde{f}(k,z)$.\\
According to Theorem \ref{thm:recspecZ}, it follows $\tilde{f}=\gamma f$, and thus $g$ does phase retrieval by Lemma \ref{lm:Vff}.
\end{proof}

In the following example, we present a window that satisfies the assumptions of Corollary \ref{cor:prspecZ}, but not those of Theorem \ref{thm:Zgeneral}. In particular, we will see that Corollary \ref{cor:prZnec} doesn't hold true for general windows.

\begin{example}
Let $g\in\ell^2(\Z)$ be given by
\begin{equation*}
g_j:=\begin{cases} \frac{(-1)^k}{\sqrt{\pi}}\cdot\frac{1}{j+1},&\text{ if } j=2k \text{ for some } k\in\Z,\\ \frac{i}{4},&\text{if } j\in\{1,-3\},\\1,&\text{if } j=-1,\\ 0,&\text{otherwise}.\end{cases}
\end{equation*}
Then the following hold true.
\begin{enumerate}[a)]
\item It holds $V_gg(2,z)=0$ for every $z\in\T$ satisfying $\operatorname{Im}(z)<0$.
\item $g$ does phase retrieval.
\end{enumerate}
\end{example}

\begin{proof}
As a preliminary result, we will first prove the identity
\begin{align}
\begin{split}\label{serLog}
&~\sum\limits_{m\in\Z} \frac{\overline{z}^{2m}}{(2m+1)\cdot(2m-2l+1)}\\
=&~\frac{\overline{z}^{2l-1}-z}{4l}\cdot\left(\operatorname{Log}(1+z^{-1})-\operatorname{Log}(1-z^{-1})-\operatorname{Log}(1+z)+\operatorname{Log}(1-z)\right)\\
=&~\frac{z-\overline{z}^{2l-1}}{2l}\cdot i\cdot\operatorname{Im}\left(\operatorname{Log}(1+z)-\operatorname{Log}(1-z)\right)
\end{split}
\end{align}
for almost every $z\in\T$, where $l\in\N$ is fixed.\\
In order to show this, recall that the series
\begin{equation*}
\sum\limits_{m\in\Z}\frac{\overline{z}^{2m}}{2m+1}\quad\quad\text{and}\quad\quad\sum\limits_{m\in\Z}\frac{\overline{z}^{2m}}{2m-2l+1}
\end{equation*}
converge for every $z\in\T,~z\notin\{1,-1\}$ according to Dirichlet's test. We may therefore use the partial fractional decomposition
\begin{equation*}
\frac{1}{(2m+1)\cdot(2m-2l+1)}=\frac{1}{2l}\cdot\frac{1}{2m-2l+1}-\frac{1}{2l}\cdot\frac{1}{2m+1}
\end{equation*}
to rewrite
\begin{align*}
\sum\limits_{m\in\Z}\frac{\overline{z}^{2m}}{(2m+1)\cdot(2m-2l+1)}	&=\frac{1}{2l}\cdot\sum\limits_{m\in\Z}\frac{\overline{z}^{2m}}{2m-2l+1}-\frac{1}{2l}\cdot\sum\limits_{m\in\Z}\frac{\overline{z}^{2m}}{2m+1}\\
											&=\frac{\overline{z}^{2l-1}-z}{2l}\cdot\sum\limits_{m\in\Z}\frac{\overline{z}^{2m+1}}{2m+1}\\
											&=\frac{\overline{z}^{2l-1}-z}{2l}\cdot\sum\limits_{m=0}^{\infty} \frac{\overline{z}^{2m+1}-z^{2m+1}}{2m+1}.
\end{align*}
for every $z\in\T,~z\notin\{1,-1\}$. When $|z|<1$, the fact that
\begin{equation*}
\operatorname{Log}(1+z)=\sum\limits_{n=1}^{\infty}\frac{(-1)^{n-1}}{n}\cdot z^n\quad\quad\text{and}\quad\quad\operatorname{Log}(1-z)=\sum\limits_{n=1}^{\infty} -\frac{1}{n}z^n
\end{equation*}
hold true, implies
\begin{equation*}
\sum\limits_{m=0}^{\infty}\frac{1}{2m+1}z^{2m+1}=\frac{1}{2}\cdot\left(\operatorname{Log}(1+z)-\operatorname{Log}(1-z)\right).
\end{equation*}
Since the coefficients of the power series on the left-hand side form a sequence in $\ell^2(\N)$, this function is an element of the Hardy space $H^2$. Using again the results in \cite{Simon3}, it follows
\begin{equation*}
\sum\limits_{m=0}^{\infty}\frac{z^{2m+1}}{2m+1}=\lim\limits_{r\uparrow 1} \frac{1}{2}\cdot\left(\operatorname{Log}(1+rz)-\operatorname{Log}(1-rz)\right)=\frac{1}{2}\cdot\left(\operatorname{Log}(1+z)-\operatorname{Log}(1-z)\right)
\end{equation*}
for almost every $z\in\T$, where we used the fact that the functions $s\mapsto \operatorname{Log}(1+s)$ and\linebreak $s\mapsto \operatorname{Log}(1-s)$ are continuous in every $z\in\T,~z\notin\{1,-1\}$.\\
Altogether, we obtain (\ref{serLog}).
\begin{enumerate}[a)]
\item For every $j\in\Z$, it holds
\begin{equation*}
g_j\overline{g_{j-2}}=\begin{cases} \frac{i}{4},&\text{if } j=1,\\ -\frac{i}{4},&\text{if } j=-1,\\ -\frac{1}{\pi}\cdot\frac{1}{(j+1)\cdot(j-1)},&\text{if } j\text{ is even},\\0,&\text{otherwise}.\end{cases}
\end{equation*}
Together with (\ref{serLog}), we obtain
\begin{align*}
V_gg(2,z)=\sum\limits_{j\in\Z}g_j\overline{g_{j-2}}\overline{z}^j	&=\frac{i}{4}\overline{z}-\frac{i}{4}z-\frac{1}{\pi}\sum\limits_{m\in\Z}\frac{\overline{z}^{2m}}{(2m+1)\cdot(2m-1)}\\
										&=\frac{1}{2}\operatorname{Im}(z)-\frac{1}{\pi}\cdot i\cdot\operatorname{Im}(z)\cdot i\cdot\operatorname{Im}\left(\operatorname{Log}(1+z)-\operatorname{Log}(1-z)\right)\\
										&=\frac{1}{2}\operatorname{Im}(z)+\frac{1}{\pi}\cdot\operatorname{Im}(z)\cdot\left(\operatorname{arg}(1+z)-\operatorname{arg}(1-z)\right)
\end{align*}
for almost every $z\in\T$. Plugging in $z=e^{i\varphi}$ for some $\varphi\in(0,\pi)\cup(\pi,2\pi)$, yields
\begin{equation*}
\operatorname{arg}(1+z)-\operatorname{arg}(1-z)=\begin{cases}\arccos\left(\frac{1+\cos(\varphi)}{|1+z|}\right)+\arccos\left(\frac{1-\cos(\varphi)}{|1-z|}\right),&\text{if } \varphi\in(0,\pi),\\ -\arccos\left(\frac{1+\cos(\varphi)}{|1+z|}\right)-\arccos\left(\frac{1-\cos(\varphi)}{|1-z|}\right),&\text{if } \varphi\in(\pi,2\pi),\end{cases}
\end{equation*}
where we can compute
\begin{equation*}
\arccos\left(\frac{1+\cos(\varphi)}{|1+z|}\right)=\arccos\left(\sqrt{\frac{1+\cos(\varphi)}{2}}\right)=\arccos\left(\sqrt{\cos^2\left(\frac{\varphi}{2}\right)}\right)
\end{equation*}
as well as
\begin{align*}
\arccos\left(\frac{1-\cos(\varphi)}{|1-z|}\right)	&=\arccos\left(\sqrt{\frac{1-\cos(\varphi)}{2}}\right)=\arccos\left(\sqrt{\sin^2\left(\frac{\varphi}{2}\right)}\right)\\
								&=\arccos\left(\sin\left(\frac{\varphi}{2}\right)\right)=\frac{\pi}{2}-\arcsin\left(\sin\left(\frac{\varphi}{2}\right)\right),
\end{align*}
since $\frac{\varphi}{2}\in(0,\pi)$.\\
For $\varphi\in(0,\pi)$, we obtain
\begin{equation*}
\arccos\left(\sqrt{\cos^2\left(\frac{\varphi}{2}\right)}\right)=\arccos\left(\cos\left(\frac{\varphi}{2}\right)\right)=\frac{\varphi}{2}
\end{equation*}
and
\begin{equation*}
\arcsin\left(\sin\left(\frac{\varphi}{2}\right)\right)=\frac{\varphi}{2}.
\end{equation*}
For $\varphi\in(\pi,2\pi)$, it holds
\begin{equation*}
\arccos\left(\sqrt{\cos^2\left(\frac{\varphi}{2}\right)}\right)=\arccos\left(-\cos\left(\frac{\varphi}{2}\right)\right)=\pi-\arccos\left(\cos\left(\frac{\varphi}{2}\right)\right)=\pi-\frac{\varphi}{2}
\end{equation*}
and
\begin{equation*}
\arcsin\left(\sin\left(\frac{\varphi}{2}\right)\right)=\pi-\frac{\varphi}{2}.
\end{equation*}
Thus, in both cases
\begin{equation*}
\arccos\left(\frac{1+\cos(\varphi)}{|1+z|}\right)+\arccos\left(\frac{1-\cos(\varphi)}{|1-z|}\right)=\frac{\pi}{2}.
\end{equation*}
Altogether, it follows
\begin{align*}
V_gg(2,z)	&=\begin{cases} \operatorname{Im}(z),&\text{if }\operatorname{Im}(z)>0,\\ 0,&\text{if }\operatorname{Im}(z)<0\end{cases}\\
		&=\sin(\varphi)\cdot\chi_{[0,\pi]}(\varphi)\quad\quad\quad\quad \left(z=e^{i\varphi}\right)
\end{align*}
for almost every $z\in\T$. By continuity of $V_gg$, the identity holds for every $z\in\T$.
\item In order to apply Corollary \ref{cor:prspecZ}, we will first show that
\begin{equation*}
0\neq V_gg(0,z)=\sum\limits_{j\in\Z} |g_j|^2\overline{z}^j=\frac{1}{16}\overline{z}+\frac{1}{16}z^3+z+\frac{1}{\pi}\cdot\sum\limits_{m\in\Z}\frac{\overline{z}^{2m}}{(2m+1)^2}
\end{equation*}
holds for every $z\in\T$. This follows from the estimate
\begin{align*}
\left|V_gg(0,z)\right|=\left|V_gg(0,z)\cdot\overline{z}\right|	&=\left|\frac{1}{16}\cdot\left(z^2+\overline{z}^2\right)+1+\frac{1}{\pi}\cdot\sum\limits_{m\in\Z}\frac{\overline{z}^{2m+1}}{(2m+1)^2}\right|\\
										&\geq1-\frac{1}{8}\cdot\left|\operatorname{Re}\left(z^2\right)\right|-\frac{1}{\pi}\cdot\left|\sum\limits_{m\in\Z}\frac{\overline{z}^{2m+1}}{(2m+1)^2}\right|\\
										&\geq\frac{7}{8}-\frac{1}{\pi}\cdot\sum\limits_{m\in\Z}\frac{1}{(2m+1)^2}=\frac{7}{8}-\frac{\pi}{4}>0.
\end{align*}
Since $D_g=\Z$ clearly holds true, we obtain that $g$ does phase retrieval as soon as we can show that every $k>4$ satisfies $V_gg(k,z)\neq 0$ for almost every $z\in\T$.\\
When $k$ is odd, this follows from the fact that $V_gg(k,\cdot)$ is a polynomial. When $k=2l\geq 6$ is even, (\ref{serLog}) implies
\begin{align*}
V_gg(k,z)	&=\sum\limits_{j\in\Z} g_j\overline{g_{j-2l}}\overline{z}^j=\sum\limits_{m\in\Z} g_{2m}\overline{g_{2m-2l}}\overline{z}^{2m}\\
		&=\frac{(-1)^l}{\pi}\cdot\sum\limits_{m\in\Z}\frac{1}{(2m+1)\cdot(2m-2l+1)}\overline{z}^{2m}\\
		&=(-1)^l\cdot\frac{z-\overline{z}^{2l-1}}{4l\pi}\cdot\underbrace{\left(\operatorname{Log}(1+z^{-1})-\operatorname{Log}(1-z^{-1})-\operatorname{Log}(1+z)+\operatorname{Log}(1-z)\right)}_{=:f(z)}
\end{align*}
for almost every $z\in\T$.\\
Now, suppose that $V_gg(k,z)=0$ holds for uncountably many $z\in\T$. In particular, this would imply that for some $n\in\N$, the set
\begin{equation*}
\left\{z\in\T~\middle|~\left|\operatorname{Im}(z)\right|\geq\frac{1}{n}, f(z)=0\right\}
\end{equation*}
is uncountable, and thus has a limit point within the compact set $\left\{z\in\T~\middle|~\left|\operatorname{Im}(z)\right|\geq\frac{1}{n}\right\}$. However, this implies that the holomorphic function
\begin{equation*}
f_n:\left\{z\in\C~\middle|~\left|\operatorname{Im}(z)\right|>\frac{1}{n+1}\right\}\to\C,\quad f_n(z):=f(z)
\end{equation*}
vanishes identically, which is a contradiction.\\
Hence, $V_gg(k,\cdot)$ may only have finitely many zeros and therefore, $g$ does phase retrieval.
\end{enumerate}
\end{proof}

On the other hand, it remains an open question whether there exists a window $g$ not doing phase retrieval despite satisfying $D_g=\Z$.

\section{STFT phase retrieval on $\Z_d$}\label{sec:Zd}

We proceed by considering the cyclic group $G=\Z_d=\Z/d\Z$ for $d\geq 2$. Throughout this section, all indices as well as their sums and products are to be understood modulo $d$, unless mentioned otherwise.\\
Since $\widehat{G}\cong G$, where the characters are given by
\begin{equation*}
\chi_l:G\to\T,\quad j\mapsto e^{\nicefrac{2\pi ijl}{d}},
\end{equation*}
we obtain
\begin{equation*}
\widehat{h}(l)=\sum\limits_{j=0}^{d-1} h_j e^{-\nicefrac{2\pi ijl}{d}}\quad (0\leq l\leq d-1)
\end{equation*}
for the Fourier transform and
\begin{equation*}
V_gf(k,l)=\sum\limits_{j=0}^{d-1} f_j\overline{g_{j-k}}e^{-\nicefrac{2\pi ijl}{d}}\quad(0\leq k,l\leq d-1)
\end{equation*}
for the STFT.

Since obviously $\overline{\chi_l}=\chi_{-l}$, the symmetry relation from Lemma \ref{lm:sym} becomes
\begin{equation}\label{symZd}
V_hh(-k,-l)=e^{-\nicefrac{2\pi ikl}{d}}\cdot\overline{V_hh(k,l)}.
\end{equation}

In contrast to the previous section, we will begin by analyzing necessity of the condition\linebreak $\Omega(g)=\Z_d\times\Z_d$ for phase retrieval, in order to illustrate the wider range of outcomes in this finite setting.
Afterwards, we will consider a special class of short windows and discuss whether there exist analogs to Theorem \ref{thm:ZlengthL+1}.

\subsection{Phase retrieval and non-vanishing ambiguity functions}

Recall that Theorem \ref{thm:Vggnonvan} ensures that every window $g\in\C^d$ satisfying $\Omega(g)=G\times\widehat{G}=\Z_d\times\Z_d$, does phase retrieval.
It is also well-known (see e.g. \cite[Proposition 2.1]{BF}) that this condition is fulfilled with probability $1$ when the window $g$ is picked randomly.

\begin{theorem}\label{thm:nonvangen}
There exists a measurable set $E\subseteq\C^d$ of measure zero such that every\linebreak $g\in\C^d\setminus E$ fulfills the condition $\Omega(g)=\Z_d\times\Z_d$ and therefore does phase retrieval.
\end{theorem}

Since $d^2\geq 4d-4$ always holds true, this appears to be in line with a general property of phase retrieval, presented in \cite[Theorem 1.1]{CEHV}. Note however, that this result doesn't \textit{imply} Theorem \ref{thm:nonvangen}, since there is no evidence on whether a frame corresponding to STFT phase retrieval is generic in the sense of \cite{CEHV}.\\

However, Theorem \ref{thm:nonvangen} states nothing about phase retrieval for a \textit{fixed} window $g\in\C^d$.\\
In this subsection, we want to discuss whether the condition $\Omega(g)=\Z_d\times\Z_d$ is also necessary for $g$ to do phase retrieval. We will see that this is indeed the case when $d\in\{2,3\}$, but -- more interestingly -- not when $d\geq 4$.
Note that we already know by \cite[Proposition 2.3]{BF} that for each $0\leq k\leq d-1$, there has to be at least one $0\leq l\leq d-1$ (and $1\leq l\leq d-1$ when $k=0$) such that $(k,l)\in\Omega(g)$, in order for $g$ to do phase retrieval.\\
Intuitively, one would also like to use general results like \cite[Proposition 3.6]{BCE} or \cite[Theorem 5.2]{CEHV}, in order to establish a lower bound on $|\Omega(g)|$ for phase retrieval. Note however, that results of this type don't apply, since we leave classical phase retrieval territory when considering recovery from $V_ff\vert_{\Omega(g)}$.\\

First, we provide a helpful observation which seems unique to the case $d=3$.

\begin{lemma}\label{lm:ambigshift}
Consider a signal $f\in\C^3$ as well as
\begin{equation*}
f^{\prime}:=\vek{f_0}{e^{\nicefrac{2\pi i}{3}}f_1}{f_2}\quad\text{and}\quad f^{\prime\prime}:=\vek{f_0}{e^{\nicefrac{-2\pi i}{3}}f_1}{f_2}.
\end{equation*}
Then,
\begin{equation*}
V_{f^{\prime}}f^{\prime}(0,l)=V_{f^{\prime\prime}}f^{\prime\prime}(0,l)=V_ff(0,l)
\end{equation*}
as well as
\begin{equation*}
V_{f^{\prime}}f^{\prime}(1,l)=V_ff(1,l-1)\quad\text{and}\quad V_{f^{\prime\prime}}f^{\prime\prime}(1,l)=V_ff(1,l+1)
\end{equation*}
hold for every $l\in\{0,1,2\}$.
\end{lemma}

The proof is a straight forward computation. An application of (\ref{symZd}) additionally yields
\begin{equation*}
V_{f^{\prime}}f^{\prime}(2,l)=e^{-\nicefrac{2\pi i}{3}}V_ff(2,l+1)\quad\text{and}\quad V_{f^{\prime\prime}}f^{\prime\prime}(2,l)=e^{\nicefrac{2\pi i}{3}}V_ff(2,l-1)
\end{equation*}
for every $l\in\{0,1,2\}$.

Because of Lemma \ref{lm:ambigshift}, we basically need just one main signal to prove our conjecture when $d=3$.

\begin{lemma}\label{lm:signal1122}
Consider the signals
\begin{equation*}
f:=\vek{1}{\frac{1}{\sqrt{3}}e^{-\nicefrac{\pi i}{6}}}{\frac{1}{\sqrt{3}}e^{\nicefrac{5\pi i}{6}}}\quad\text{and}\quad\tilde{f}:=\vek{1}{\frac{1}{\sqrt{3}}e^{-\nicefrac{\pi i}{2}}}{\frac{1}{\sqrt{3}}e^{-\nicefrac{5\pi i}{6}}}.
\end{equation*}
It holds $V_ff(k,l)=V_{\tilde{f}}{\tilde{f}}(k,l)$ for every $(k,l)\in\{0,1,2\}^2$ satisfying $(k,l)\notin\{(1,1),(2,2)\}$.
\end{lemma}

This is proven once more by directly computing the entries $V_ff(k,l)$, where the only effort lies within $k=1$: Since $|f|=\left|\tilde{f}\right|$, the statement is obvious for $k=0$. For $k=2$, one may use (\ref{symZd}).

We are now able to prove that, when $d\in\{2,3\}$, phase retrieval is equivalent to a nowhere-vanishing ambiguity function $V_gg$.

\begin{theorem}\label{thm:PRVgg23}
Let $d\in\{2,3\}$. If $g$ does phase retrieval, then $\Omega(g)=\Z_d\times\Z_d$.
\end{theorem}

\begin{proof}
First, assume $d=2$. We will show for every $(k,l)\in\{0,1\}^2$ that $V_gg(k,l)=0$ implies that $g$ can't do phase retrieval.
\begin{enumerate}[1)]
\item If $V_gg(0,0)=0$, then clearly $g=0$ and thus, $g$ can't do phase retrieval.
\item If $V_gg(0,1)=0$, consider
\begin{equation*}
f:=\vek{2}{1}{}\quad\text{and}\quad \tilde{f}:=\vek{1}{2}{}.
\end{equation*}
One directly computes $V_ff(k,l)=V_{\tilde{f}}{\tilde{f}}(k,l)$ whenever $(k,l)\neq(0,1)$ and therefore, $g$ can't do phase retrieval.
\item If $V_gg(1,0)=0$, consider
\begin{equation*}
f:=\vek{2}{1}{}\quad\text{and}\quad \tilde{f}:=\vek{2}{-1}{}.
\end{equation*}
One directly computes $V_ff(k,l)=V_{\tilde{f}}{\tilde{f}}(k,l)$ whenever $(k,l)\neq(1,0)$ and therefore, $g$ can't do phase retrieval.
\item If $V_gg(1,1)=0$, consider
\begin{equation*}
f:=\vek{2}{i}{}\quad\text{and}\quad \tilde{f}:=\vek{2}{-i}{}.
\end{equation*}
One directly computes $V_ff(k,l)=V_{\tilde{f}}{\tilde{f}}(k,l)$ whenever $(k,l)\neq(1,1)$ and therefore, $g$ can't do phase retrieval.
\end{enumerate}

Now, let $d=3$. We will once again analyze the vanishing of every single $V_gg(k,l)$ separately, where we may identify pairs thanks to (\ref{symZd}).
\begin{enumerate}[1)]
\item If $V_gg(0,0)=0$, then clearly $g=0$ and thus, $g$ can't do phase retrieval.
\item Let $V_gg(0,1)=0$. This is equivalent to $V_gg(0,2)=0$. Consider
\begin{equation*}
h:=\vek{1}{0}{0}\quad\text{and}\quad\tilde{h}:=\vek{0}{1}{0}.
\end{equation*}
Then, $V_hh(0,0)=V_{\tilde{h}}\tilde{h}(0,0)=1$ as well as $V_hh(k,l)=V_{\tilde{h}}\tilde{h}(k,l)=0$ for $k\neq 0$, and therefore, $g$ can't do phase retrieval.
\item Let $V_gg(1,1)=0$. This is equivalent to $V_gg(2,2)=0$. Thanks to the signals $f$ and $\tilde{f}$ from Lemma \ref{lm:signal1122}, it immediately follows that $g$ can't do phase retrieval.
\item Let $V_gg(1,2)=0$. This is equivalent to $V_gg(2,1)=0$. With $f$ and $\tilde{f}$ as above, Lemmas \ref{lm:signal1122} and \ref{lm:ambigshift} yield $V_{f^{\prime}}f^{\prime}(k,l)=V_{\tilde{f}^{\prime}}\tilde{f}^{\prime}(k,l)$ for every $(k,l)\in\{0,1,2\}^2\setminus\{(1,2),(2,1)\}$, and therefore, $g$ can't do phase retrieval.
\item Let $V_gg(1,0)=0$. This is equivalent to $V_gg(2,0)=0$. With $f$ and $\tilde{f}$ as above, Lemmas \ref{lm:signal1122} and \ref{lm:ambigshift} yield $V_{f^{\prime\prime}}f^{\prime\prime}(k,l)=V_{\tilde{f}^{\prime\prime}}\tilde{f}^{\prime\prime}(k,l)$ for every $(k,l)\in\{0,1,2\}^2\setminus\{(1,0),(2,0)\}$, and therefore, $g$ can't do phase retrieval.
\end{enumerate}
\end{proof}

For larger dimensions, the analog of Theorem \ref{thm:PRVgg23} doesn't hold true. In order to give an example of a window that does phase retrieval even though having one ore more zeros on its ambiguity function, we will first prove the following reconstruction result.

\begin{theorem}\label{thm:recwod2}
Let $d\geq 4$ be even and let $f,\tilde{f}\in\C^d$ be such that $V_ff(k,l)=V_{\tilde{f}}\tilde{f}(k,l)$ for every $(k,l)\in(\Z_d\times\Z_d)\setminus\left\{\left(\frac{d}{2},\frac{d}{2}\right)\right\}$.
Then, there exists $\gamma\in\T$ such that $\tilde{f}=\gamma f$.
\end{theorem}

\begin{proof}
First, the assumption implies $f\cdot \overline{T_k f}=\tilde{f}\cdot \overline{T_k\tilde{f}}$ for every $0\leq k\leq d-1, k\neq \frac{d}{2}$. Choosing $k=0$, this shows $|f|=\left|\tilde{f}\right|$.
We will distinguish various cases depending on the elements of $\operatorname{supp}(f)=\operatorname{supp}(\tilde{f})$.
\begin{enumerate}[1)]
\item If $|\operatorname{supp}(f)|\leq 1$, then $\tilde{f}=\gamma f$ follows immediately for some $\gamma\in\T$.
\item If $|\operatorname{supp}(f)|\geq 3$, pick $j\in\operatorname{supp}(f)$ and choose $\gamma:=\frac{\tilde{f}_j}{f_j}$. The fact that $f\cdot \overline{T_k f}=\tilde{f}\cdot \overline{T_k \tilde{f}}$ holds for every $0\leq k\leq d-1, k\neq \frac{d}{2}$ then implies $\tilde{f}_m=\gamma f_m$ for every $m\neq j+\frac{d}{2}$.
Since $|\operatorname{supp}(f)|\geq 3$, there exists $m\notin\left\{j,j+\frac{d}{2}\right\}$ satisfying $0\neq\tilde{f}_m=\gamma f_m$. Since $j+\frac{d}{2}-m\neq\frac{d}{2}$, we obtain $f_{j+\frac{d}{2}}\overline{f_m}=\tilde{f}_{j+\frac{d}{2}}\overline{\tilde{f}_m}$ and thus $\tilde{f}=\gamma f$.
\item If $\operatorname{supp}(f)=\{j,m\}$ with $m\notin\{j,j+\frac{d}{2}\}$, it follows $\tilde{f}=\gamma f$ as in 2).
\item The remaining case is given by $\operatorname{supp}(f)=\left\{j,j+\frac{d}{2}\right\}$ for some $0\leq j\leq d-1$. W.l.o.g assume that $j=0$ and let $\gamma:=\frac{\tilde{f}_0}{f_0}$ It then follows
\begin{equation*}
V_ff\left(\frac{d}{2},l\right)=\sum\limits_{j=0}^{d-1}f_j f_{j-\frac{d}{2}}e^{-\nicefrac{2\pi ijl}{d}}=f_0\overline{f_{\frac{d}{2}}}+f_{\frac{d}{2}}\overline{f_0}e^{-\pi il}
\end{equation*}
for every $0\leq l\leq d-1$ and analogously $V_{\tilde{f}}\tilde{f}\left(\frac{d}{2},l\right)=\tilde{f}_0\overline{\tilde{f}_{\frac{d}{2}}}+\tilde{f}_{\frac{d}{2}}\overline{\tilde{f}_0}e^{-\pi il}$.
This implies
\begin{equation*}
2\tilde{f}_0\overline{f_{\frac{d}{2}}}=V_{\tilde{f}}\tilde{f}\left(\frac{d}{2},0\right)+V_{\tilde{f}}\tilde{f}\left(\frac{d}{2},1\right)=V_ff\left(\frac{d}{2},0\right)+V_ff\left(\frac{d}{2},1\right)=2f_0\overline{f_{\frac{d}{2}}}
\end{equation*}
and thus $\tilde{f}=\gamma f$.
\end{enumerate}
\end{proof}

It remains to show that there is in fact a window $g\in\C^d$ satisfying $\Omega(g)=(\Z_d\times\Z_d)\setminus\left\{\left(\frac{d}{2},\frac{d}{2}\right)\right\}$.

\begin{example}\label{ex:PRVgg}
Let $d\geq 4$ be even. If $d$ is a multiple of $4$, consider $g\in\C^d$ given by
\begin{equation*}
g_j:=\begin{cases} 2^j,&\text{if } 0\leq j\leq\frac{d}{2}-1,\\
			\frac{1}{2^{j-\nicefrac{d}{2}}}\left(1+\sqrt{2^{4j-d}-1}\cdot i\right),&\text{if } j\in\left\{\frac{d}{2},\frac{d}{2}+1\right\},\\
			\frac{1}{2^{j-\nicefrac{d}{2}}}\left(2+\sqrt{2^{4j-d}-4}\cdot i\right),&\text{if } \frac{d}{2}+2\leq j\leq d-1.
\end{cases}
\end{equation*}
Otherwise, let
\begin{equation*}
g_j:=\begin{cases} 2^j,&\text{if } 0\leq j\leq\frac{d}{2}-1,\\
			\frac{1}{2^{j-\nicefrac{d}{2}}}\left(\sqrt{2^{4j-d}-1}+i\right),&\text{if } j\in\left\{\frac{d}{2},\frac{d}{2}+1\right\},\\
			\frac{1}{2^{j-\nicefrac{d}{2}}}\left(\sqrt{2^{4j-d}-4}+2i\right),&\text{if } \frac{d}{2}+2\leq j\leq d-1.
\end{cases}
\end{equation*}
Then, $V_gg\left(\frac{d}{2},\frac{d}{2}\right)=0$, but $g$ does phase retrieval.
\end{example}

\begin{proof}
We will only consider the case that $d$ is a multiple of $4$ (and thus equivalently $\frac{d}{2}$ is even). The other case is treated entirely analogously.

First, we compute
\begin{align*}
V_gg\left(\frac{d}{2},l\right)	&=\sum\limits_{j=0}^{\frac{d}{2}-1} g_j\overline{g_{j+\frac{d}{2}}}e^{-\nicefrac{2\pi ijl}{d}}~+~\sum\limits_{j=\frac{d}{2}}^{d-1}g_j\overline{g_{j-\frac{d}{2}}}e^{-\nicefrac{2\pi ijl}{d}}\\
					&=\sum\limits_{j=0}^{\frac{d}{2}-1} g_j\overline{g_{j+\frac{d}{2}}}e^{-\nicefrac{2\pi ijl}{d}}~+~\sum\limits_{j=0}^{\frac{d}{2}-1}g_{j+\frac{d}{2}}\overline{g_j}e^{-\nicefrac{2\pi ijl}{d}}e^{-\pi il}\\
					&=\sum\limits_{j=0}^{\frac{d}{2}-1} g_j\overline{g_{j+\frac{d}{2}}}e^{-\nicefrac{2\pi ijl}{d}}+(-1)^l\overline{g_j\overline{g_{j+\frac{d}{2}}}}e^{-\nicefrac{2\pi ijl}{d}}\\
					&=\begin{cases} 2\sum\limits_{j=0}^{\frac{d}{2}-1} e^{-\nicefrac{2\pi ijl}{d}}\operatorname{Re}(g_j\overline{g_{j+\frac{d}{2}}}),&\text{if } l \text{ is even},\\
								2i\sum\limits_{j=0}^{\frac{d}{2}-1} e^{-\nicefrac{2\pi ijl}{d}}\operatorname{Im}(g_j\overline{g_{j+\frac{d}{2}}}),&\text{if } l \text{ is odd}.\end{cases}
\end{align*}
For $0\leq j\leq\frac{d}{2}-1$, we obtain
\begin{equation*}
g_j\overline{g_{j+\frac{d}{2}}}=\begin{cases}1-i\cdot\sqrt{2^{4j+d}-1},&\text{if } j\in\{0,1\},\\
							2-i\cdot\sqrt{2^{4j+d}-4},&\text{otherwise}.
							\end{cases}
\end{equation*}
When $l$ is even, we therefore obtain
\begin{align*}
\frac{1}{2}V_gg\left(\frac{d}{2},l\right)&=\sum\limits_{j=2}^{\frac{d}{2}-1} 2e^{-\nicefrac{2\pi ijl}{d}} + 1+e^{-\nicefrac{2\pi il}{d}}\\
						&=2\sum\limits_{j=0}^{\frac{d}{2}-1} e^{-\nicefrac{2\pi ijl}{d}} -1-e^{-\nicefrac{2\pi il}{d}}\\
						&=2\frac{e^{-\pi il}-1}{e^{-\nicefrac{2\pi il}{d}}-1}-1-e^{-\nicefrac{2\pi il}{d}}=-1-e^{-\nicefrac{2\pi il}{d}},
\end{align*}
which yields
\begin{equation*}
V_gg\left(\frac{d}{2},l\right)=0\quad\Leftrightarrow\quad e^{-\nicefrac{2\pi il}{d}}=-1\quad\Leftrightarrow\quad l=\frac{d}{2}.
\end{equation*}
On the other hand, we compute
\begin{equation*}
|\operatorname{Im}(g_j\overline{g_{j+\frac{d}{2}}})|=\begin{cases}\sqrt{2^{4j+d}-1},&\text{if } j\in\{0,1\},\\\sqrt{2^{4j+d}-4},&\text{if } 2\leq j\leq\frac{d}{2}-1.\end{cases}
\end{equation*}
Hence, when $l$ is odd, we can use the inequality $\sqrt{a-b}\geq\sqrt{a}-\sqrt{b}$ for $a>b>0$ to obtain
\begin{align*}
\frac{1}{2}\left|V_gg\left(\frac{d}{2},l\right)\right|&\geq\sqrt{2^{3d-4}-4}-\sum\limits_{j=0}^{\frac{d}{2}-2}\sqrt{2^{4j+d}-1}\geq\sqrt{2^{3d-4}-4}-\sum\limits_{j=0}^{\frac{d}{2}-2}2^{2j+\nicefrac{d}{2}}\\
								&=\sqrt{2^{3d-4}-4}-2^{\nicefrac{d}{2}}\cdot\frac{4^{\nicefrac{d}{2}-1}-1}{3}=\sqrt{2^{3d-4}-4}-\frac{1}{3}\cdot 2^{\nicefrac{3d}{2}-2}+\frac{1}{3}\cdot 2^{\nicefrac{d}{2}}\\
								&\geq 2^{\nicefrac{3d}{2}-2}-2-\frac{1}{3}\cdot 2^{\nicefrac{3d}{2}-2}+\frac{1}{3}\cdot 2^{\nicefrac{d}{2}}=\frac{2}{3}\cdot 2^{\nicefrac{3d}{2}-2}+\frac{1}{3}\cdot 2^{\nicefrac{d}{2}}-2\\
								&\geq\frac{2}{3}\cdot 2^4+\frac{4}{3}-2=10.
\end{align*}
Therefore, it follows $V_gg\left(\frac{d}{2},l\right)\neq 0$ for every $l\neq\frac{d}{2}$.

It remains to show that $V_gg(k,l)\neq 0$ whenever $k\neq \frac{d}{2}$. Then, Theorem \ref{thm:recwod2} will imply that $g$ does phase retrieval.\\
Directly from the definition of $g$, it follows $|g_j|=2^j$ for every $0\leq j\leq d-1$.
For $k<\frac{d}{2}$, we therefore obtain
\begin{align*}
|V_gg(k,l)|&\geq |g_{d-1}\overline{g_{d-1-k}}|-\sum\limits_{j=0}^{d-2}|g_j\overline{g_{j-k}}|\\
		&=2^{2d-2-k}-\sum\limits_{j=0}^{k-1} 2^{2j-k+d}-\sum\limits_{j=k}^{d-2} 2^{2j-k}\\
		&=2^{2d-2-k}-\sum\limits_{j=0}^{k-1} 2^{2j-k+d}-\sum\limits_{j=0}^{d-k-2} 2^{2j+k}\\
		&=2^{2d-2-k}-2^{d-k}\sum\limits_{j=0}^{k-1}4^j-2^k\sum\limits_{j=0}^{d-k-2}4^j\\
		&=2^{-k}\cdot 4^{d-1}-2^{d-k}\cdot\frac{4^k-1}{3}-2^k\cdot\frac{4^{d-k-1}-1}{3}.
\end{align*}
The last expression is positive if and only if
\begin{align*}
&3\cdot 4^{d-1}-4^{\nicefrac{d}{2}}\cdot(4^k-1)-4^k\cdot(4^{d-k-1}-1)>0\\
\Leftrightarrow\quad&3\cdot 4^{d-1}-4^{\nicefrac{d}{2}+k}+4^{\nicefrac{d}{2}}-4^{d-1}+4^k>0\\
\Leftrightarrow\quad&2\cdot 4^{d-1}+4^{\nicefrac{d}{2}}+4^k-4^{\nicefrac{d}{2}+k}>0.
\end{align*}
In order to show that the last inequality holds true, we will use the assumption that $k<\frac{d}{2}$, which implies
\begin{equation*}
2\cdot 4^{d-1}+4^{\nicefrac{d}{2}}+4^k-4^{\nicefrac{d}{2}+k}\geq 4^{d-1}+4^{\nicefrac{d}{2}}+4^k>0.
\end{equation*}
Finally, it follows $V_gg(k,l)\neq 0$ for $k>\frac{d}{2}$ by (\ref{symZd}), which concludes the proof.
\end{proof}

When $d\geq 5$ is odd, we will still be able to prove that there are windows that do phase retrieval but whose ambiguity functions have zeros. However, we will not construct a concrete window as in Example \ref{ex:PRVgg}. Instead, we will see that there are actually plenty of windows satisfying this property.
As before, we will first prove a reconstruction result. Unlike for Theorem \ref{thm:recwod2}, the proof is already non-algorithmic.

\begin{theorem}\label{thm:recwo0l}
Let $d\geq 5$ and $1\leq l^{\ast}\leq d-1$ be such that $l^{\ast}$ and $d$ are coprime or $d=6$ and $l^{\ast}=2$. Furthermore, consider $f,\tilde{f}\in\C^d$ such that $V_ff(k,l)=V_{\tilde{f}}\tilde{f}(k,l)$ holds for every $(k,l)\in(\Z_d\times\Z_d)\setminus\{(0,l^{\ast}),(0,-l^{\ast})\}$. Then there exists $\gamma\in\T$ such that $\tilde{f}=\gamma f$.
\end{theorem}

\begin{proof}
We need to show that $V_ff(0,l^{\ast})=V_{\tilde{f}}\tilde{f}(0,l^{\ast})$. In this case, we obtain by (\ref{symZd}) that\linebreak $V_ff(k,l)=V_{\tilde{f}}\tilde{f}(k,l)$ holds for all $(k,l)\in\Z_d\times\Z_d$ and exactly as in Theorem \ref{thm:Vggnonvan}, it follows $\tilde{f}=\gamma f$ for some $\gamma\in\T$.\\
Let $c:=V_ff(0,l^{\ast})-V_{\tilde{f}}\tilde{f}(0,l^{\ast})$ and assume $c\neq 0$. By assumption and (\ref{symZd}), it follows
\begin{equation*}
V_ff(0,\cdot)-V_{\tilde{f}}\tilde{f}(0,\cdot)=c\cdot\delta_{l^{\ast}}+\overline{c}\cdot\delta_{-l^{\ast}}.
\end{equation*}
Fourier inversion implies
\begin{align*}
\left|f_j\right|^2-\left|\tilde{f}_j\right|^2	&=\left(\mathcal{F}^{-1}(c\cdot\delta_{l^{\ast}}+\overline{c}\cdot\delta_{-l^{\ast}})\right)(j)\\
							&=\frac{1}{d}\left(c\cdot e^{\nicefrac{2\pi ijl^{\ast}}{d}}+\overline{c}\cdot e^{-\nicefrac{2\pi ijl^{\ast}}{d}}\right)\\
							&=\frac{2}{d}\cdot\operatorname{Re}\left(c\cdot e^{\nicefrac{2\pi ijl^{\ast}}{d}}\right).
\end{align*}
Since $l^{\ast}$ and $d$ are coprime, the numbers $e^{\nicefrac{2\pi ijl^{\ast}}{d}}$ are all distinct for $0\leq j\leq d-1$. Thus, we obtain
\begin{equation}\label{max2id}
\left|\left\{0\leq j\leq d-1~\middle|~\left|f_j\right|=\left|\tilde{f}_j\right|\right\}\right|\leq 2
\end{equation}
(and even $\leq 1$ when $d$ is odd).\\
When $d=6$ and $l^{\ast}=2$, then $\operatorname{Re}\left(c\cdot e^{\nicefrac{2\pi ijl^{\ast}}{d}}\right)=0$ can hold for at most one $0\leq j\leq 2$ and at most one $3\leq j\leq 5$. Thus, (\ref{max2id}) follows as well.\\
Now, let $0\leq j\leq d-1$ and assume $f_j=0$. This implies $0=f_j\overline{f_{j-k}}=\tilde{f}_j\overline{\tilde{f}_{j-k}}$ for every $k\neq 0$. Hence, it holds either $\tilde{f}_j=0$ or $\tilde{f}_m=0$ for every $m\neq j$.\\
Suppose that $\tilde{f}_j\neq 0$, i.e. $\tilde{f}_m=0$ for every $m\neq j$. If there were two distinct indices $m_1,m_2\neq j$ such that both $f_{m_1}\neq 0$ and $f_{m_2}\neq 0$, it would follow
\begin{equation*}
0=\tilde{f}_{m_1}\overline{\tilde{f}_{m_2}}=f_{m_1}\overline{f_{m_2}}\neq 0.
\end{equation*}
Therefore, both $f$ and $\tilde{f}$ can each have at most one non-zero entry, which contradicts (\ref{max2id}). Thus, we obtain $\tilde{f}_j=0$.\\
Using the symmetry of the problem, we obtain the equivalence
\begin{equation*}
f_j=0\quad\Leftrightarrow\quad \tilde{f}_j=0.
\end{equation*}
Together with (\ref{max2id}) and $d\geq 5$, it follows that there are at least three distinct entries\linebreak $0\leq j_1,j_2,j_3\leq d$ satisfying $0\neq\left|f_{j_n}\right|\neq\left|\tilde{f}_{j_n}\right|\neq 0$ for every $n\in\{1,2,3\}$.\\
In particular, it holds
\begin{equation*}
\alpha:=\frac{\tilde{f}_{j_1}}{f_{j_1}}\notin\T.
\end{equation*}
For $m\in\{j_2,j_3\}$, we obtain
\begin{equation*}
f_m\overline{f_{j_1}}=\tilde{f}_m\overline{\tilde{f}_{j_1}}\quad\Rightarrow\quad \tilde{f}_m=\frac{1}{\overline{\alpha}}\cdot f_m.
\end{equation*}
On the other hand,
\begin{equation*}
f_{j_3}\overline{f_{j_2}}=\tilde{f}_{j_3}\overline{\tilde{f}_{j_2}}\quad\Rightarrow\quad\tilde{f}_{j_3}=f_{j_3}\cdot\frac{\overline{f_{j_2}}}{\overline{\tilde{f}_{j_2}}}=\alpha\cdot f_{j_3}.
\end{equation*}
Together, this implies $\alpha=\frac{1}{\overline{\alpha}}$, which is a contradiction since $\alpha\notin\T$.\\
Therefore, it must hold $c=0$, which concludes the proof.
\end{proof}

As before, we will now show that there are indeed windows $g\in\C^d$ satisfying\linebreak $\Omega(g)=(\Z_d\times\Z_d)\setminus\{(0,l^{\ast}),(0,-l^{\ast})\}$ for some $1\leq l^{\ast}\leq d-1$ which is coprime to $d$.
We begin with a lemma concerning the choice of $l^{\ast}$.

\begin{lemma}\label{lm:last}
Let $d\geq 5$ and define
\begin{equation*}
l^{\ast}:=\begin{cases}2,&\text{if } d=6,\\ \frac{d-1}{2},&\text{if } d\text{ is odd},\\ \frac{d}{2}-1,&\text{if } 4\mid d,\\ \frac{d}{2}-2,&\text{otherwise}.\end{cases}
\end{equation*}
Then the following hold true.
\begin{enumerate}[a)]
\item It holds $\frac{d}{4}<l^{\ast}<\frac{3d}{4}$.
\item When $d\neq 6$, then $l^{\ast}$ and $d$ are coprime.
\end{enumerate}
\end{lemma}

\begin{proof}\mbox{}
\begin{enumerate}[a)]
\item Obviously, it holds $l^{\ast}<\frac{d}{2}<\frac{3d}{4}$ in all cases. When $d=6$, we obtain $l^{\ast}=2>\frac{3}{2}=\frac{d}{4}$. When $d\in\{5,7,8\}$ is odd, it holds
\begin{equation*}
l^{\ast}>\frac{d}{2}-1=\frac{d}{2}-\frac{4}{4}>\frac{d}{2}-\frac{d}{4}=\frac{d}{4}.
\end{equation*}
Finally, when $d>8$, it follows
\begin{equation*}
l^{\ast}\geq\frac{d}{2}-2=\frac{d}{2}-\frac{8}{4}>\frac{d}{2}-\frac{d}{4}=\frac{d}{4}.
\end{equation*}
\item First, assume that $d$ is odd. If $r$ is a common divisor of $\frac{d-1}{2}$ and $d$, it follows that $r$ is a common divisor of $d-1$ and $d$, which implies $r=1$.\\
Next, we consider the case that $d$ is a multiple of $4$ and let $r$ be a common divisor of $\frac{d}{2}-1$ and $d$. Then, $r$ is a common divisor of $d-2$ and $d$, which yields $r\in\{1,2\}$. Since $\frac{d}{2}-1$ is odd, we obtain $r=1$.\\
Finally, let ($d\neq 6$ and) $d\equiv 2~\operatorname{mod} 4$ and suppose that $r$ is a common divisor of $\frac{d}{2}-2$ and $d$. In this case, $r$ is a common divisor of $d-4$ and $d$, which gives $r\in\{1,2,4\}$. Once again, it follows $r=1$ because $\frac{d}{2}-2$ is odd.
\end{enumerate}
\end{proof}

Lemma \ref{lm:last} now allows us to prove the existence of windows with the desired properties.

\begin{theorem}\label{thm:PRVggodd}
Let $d\geq 5$ and $l^{\ast}$ as in Lemma \ref{lm:last}. Then there exists a window $g\in\C^d$ satisfying the following conditions.
\begin{enumerate}[(i)]
\item $g$ does phase retrieval.
\item It holds $\Omega(g)=(\Z_d\times\Z_d)\setminus\{(0,l^{\ast}),(0,-l^{\ast})\}$.
\item It holds $\operatorname{supp}(g)\subseteq\left\{0,\dots,\left\lfloor\frac{d}{2}\right\rfloor\right\}$.
\end{enumerate}
\end{theorem}

\begin{proof}
Let $m:=\left\lfloor\frac{d}{2}\right\rfloor$. First, we compute
\begin{equation*}
\left(z-e^{-\nicefrac{2\pi il^{\ast}}{d}}\right)\cdot\left(z-e^{\nicefrac{2\pi il^{\ast}}{d}}\right)=z^2-\left(e^{\nicefrac{2\pi il^{\ast}}{d}}+e^{-\nicefrac{2\pi il^{\ast}}{d}}\right)+1=z^2-2\operatorname{Re}\left(e^{\nicefrac{2\pi il^{\ast}}{d}}\right)\cdot z+1.
\end{equation*}
Since $\frac{d}{4}<l^{\ast}<\frac{3d}{4}$ by Lemma \ref{lm:last}, all coefficients of the polynomial
\begin{equation*}
p(z):=\sum\limits_{j=0}^m a_j z^j:=\left(z-e^{-\nicefrac{2\pi il^{\ast}}{d}}\right)\cdot\left(z-e^{\nicefrac{2\pi il^{\ast}}{d}}\right)\cdot(z+2)^{m-2}
\end{equation*}
are thus real and positive. Furthermore, $e^{-\nicefrac{2\pi il^{\ast}}{d}}$ and $e^{\nicefrac{2\pi il^{\ast}}{d}}$ are the only zeros of $p$ in $\T$. By letting $c_j:=\sqrt{a_j}$ for every $0\leq j\leq m$, any window $g$ satisfying
\begin{equation*}
|g_j|=\begin{cases} c_j,&\text{if } 0\leq j\leq m,\\ 0,&\text{otherwise},\end{cases}
\end{equation*}
already fulfills $V_gg(0,l^{\ast})=V_gg(0,-l^{\ast})=0$ and $V_gg(0,l)\neq 0$ for every $l\in\{0,\dots,d-1\}\setminus\{l^{\ast},-l^{\ast}\}$. Since (ii) implies (i) by Theorem \ref{thm:recwo0l} and Lemma \ref{lm:Vff}, it remains to show that $V_gg(k,l)\neq 0$ whenever $k\neq 0$. By (\ref{symZd}), it suffices to prove this for $k\leq\frac{d}{2}$.\\
Define $g_j:=c_j$ for every $0\leq j\leq m-1$ and let $g_{m}\in\C$ be arbitrary.\\
When $d$ is even, i.e. $m=\frac{d}{2}$, we obtain
\begin{equation*}
V_gg\left(\frac{d}{2},l\right)=g_0\overline{g_m}+(-1)^l\cdot g_m\overline{g_0}\in\{2\operatorname{Re}(g_m\overline{g_0}),-2i\operatorname{Im}(g_m\overline{g_0})\}.
\end{equation*}
Since $g_0>0$, it may only hold $V_gg\left(\frac{d}{2},l\right)=0$, when $g_m\in\R\cup i\R$.\\
Now, consider $k<\frac{d}{2}$. In this case, we obtain
\begin{equation*}
0=V_gg(k,l)=\sum\limits_{j=k}^{m} g_j\overline{g_{j-k}} e^{-\nicefrac{2\pi ijl}{d}}\quad\Leftrightarrow\quad g_m=\frac{e^{\nicefrac{2\pi ilm}{d}}}{\overline{g_{m-k}}}\cdot\left(\sum\limits_{j=k}^{m-1}g_j\overline{g_{j-k}}e^{-\nicefrac{2\pi ijl}{d}}\right)=:s_{k,l}
\end{equation*}
Obviously, the set
\begin{equation*}
M:=\{s\in\C~|~|s|=c_m\}\setminus\left(\{c_m,-c_m,ic_m,-ic_m\}\cup\left\{s_{k,l}~|~1\leq k<\frac{d}{2},0\leq l\leq d-1\right\}\right)
\end{equation*}
is non-empty. By the previous calculations, any choice of $g_m\in M$ results in a window $g$ satisfying the desired properties.
\end{proof}

Note that (when $d$ is even) the windows obtained in the proof of Theorem \ref{thm:PRVggodd} are clearly distinct from Example \ref{ex:PRVgg}. Thus, the Theorem is not only valuable when $d$ is odd.
In fact, windows of finite length are interesting on their own, with respect to phase retrieval. We will analyze their behaviour in the following subsections.

The following corollary summarizes the results of the current subsection.

\begin{corollary}\label{cor:Vgg}
The following are equivalent.
\begin{enumerate}[(i)]
\item There exists a window $g\in\C^d$ satisfying $\Omega(g)\subsetneq\Z_d\times\Z_d$ that does phase retrieval.
\item It holds $d\geq 4$.
\end{enumerate}
\end{corollary}

Intuitively, we would expect the loss of a single pair $V_ff(k,l), V_ff(-k,-l)$ in the finite setting to be more severe than the loss of $V_ff\vert_{\{k\}\times\Omega\,\cup\,\{-k\}\times\overline{\Omega}}$ in the infinite-dimensional setting, where $\Omega\subseteq\T$ is a set of positive, but potentially small measure.\\
However, a comparison between Corollaries \ref{cor:prZnec} and \ref{cor:Vgg} suggests otherwise. This emphasizes the fact that Corollary \ref{cor:prZnec} is more about the behaviour of power series (and thus the possible values of $V_gg$ for one-sided $g$) than about actual recovery requirements for the phase retrieval problem associated with $g$.

\subsection{Windows of length $L+1$: The generic case}

For the remainder of the section, we will step aside from the condition $\Omega(g)=\Z_d\times\Z_d$ and will instead focus on short windows which we have seen to allow a very straight-forward analysis of phase retrieval in the infinite-dimensional setting.\\
We fix $L<\frac{d}{2}$ and make the following definition.

\begin{definition}\label{def:lengthL+1Zd}
The set of all windows of length $L+1$ is defined as
\begin{equation*}
\C_L^d:=\left\{g\in\C^d~\middle|~g_j\neq 0~\Leftrightarrow~0\leq j\leq L\right\}.
\end{equation*}
\end{definition}

As before, fixing the support to $\{0,\dots,L\}$ is just for convenience. Lemma \ref{lm:windowshift} ensures that the results of this section also apply when $\operatorname{supp}(g)\subseteq\{j_0,\dots,j_0+L\}$ holds for some $0\leq j_0\leq d-1$.\\
The main benefit of the class $\C_L^d$ is that there is no ``cyclic overlap'', i.e. we can write
\begin{equation*}
V_gg(k,l)=\sum\limits_{j=k}^{L} g_j\overline{g_{j-k}}e^{\nicefrac{-2\pi ijl}{d}}
\end{equation*}
for every $g\in\C_L^d$ and all $0\leq k\leq L$, $0\leq l\leq d-1$, without considering the indices modulo $d$.\\

Short windows $g\in\C_L^d$ have been studied e.g. in \cite[Theorem III.1]{EHJ}, where it is shown that almost all non-vanishing signals (i.e. signals $f\in\C^d$ fulfilling $f_j\neq 0$ for every $0\leq j\leq d-1$) are phase retrievable w.r.t. a fixed window $g\in\C_L^d$.\\
In this section, we will see that it is possible to ``flip'' the ``almost every'' part of the statement to the window, i.e. that almost every window $g\in\C_L^d$ does phase retrieval (for \textit{all} non-vanishing signals $f\in\C^d$).
This follows as a special instance of Theorem \ref{thm:Cdgeneric}, when combined with Lemma \ref{lm:Omegag}, and can be seen as a generalization of Theorem \ref{thm:nonvangen}.\\
In \cite[Corollary III.1]{EHJ}, the authors also prove a result on so-called sparse signals (i.e. signals that are not non-vanishing). Note however, that this result doesn't apply to our setting since we only consider the case $L=1$ (where $L$ is the separation parameter from \cite{EHJ}).\\

Because of Lemma \ref{lm:Vff}, phase retrieval is fully characterized by $\Omega(g)$. We will thus take a close look at this set for windows of length $L+1$.

\begin{lemma}\label{lm:Omegag}
Let $g\in\C_L^d$.
\begin{enumerate}[a)]
\item It holds $\Omega(g)\subseteq\Omega_L^d:=\{(k,l)\in\Z_d\times\Z_d~|~0\leq k\leq L~\text{or}~d-L\leq k\leq d-1\}$
\item It holds $\Omega(g)=\Omega_L^d$ for almost every $g\in\C_L^d$, i.e. there exists a set $E\subseteq\C^{L+1}$ of $(L+1)$-dimensional Lebesgue measure zero such that $\Omega(g)=\Omega_L^d$ holds for all $g\in\C_L^d$ satisfying\linebreak $g\vert_{\{0,\dots,L\}}\in\C^{L+1}\setminus E$.
\end{enumerate}
\end{lemma}

\begin{proof}\mbox{}
\begin{enumerate}[a)]
\item Let $L<k<d-L$ and $0\leq j\leq d-1$. Since $L<\frac{d}{2}$, either $g_j=0$ or $g_{j-k}=0$. Thus,
\begin{equation*}
V_g g(k,l)=\sum\limits_{j=0}^{d-1} g_j\overline{g_{j-k}} e^{\nicefrac{2\pi i jl}{d}}=0
\end{equation*}
for all $0\leq l\leq d-1$.
\item We adapt the proof of \cite[Proposition 2.1]{BF}, which we have already referred to for Theorem \ref{thm:nonvangen}.
It suffices to show that
\begin{equation*}
E_{k,l}:=\{g\vert_{\{0,\dots,d\}} ~|~ g\in\C_L^d, V_g g(k,l)=0\}
\end{equation*}
has measure zero for every $(k,l)\in\Omega_L^d$. Using (\ref{symZd}), we may w.l.o.g. assume $0\leq k\leq L$.\\ It follows
\begin{align*}
E_{k,l}	&=\left\{g\vert_{\{0,\dots,L\}} ~\middle|~ g\in\C_L^d, V_g g(k,l)=0\right\}\\
		&=\left\{\tilde{g}\in\C^{L+1}~\middle|~\sum\limits_{j=k}^L \tilde{g}_j \overline{\tilde{g}_{j-k}}e^{\nicefrac{-2\pi ijl}{d}}\right\}\\
		&=\left\{\tilde{g}\in\C^{L+1}~\middle|~\tilde{g}_k=0\right\}\cup\left\{\tilde{g}\in\C^{L+1}~\middle|~\tilde{g}_k\neq 0, \tilde{g}_0=-\frac{1}{\overline{\tilde{g}_k}}\sum\limits_{j=k+1}^L \overline{\tilde{g}_j} \tilde{g}_{j-k}e^{\nicefrac{2\pi ijl}{d}}\right\}.
\end{align*}
Thus, $E_{k,l}$ is the union of two sets of measure zero and therefore has measure zero itself.
\end{enumerate}
\end{proof}

Part b) of Lemma \ref{lm:Omegag} tells us that $\Omega(g)=\Omega_L^d$ should be considered to be the ``generic case": Picking $g\in\C_L^d$ randomly, guarantees $g$ to fulfill this property with probability $1$. We will present a very straight-forward example for a generic window.

\begin{example}\label{ex:genericg}
Let $g\in\C_L^d$ be such that
\begin{equation*}
|g_j|:=\begin{cases} 2^j,&\text{if }0\leq j\leq L,\\0,&\text{otherwise}.\end{cases}
\end{equation*}
Then, it holds $\Omega(g)=\Omega_L^d$.
Additionally, when $d$ is odd and $L=\frac{d-1}{2}$, then $g$ does phase retrieval. In particular, by choosing $g_j:=2^j$ when $0\leq j\leq\frac{d-1}{2}$, one obtains a real-valued window doing phase retrieval.
\end{example}

\begin{proof}
In order to show the first statement, it suffices by (\ref{symZd}) to show that $V_gg(k,l)\neq 0$ holds true whenever $(k,l)\in\Omega_L^d$ and $0\leq k\leq\frac{d}{2}$, i.e. $0\leq k\leq L$. This follows from
\begin{align*}
|V_gg(k,l)|	&=\left|\sum\limits_{j=0}^{d-1} g_j\overline{g_{j-k}} e^{-\nicefrac{2\pi ijl}{d}}\right|=\left|\sum\limits_{j=k}^L g_j\overline{g_{j-k}}e^{-\nicefrac{2\pi ijl}{d}}\right|\\
		&\geq|g_L|\cdot|g_{L-k}|-\sum\limits_{j=k}^{L-1} |g_j|\cdot|g_{j-k}|\\
		&=2^{2L-k}-\sum\limits_{j=k}^{L-1}2^{2j-k}\\
		&\geq\frac{1}{2^k}\left(4^L-\sum\limits_{j=0}^{L-1} 4^j\right)\\
		&=\frac{1}{2^k}\left(4^L-\frac{4^L-1}{3}\right)\geq\frac{1}{2^k}>0.
\end{align*}
When $d$ is odd and $L=\frac{d-1}{2}$, it follows $d-L=\frac{d+1}{2}=L+1$ and thus, $\Omega_L^d=\Z_d\times\Z_d$. Now, Theorem \ref{thm:Vggnonvan} implies that $g$ does phase retrieval.
\end{proof}

The question whether there exists a real-valued window doing phase retrieval, is interesting on its own. Recall that it is easy to construct such windows in the infinite-dimensional setting (by applying Theorem \ref{thm:Zoneside} or Theorem \ref{thm:Zgeneral}). For finite phase retrieval, it turns out that this is exclusive to odd dimensions.

\begin{theorem}\label{thm:realwindowdeven}
If $d$ is even, there exists no real-valued window $g\in\R^d$ doing phase retrieval.
\end{theorem}

\begin{proof}
Let $d$ be even, $g\in\R^d$ and assume that $0\leq l\leq d-1$ is odd. Then,
\begin{align*}
V_gg\left(\frac{d}{2},l\right)	&=\sum\limits_{j=0}^{d-1}g_jg_{j-\frac{d}{2}}e^{-\nicefrac{2\pi ijl}{d}}\\
					&=\sum\limits_{j=0}^{\frac{d}{2}-1}g_jg_{j+\frac{d}{2}}e^{-\nicefrac{2\pi ijl}{d}}~+~\sum\limits_{j=\frac{d}{2}}^{d-1} g_jg_{j-\frac{d}{2}}e^{-\nicefrac{2\pi ijl}{d}}\\
					&=\sum\limits_{j=0}^{\frac{d}{2}-1} \left(1+(-1)^l\right)\cdot g_j g_{j+\frac{d}{2}}e^{-\nicefrac{2\pi ijl}{d}}\\
					&=0.
\end{align*}
Thus, it holds $\left(\frac{d}{2},l\right)\notin\Omega(g)$ whenever $l$ is odd. Now, consider $f,\tilde{f}\in\C^d$ given by
\begin{equation*}
f_j:=\begin{cases} 1,&\text{if } j=0,\\ 1+i,&\text{if } j=\frac{d}{2},\\0,&\text{otherwise}\end{cases}\quad\quad\text{and}\quad\quad\tilde{f}_j:=\begin{cases} 1,&\text{if } j=0,\\1-i,&\text{if } j=\frac{d}{2},\\0,&\text{otherwise}.\end{cases}
\end{equation*}
Clearly, it holds $\tilde{f}\notin\T\cdot f$.\\
Since $|f|=\left|\tilde{f}\right|$, we obtain $V_ff(0,l)=V_{\tilde{f}}\tilde{f}(0,l)$ for every $0\leq l\leq d-1$. Moreover, it holds $V_ff(k,l)=V_{\tilde{f}}\tilde{f}(k,l)=0$ whenever $k\notin\left\{0,\frac{d}{2}\right\}$. Finally, let $0\leq l\leq d-1$ be even.
It follows
\begin{equation*}
V_ff\left(\frac{d}{2},l\right)=f_0\overline{f_{\frac{d}{2}}}+\overline{f_0}f_{\frac{d}{2}}=1-i+1+i=2=1+i+1-i=\tilde{f}_0\overline{\tilde{f}_{\frac{d}{2}}}+\overline{\tilde{f}_0}\tilde{f}_{\frac{d}{2}}.
\end{equation*}
Hence, it holds $V_ff(k,l)=V_{\tilde{f}}\tilde{f}(k,l)$ for every $(k,l)\in\Omega(g)$ and therefore, $g$ can't do phase retrieval by Lemma \ref{lm:Vff}.
\end{proof}

We will now return to windows of length $L+1$ and prove an easy characterization for phase retrieval in the generic case, which is very similar to Theorem \ref{thm:ZlengthL+1}. First, we need a cyclic version of $L$-connectivity.

\begin{definition}\label{def:LconnecZd} Consider a subset $M\subseteq\Z_d$ and a signal $f\in\C^d$.
\begin{enumerate}[a)]
\item Two indices $j,k\in M$ are called $L$-close $(\operatorname{mod} d)$, iff
\begin{equation*}
\operatorname{dist}_d(j,k):=\min\limits_{n\in\Z} |j-(k+nd)|\leq L.
\end{equation*}
\item Two indices $j,k\in M$ are called $L$-connected $(\operatorname{mod} d)$ and we write $j\sim_{d,L} k$, iff there exist indices $l_0,\dots,l_n\in M$ satisfying $l_0=j$ and $l_n=k$, such that $l_m$ and $l_{m+1}$ are $L$-close $(\operatorname{mod} d)$ for every $0\leq m<n$.
\item Clearly, $\sim_{d,L}$ defines an equivalence relation on $M$. The equivalence classes of $M$ are called $L$-connectivity components of $M$, and $M$ is called $L$-connected $(\operatorname{mod} d)$, iff there exists at most one connectivity component.
\item $f$ is called $L$-connected $(\operatorname{mod} d)$, iff $\operatorname{supp}(f)$ is $L$-connected $(\operatorname{mod} d)$. The $L$-connectivity components of $\operatorname{supp}(f)$ are called $L$-connectivity components of $f$.
\end{enumerate}
\end{definition}

Again, a set (or a signal) is thus $L$-connected, iff it is connected up to ``holes'' of length at most $L-1$, when considered modulo $d$. The ``modulo $d$'' part is obviously important. Note e.g. that $(0,1,0,0,1)^t\in\C^5$ is $2$-connected by definition.\\
Definition \ref{def:LconnecZd} allows us to prove a version of Theorem \ref{thm:ZlengthL+1} for generic short windows in $\C_L^d$. The necessary conditions for phase retrieval also hold for general windows of length $L+1$ (cf. subsection \ref{subs:L+1general}) and have already been observed in \cite[Theorem 2]{BCEMS}.

\begin{theorem}\label{thm:Cdgeneric}
Let $g\in\C_L^d$ such that $\Omega(g)=\Omega_L^d$. Furthermore, consider a signal $f\in\C^d$ and let $C_1,\dots,C_n$ be the $L$-connectivity components of $f$. Then, the following statements hold true.
\begin{enumerate}[a)]
\item A signal $\tilde{f}\in\C^d$ satisfies $\left|V_g f\right|=\left|V_g\tilde{f}\right|$, if and only if $\operatorname{supp}\left(\tilde{f}\right)=\operatorname{supp}(f)$ and for every $1\leq m\leq n$, there exists $\gamma_m\in\T$ such that $\tilde{f}\vert_{C_m}=\gamma_mf\vert_{C_m}$.
\item $f$ is phase retrievable w.r.t. $g$ if and only if $f$ is $L$-connected.
\end{enumerate}
\end{theorem}

\begin{proof}\mbox{}
This is proved analogously to Theorem \ref{thm:ZlengthL+1}.
\end{proof}

Once again, the measurement $|V_gf|$ always tells us whether the underlying signal is phase retrievable (cf. Remark \ref{rem:ZlengthL+1}).\\

For $L=2$, the sufficient condition for phase retrieval in Theorem \ref{thm:Cdgeneric} is quite similar to\linebreak \cite[Theorem 2.4]{BF}. Note however, that $1$-connectivity is strictly weaker than assuming the signal to vanish nowhere. On the other hand, the result in \cite{BF} is not restricted to short windows.\\
Whenever $d$ and $L$ are coprime, it is shown in \cite[Theorem 1]{BCEMS} that it suffices to assume $V_gg(0,l)\neq 0$ for all $0\leq l\leq d-1$, in order to guarantee phase retrieval for all non-vanishing signals. In \cite[Corollary 3.2]{CHLSS}, a similar result has been proved for windows $g$ satisfying the weaker assumption $\operatorname{supp}(g)\subseteq\{0,\dots,L\}$.\\
On the other hand, Theorem \ref{thm:PRVggodd} shows that (at least when $d\geq 5$) $V_gg(0,l)\neq 0$ for all $0\leq l\leq d-1$ is not even a necessary condition for (global) phase retrieval.

\subsection{Windows of length $L+1$: The general case}\label{subs:L+1general}

Our intuition in understanding Theorem \ref{thm:Cdgeneric} is that phase retrieval depends on whether or not $g$ is ``long enough'' to connect certain areas of $f$. Since $\Omega(g)\subseteq\Omega_L^d$ holds for all $g\in\C_L^d$, this remains a necessary condition for the general case, but is no longer sufficient, as shows the following example.

\begin{example}\label{ex:cexprshift}
Let $r\geq 2$ be a common divisor of $d$ and $L+1$. Consider
\begin{equation*}
g_j:=\begin{cases} 1,&\text{if }0\leq j\leq L,\\0,&\text{otherwise},\end{cases}
\end{equation*}
as well as
\begin{equation*}
f_j:=\begin{cases} 1,&\text{if } j\equiv 0~ (\operatorname{mod } r),\\0,&\text{otherwise}.\end{cases}
\end{equation*}
For $0\leq m<r$, let $f^{(m)}:=T_m f$.\\
Then, it holds $g\in\C_L^d$ as well as
\begin{equation*}
V_{f^{(m)}}f^{(m)}\vert_{\Omega(g)}=V_{f^{(n)}}f^{(n)}\vert_{\Omega(g)}
\end{equation*}
for all $0\leq m<n<r$, and there is no $\gamma\in\T$ such that $f^{(n)}=\gamma f^{(m)}$. Thus, $f=f^{(0)}$ is not phase retrievable w.r.t. $g$.
However, if $r<L+1$, $f$ is $L$-connected. Choosing e.g. $d=8$, $L=3$ and $r=2$ shows that Theorem \ref{thm:Cdgeneric} does not hold true for general windows $g\in\C_L^d$.
\end{example}

\begin{proof}
For every $0\leq m<r$ and all $0\leq k,l\leq d-1$, it holds
\begin{equation*}
V_{f^{(m)}}f^{(m)}(k,l)=\sum\limits_{j=0}^{d-1}f^{(m)}_j\overline{f^{(m)}_{j-k}}e^{-\nicefrac{2\pi ijl}{d}}=\sum\limits_{j=0}^{d-1} f_{j-m}\overline{f_{j-m-k}} e^{-\nicefrac{2\pi ijl}{d}}=e^{-\nicefrac{2\pi iml}{d}}V_ff(k,l).
\end{equation*}
Therefore, we obtain
\begin{equation*}
V_{f^{(m)}}f^{(m)}\vert_{\Omega(g)}=V_{f^{(n)}}f^{(n)}\vert_{\Omega(g)},
\end{equation*}
as soon as either $l=0$ or $V_ff(k,l)=0$ holds for every $(k,l)\in\Omega(g)$. Because of (\ref{symZd}), we only have to consider the case $0\leq k\leq L$.  Clearly, it holds
\begin{equation*}
f_j\overline{f_{j-k}}=\begin{cases}1,&\text{if } r\mid j\text{ and }r\mid k,\\0,&\text{otherwise,}\end{cases}
\end{equation*}
and therefore $V_ff(k,l)=0$ whenever $k\nmid r$. When $k\mid r$, it follows
\begin{align*}
V_ff(k,l)	&=\sum\limits_{j=0}^{d-1}f_j\overline{f_{j-k}}e^{-\nicefrac{2\pi ijl}{d}}=\sum\limits_{s=0}^{\frac{d}{r}-1}e^{-\nicefrac{2\pi isrl}{d}}\\
		&=\sum\limits_{s=0}^{\frac{d}{r}-1}\left(e^{-\nicefrac{2\pi ilr}{d}}\right)^s=0,
\end{align*}
unless $\frac{d}{r}\mid l$. It remains to show that $(k,l)\notin\Omega(g)$ whenever $k\leq L$, $k\mid r$ and $\frac{d}{r}\mid l\neq 0$. Under these assumptions, let
\begin{equation*}
L+1=tr,\quad\quad k=sr,\quad\quad\text{and}\quad\quad l=v\frac{d}{r}.
\end{equation*}
It follows $0<v<r$ and thus
\begin{align*}
V_gg(k,l)	&=\sum\limits_{j=0}^{d-1} g_j\overline{g_{j-k}}e^{-\nicefrac{2\pi ijl}{d}}=\sum\limits_{j=k}^L g_j\overline{g_{j-k}}\left(e^{-\nicefrac{2\pi iv}{r}}\right)^j\\
		&=e^{-\nicefrac{2\pi ikv}{r}}\cdot\sum\limits_{j=0}^{L-k}\left(e^{-\nicefrac{2\pi iv}{r}}\right)^j=e^{-\nicefrac{2\pi ikv}{r}}\cdot\sum\limits_{j=0}^{(t-s)r-1}\left(e^{-\nicefrac{2\pi iv}{r}}\right)^j\\
		&=e^{-\nicefrac{2\pi ikv}{r}}\cdot\sum\limits_{n=0}^{(t-s)-1}\sum\limits_{j=nr}^{(n+1)r-1}\left(e^{-\nicefrac{2\pi iv}{r}}\right)^j=e^{-\nicefrac{2\pi ikv}{r}}\cdot\sum\limits_{n=0}^{(t-s)-1}\sum\limits_{\alpha=0}^{r-1}\left(e^{-\nicefrac{2\pi iv}{r}}\right)^{nr+\alpha}\\
		&=e^{-\nicefrac{2\pi ikv}{r}}\cdot(t-s)\cdot\sum\limits_{\alpha=0}^{r-1}\left(e^{-\nicefrac{2\pi iv}{r}}\right)^{\alpha}=0.
\end{align*}
The fact that there is no $\gamma\in\T$ such that $f^{(n)}=\gamma f^{(m)}$ follows immediately from\linebreak $\operatorname{supp}\left(\tilde{f}\right)\neq\operatorname{supp}(f)$ when $0\leq m<n<r$.
\end{proof}

Our goal is to find a large class of signals for which we can once more link phase retrievability to $L$-connectivity.
Since it makes essentially no difference for the calculations below, we will return to directly using the measurements $|V_g f|$ instead of $V_ff\vert_{\Omega(g)}$. Note that the computation and subsequent use of entries of $f\cdot\overline{T_k f}$ is still our main tool, and thus the basic idea stays the same.

We fix a window $g\in\C_L^d$, a signal $f\in\C^d$ and a measurement $X\in\C^{d\times d}$. For $0\leq k\leq L$ and $0\leq j\leq d-1$, consider the coefficients
\begin{align*}
a_j^{(k)}&:=f_j\overline{f_{j-k}},\\
c_j^{(k)}&:=g_j\overline{g_{j-k}}\quad\text{and}\\
b_j^{(k)}&:=\frac{1}{d^2}\sum\limits_{l=0}^{d-1} \hat{X}(l,-k)e^{\nicefrac{2\pi ijl}{d}}.
\end{align*}

Note that all $c^{(k)}_j$ are known and all $b^{(k)}_j$ can be computed based on the measurement, while it is our goal to recover all (or at least some of the) $a^{(k)}_j$.\\
In order to establish a relation between the various coefficients, we will use the Fourier inversion formula. It is therefore once again useful to compute the constant $c_G$.

\begin{lemma}\label{lm:cGZd}
When $G=\Z_d$, it holds $c_G=d$.
\end{lemma}

\begin{proof}
As in the proof of Lemma \ref{lm:cGZ}, it is enough to compute $\frac{\left\|\widehat{h}\right\|_2}{\|h\|_2}$ for $h=\delta_0$. We obtain $\widehat{h}(l)=1$ for every $0\leq l\leq 1$, which implies
\begin{equation*}
\left\|\widehat{h}\right\|_2=\sqrt{d}=\sqrt{d}\cdot\|h\|_2
\end{equation*}
and therefore $c_G=d$.
\end{proof}

The identity from Theorem \ref{thm:VgfvsVff} therefore becomes (cf. e.g. \cite[Proposition 3.13]{GKR})
\begin{equation}\label{VgfvsVffZd}
\widehat{\left|V_gf\right|^2}(l,k)=d\cdot V_ff(-k,l)\cdot\overline{V_gg(-k,l)}.
\end{equation}

This enables us to express the coefficients $b^{(k)}_j$ in terms of $a^{(k)}$ and $c^{(k)}$.

\begin{lemma}\label{lm:lincomb}
Assume that $|V_gf|^2=X$. Then,
\begin{equation*}
b_j^{(k)}=\sum\limits_{m=k}^L \overline{c_m^{(k)}}a_{m+j}^{(k)}
\end{equation*}
holds for every $0\leq k\leq L$ and $0\leq j\leq d-1$.
\end{lemma}

\begin{proof}
Since $\widehat{G}=\Z_d$ is finite, it obviously holds $V_ff(k,\cdot)\in L^1\left(\widehat{G}\right)$ for every $0\leq k\leq d-1$.
Using the identity (\ref{VgfvsVffZd}) as well as the inversion formula from Corollary \ref{cor:planch}, we obtain
\begin{align*}
b_j^{(k)}	&=\frac{1}{d^2}\cdot\sum\limits_{l=0}^{d-1} \widehat{|V_g f|^2}(l,-k)e^{\nicefrac{2\pi ijl}{d}}\\
		&=\frac{1}{d}\cdot\sum\limits_{l=0}^{d-1} V_ff(k,l)\cdot\overline{V_gg(k,l)}e^{\nicefrac{2\pi ijl}{d}}\\
		&=\frac{1}{d}\cdot\sum\limits_{l=0}^{d-1} V_ff(k,l)\cdot e^{\nicefrac{2\pi ijl}{d}}\cdot\left(\sum\limits_{m=k}^L \overline{c_m^{(k)}}e^{\nicefrac{2\pi iml}{d}}\right)\\
		&=\sum\limits_{m=k}^L\overline{c_m^{(k)}}\cdot\frac{1}{d}\cdot\sum\limits_{l=0}^{d-1} V_ff(k,l)e^{\nicefrac{2\pi i(j+m)l}{d}}\\
		&=\sum\limits_{m=k}^L\overline{c_m^{(k)}}a_{m+j}^{(k)}.
\end{align*}
\end{proof}

Lemma \ref{lm:lincomb} gives us $L+1$ systems of linear equations which may or may not allow us to recover $a^{(k)}$, i.e. $f\cdot\overline{T_kf}$ for a given $0\leq k\leq L$.

Fixing $k$, we observe that the right-hand side of every equation is given by a linear combination of $L-k+1$ consecutive entries of $a^{(k)}$ with non-zero coefficients $c_m^{(k)}$. Knowing $L-k$ of those $L-k+1$ entries for just a single index $j(k)$ therefore enables us to iteratively compute $a^{(k)}$.

Recovering $a^{(k)}=f\cdot\overline{T_kf}$ for \textit{every} $0\leq k\leq L$ would then once more allow us to compute the signal $f$ on each of its $L$-connectivity components.
This leads to the following corollary.

\begin{corollary}\label{cor:LinCombIt}
Let $\tilde{f}\in\C^d$ be such that $\left|V_g f\right|=\left|V_g\tilde{f}\right|$.\\
If for every $0\leq k\leq L$ there exists an index $j^{\ast}(k)\in\{0,\dots,d-1\}$, such that
\begin{equation*}
\left|\left\{j\in\{j^{\ast}(k),\dots,j^{\ast}(k)+L-k\}~\middle|~a_j^{(k)}\neq\tilde{f}_j\overline{\tilde{f}_{j-k}}\right\}\right|\leq 1,
\end{equation*}
then $\operatorname{supp}\left(\tilde{f}\right)=\operatorname{supp}(f)$ and for every $1\leq m\leq n$, there exists $\gamma_m\in\T$ such that\linebreak $\tilde{f}\vert_{C_m}=\gamma_mf\vert_{C_m}$.
\end{corollary}

Since $c_m^{(0)}=|g_m|^2> 0$ for every $0\leq m\leq L$, things are especially convenient when $k=0$: Here, the existence of $0\leq j\leq d-1$ such that $b_j^{(0)}=0$ is equivalent to $a_j^{(0)}=\dots=a_{j+d}^{(0)}=0$ and thus $f_j=\dots=f_{j+d}=0$.
Together with Corollary \ref{cor:LinCombIt}, we obtain the following version of Theorem \ref{thm:Cdgeneric} for general windows and signals with at least $L+1$ consecutive zeros.

\begin{theorem}\label{thm:CdHoleL+1}
Let there be $j^{\ast}\in\{0,\dots,d-1\}$ such that
\begin{equation*}
\operatorname{supp}(f)\cap\{j^{\ast},\dots,j^{\ast}+L\}=\emptyset,
\end{equation*}
and let $C_1,\dots,C_n$ be the $L$-connectivity components of $f$. Then, the following statements hold true.
\begin{enumerate}[a)]
\item A signal $\tilde{f}\in\C^d$ satisfies $\left|V_g f\right|=\left|V_g\tilde{f}\right|$, if and only if $\operatorname{supp}\left(\tilde{f}\right)=\operatorname{supp}(f)$ and for every $1\leq m\leq n$, there exists $\gamma_m\in\T$ such that $\tilde{f}\vert_{C_m}=\gamma_mf\vert_{C_m}$.
\item $f$ is phase retrievable w.r.t. $g$ if and only if $f$ is $L$-connected.
\end{enumerate}
\end{theorem}

\begin{proof} Let $\tilde{f}\in\C^d$.
\begin{enumerate}[a)]
\item Suppose that $\operatorname{supp}\left(\tilde{f}\right)=\operatorname{supp}(f)$ and for every $1\leq m\leq n$, there exists $\gamma_m\in\T$ such that $\tilde{f}\vert_{C_m}=\gamma_mf\vert_{C_m}$.
By Theorem \ref{thm:Cdgeneric}, it follows $\left|V_{\tilde{g}}f\right|=\left|V_{\tilde{g}}\tilde{f}\right|$ for every generic window $\tilde{g}\in\C_L^d$, which implies $V_ff\vert_{\Omega_L^d}=V_{\tilde{f}}\tilde{f}\vert_{\Omega_L^d}$, and therefore $\left|V_g f\right|=\left|V_g\tilde{f}\right|$, because of $\Omega(g)\subseteq\Omega_L^d$.\\
Conversely, suppose that $\left|V_gf\right|=\left|V_g\tilde{f}\right|$. By assumption, it holds $b_{j^{\ast}}^{(0)}=0$, which yields $\tilde{f}_{j^{\ast}}=\cdots=\tilde{f}_{j^{\ast}+L}=0$. Now, the assumption of Corollary \ref{cor:LinCombIt} is fulfilled with the choice of $j^{\ast}(k):=j^{\ast}$ for every $k$, which proves the claim.
\item This follows immediately from a).
\end{enumerate}
\end{proof}

Just as before (see again Remark \ref{rem:ZlengthL+1}), the proof has an algorithmic nature, which allows us to determine phase retrievability of the underlying signal for any given measurement.\\

Note that part a) of Theorem \ref{thm:CdHoleL+1} doesn't remain correct if we only ask for $L$ consecutive zeros: Choosing $r:=L+1\mid d$ in Example \ref{ex:cexprshift} shows that this won't even allow us to uniquely recover $\operatorname{supp}(f)$.
However, by adding a small assumption, we will be able to exclude periodic signals of this type, and thereby obtain a similar theorem for signals with only $L$ consecutive zeros. Moreover, note that the aforementioned counterexample is neither $L$-connected nor phase retrievable and thus doesn't affect part b) of the theorem.
In fact, this statement will remain true even with the weakened assumption.

We begin by classifying a certain kind of holes of length $d$ via the coefficients $b_j^{(k)}$.

\begin{lemma}\label{lm:CdHoleL}
Let $j^{\ast}\in\{0,\dots,d-1\}$. The following are equivalent.
\begin{enumerate}[(i)]
\item It holds $f_{j^{\ast}}\neq 0,f_j=0$ for all $j\in\{j^{\ast}+1,\dots,j^{\ast}+L\}$ and there exists at least one $j\in\{j^{\ast}-L,\dots,j^{\ast}-1\}$ such that $f_j\neq 0$.
\item The following conditions hold true.
\begin{enumerate}[a)]
\item It holds $b_j^{(k)}=0$ for all $1\leq k\leq L$ and all $j\in\{j^{\ast}+1-k,\dots,j^{\ast}+k\}$.
\item There exists $1\leq k\leq L$ such that $b_{j^{\ast}-k}^{(k)}\neq 0$.
\end{enumerate}
\end{enumerate}
\end{lemma}

\begin{proof}
Suppose that (i) holds true.
\begin{enumerate}[a)]
\item Fix $1\leq k\leq L$ and $j\in\{j^{\ast}+1-k,\dots,j^{\ast}+k\}$. We compute
\begin{equation*}
b_j^{(k)}=\sum\limits_{m=k}^L \overline{c_m^{(k)}}a_{m+j}^{(k)}=\sum\limits_{m=k}^L\overline{c_m^{(k)}}f_{j+m}\overline{f_{j+m-k}}.
\end{equation*}
Therefore, it suffices to show that for every $k\leq m\leq L$ either $f_{j+m}=0$ or $f_{j+m-k}=0$. First, we may estimate
\begin{equation*}
j^{\ast}+1\leq j+m\leq j^{\ast}+k+L.
\end{equation*}
Whenever $j+m\leq j^{\ast}+L$ holds true, this implies $f_{j+m}=0$. Otherwise, it holds\linebreak $j+m\geq j^{\ast}+L+1$ and thus
\begin{equation*}
j^{\ast}+1\leq j+m-k\leq j^{\ast}+L,
\end{equation*}
which yields $f_{j+m-k}=0$.
\item Since $f_{j^{\ast}+m-k}=0$ holds for every $k+1\leq m\leq L\,(\leq L+k)$, we obtain
\begin{align*}
b_{j^{\ast}-k}^{(k)}=\sum\limits_{m=k}^L\overline{c_m^{(k)}}f_{j^{\ast}+m-k}\overline{f_{j^{\ast}+m-2k}}=\overline{c_k^{(k)}}f_{j^{\ast}}\overline{f_{j^{\ast}-k}}.
\end{align*}
By assumption, this is non-zero for at least one $1\leq k\leq L$.
\end{enumerate}
\mbox{}\\
Now, suppose that (ii) holds true. Fix $1\leq k\leq L$ satisfying
\begin{equation*}
0\neq b_{j^{\ast}-k}^{(k)}=\sum\limits_{m=k}^L \overline{c_m^{(k)}}a_{j^{\ast}+m-k}^{(k)}.
\end{equation*}
This implies that there exists $k\leq m^{\ast}\leq L$ such that $a_{j^{\ast}+m^{\ast}-k}^{(k)}\neq 0$, i.e.
\begin{equation}\label{2entneq0}
f_{j^{\ast}+m^{\ast}-k}\neq 0\neq f_{j^{\ast}+m^{\ast}-2k}.
\end{equation}
First, we will show that $|M|\leq 1$, where
\begin{equation*}
M:=\left\{1\leq l\leq L~\middle|~f_{j^{\ast}+l}\neq 0\right\}.
\end{equation*}
We may as well show the stronger claim that $|\widetilde{M}|\leq 1$, where
\begin{equation*}
\widetilde{M}:=\left\{1\leq l\leq L+1~\middle|~f_{j^{\ast}+l}\neq 0\right\}.
\end{equation*}
In th case that $\widetilde{M}\neq\emptyset$, let $l_1:=\min\widetilde{M}$ and $l_2:=\max\widetilde{M}$ and assume that $1\leq n:=l_2-l_1\,(\leq L)$. It holds
\begin{equation*}
b_{j^{\ast}+1}^{(n)}=\sum\limits_{m=n}^L\overline{c_m^{(n)}}f_{j^{\ast}+m+1}\overline{f_{j^{\ast}+m+1-n}},
\end{equation*}
where $f_{j^{\ast}+m+1}=0$ whenever $m>l_2-1$ and $f_{j^{\ast}+m+1-n}=0$ whenever $m<l_1-1+n=l_2-1$. Therefore, we obtain
\begin{equation*}
0=b_{j^{\ast}+1}^{(n)}=\overline{c_{l_2-1}^{(n)}}f_{j^{\ast}+l_2}\overline{f_{j^{\ast}+l_1}}\neq 0,
\end{equation*}
which is a contradiction. Hence, it follows $l_1=l_2$, i.e. $|M|\leq\left|\widetilde{M}\right|\leq 1$.\\
Now, assume that $m^{\ast}> k$ holds in (\ref{2entneq0}). Since $|M|\leq 1$, it follows
\begin{equation*}
f_{j^{\ast}+l}=0\quad\text{for every } l\in\{1,\dots,L\}\setminus\{m^{\ast}-k\}.
\end{equation*}
This entails
\begin{equation*}
0=b_{j^{\ast}-k+1}^{(k)}=\sum\limits_{m=k}^L\overline{c_m^{(k)}}f_{j^{\ast}+m-k+1}\overline{f_{j^{\ast}+m-2k+1}}=\overline{c_{m^{\ast}-1}^{(k)}}f_{j^{\ast}+m^{\ast}-k}\overline{f_{j^{\ast}+m^{\ast}-2k}}\neq 0,
\end{equation*}
which is again a contradiction. Thus, it follows $m^{\ast}=k$, i.e. $f_{j^{\ast}}\neq 0\neq f_{j^{\ast}-k}$.\\
Finally, assume that there exists $1\leq l\leq L$ such that $f_{j^{\ast}+l}\neq 0$. The fact that $|M|\leq 1$ then implies that
\begin{equation*}
f_{j^{\ast}+m}=0 \quad\text{for every } m\in\{1,\dots,L\}\setminus\{l\}
\end{equation*}
and therefore
\begin{equation*}
0=b_{j^{\ast}}^{(l)}=\sum\limits_{m=l}^L\overline{c_m^{(l)}}f_{j^{\ast}+m}\overline{f_{j^{\ast}+m-l}}=\overline{c_l^{(l)}}f_{j^{\ast}+l}\overline{f_{j^{\ast}}}\neq 0,
\end{equation*}
which is yet another contradiction. Thus, we have shown (i).
\end{proof}

With Lemma \ref{lm:CdHoleL} at hand, we are able to prove a phase retrieval theorem for signals with $L$ consecutive zeros.

\begin{theorem}\label{thm:CdHoleL}
Let there be $j^{\ast}\in\{0,\dots,d-1\}$ such that $f_{j^{\ast}}\neq 0$ as well as
\begin{equation*}
\operatorname{supp}(f)\cap\{j^{\ast}+1,\dots,j^{\ast}+L\}=\emptyset,
\end{equation*}
and let $C_1,\dots,C_n$ be the $L$-connectivity components of $f$. Then, the following statements hold true.
\begin{enumerate}[a)]
\item If there exists at least one $j\in\{j^{\ast}-L,\dots,j^{\ast}-1\}$ such that $f_j\neq 0$, a signal $\tilde{f}\in\C^d$ satisfies $\left|V_g f\right|=\left|V_g\tilde{f}\right|$, if and only if $\operatorname{supp}\left(\tilde{f}\right)=\operatorname{supp}(f)$ and for every $1\leq m\leq n$, there exists $\gamma_m\in\T$ such that $\tilde{f}\vert_{C_m}=\gamma_mf\vert_{C_m}$.
\item $f$ is phase retrievable w.r.t. $g$ if and only if $f$ is $L$-connected.
\end{enumerate}
\end{theorem}

\begin{proof} Let $\tilde{f}\in\C^d$.
\begin{enumerate}[a)]
\item Suppose that there exists $j\in\{j^{\ast}-L,\dots,j^{\ast}-1\}$ such that $f_j\neq 0$.\\
When $\operatorname{supp}\left(\tilde{f}\right)=\operatorname{supp}(f)$ and for every $1\leq m\leq n$, there exists $\gamma_m\in\T$ such that $\tilde{f}\vert_{C_m}=\gamma_mf\vert_{C_m}$, it follows $\left|V_gf\right|=\left|V_g\tilde{f}\right|$ exactly as in Theorem \ref{thm:CdHoleL+1}.\\
Conversely, suppose $\left|V_gf\right|=\left|V_g\tilde{f}\right|$. By assumption, statement (i) from Lemma \ref{lm:CdHoleL} holds true. This yields statement (ii) from the lemma and therefore statement (i) applied to $\tilde{f}$. In particular, it follows $\tilde{f}_{j^{\ast}+l}=0$ for every $1\leq l\leq L$.
Therefore, we obtain
\begin{equation*}
a_{j^{\ast}+l}^{(k)}=0=\tilde{f}_{j^{\ast}+l}\overline{\tilde{f}_{j^{\ast}+l-k}}
\end{equation*}
for every $0\leq k\leq L$ and $1\leq l\leq L-k$, and Corollary \ref{cor:LinCombIt} implies $\operatorname{supp}\left(\tilde{f}\right)=\operatorname{supp}(f)$ as well as the existence of phases $\gamma_m\in\T$ such that $\tilde{f}\vert_{C_m}=\gamma_m f\vert_{C_m}$.
\item When $f$ is not $L$-connected, $f$ can't be phase retrievable w.r.t. a generic window $\tilde{g}\in\C_L^d$ by Theorem \ref{thm:Cdgeneric}. Because of Lemma \ref{lm:Vff}, $f$ can then neither be phase retrievable w.r.t. $g$.\\
Now, assume that $f$ is $L$-connected. If $f_{j^{\ast}+L+1}=0$, we may apply Theorem \ref{thm:CdHoleL+1} in order to show that $f$ is phase retrievable. Otherwise, since $f$ is $L$-connected, there exists\linebreak $j\in\{j^{\ast}-L,\dots,j^{\ast}-1\}$ such that $f_j\neq 0$. Now, a) shows that $f$ is phase retrievable.
\end{enumerate}
\end{proof}

Given a measurement, we are yet again able to decide whether the underlying signal is phase retrievable (even though the necessary calculations are somewhat more involved here).\\

We have already seen a signal with $L$ consecutive zeros, to which part a) of Theorem \ref{thm:CdHoleL} can't be applied.
However, if $f$ has no $L+1$ consecutive zeros and there exists one or several $j^{\ast}\in\{0,\dots,d-1\}$ satisfying $f_{j^{\ast}}\neq 0$ and $\operatorname{supp}(f)\cap\{j^{\ast}+1,\dots,j^{\ast}+L\}=\emptyset$, but none of them allows an index $j\in\{j^{\ast}-L,\dots,j^{\ast}-1\}$ such that $f_j\neq 0$, the characteristic function $\chi_{\operatorname{supp}(f)}$ is necessarily $L+1$-periodic with one period consisting of exactly one non-zero entry and $L$ zeros.
Since such signals only exist for special combinations of dimension and window length, we may drop the additional assumption of part a) in many instances.

\begin{corollary}
Let $L+1\nmid d$ and let there be $j^{\ast}\in\{0,\dots,d-1\}$, such that
\begin{equation*}
\operatorname{supp}(f)\cap\{j^{\ast}+1,\dots,j^{\ast}+L\}=\emptyset.
\end{equation*}
Then, a signal $\tilde{f}\in\C^d$ satisfies $\left|V_g f\right|=\left|V_g\tilde{f}\right|$, if and only if $\operatorname{supp}\left(\tilde{f}\right)=\operatorname{supp}(f)$ and for every $1\leq m\leq n$, there exists $\gamma_m\in\T$ such that $\tilde{f}\vert_{C_m}=\gamma_mf\vert_{C_m}$.
\end{corollary}

\begin{proof}
When $f\equiv 0$, the statement is trivial. Otherwise, the set
\begin{equation*}
I:=\left\{0\leq j^{\ast}\leq d-1~\middle|~f_{j^{\ast}}\neq 0,\, \operatorname{supp}(f)\cap\{j^{\ast}+1,\dots,j^{\ast}+L\}=\emptyset\right\}.
\end{equation*}
is non-empty. If the set
\begin{equation*}
\tilde{I}:=\left\{j^{\ast}\in I~|~f_{j^{\ast}+L+1}=0\right\}
\end{equation*}
is non-empty as well, we may apply Theorem \ref{thm:CdHoleL+1}.\\
Otherwise, suppose that $\operatorname{supp}(f)\cap\{j^{\ast}-L,\dots,j^{\ast}-1\}=\emptyset$ for every $j^{\ast}\in I$. However, this yields $j^{\ast}-(L+1)\in I$, since $f_{j^{\ast}-(L+1)}=0$ would imply $\tilde{I}\neq\emptyset$. By induction, $j^{\ast}-n(L+1)\in I$ holds for every $n\in\N$.
Since $L+1\nmid d$, there exists $n\in\N$ satisfying $d<n(L+1)<d+L+1$. The fact that $j^{\ast}-n(L+1)\in I$ then implies $0\neq f_{j^{\ast}-n(L+1)}=f_{j^{\ast}-(n(L+1)-d)}$ which is a contradiction since $1\leq n(L+1)-d\leq L$ and $\operatorname{supp}(f)\cap\{j^{\ast}-L,\dots,j^{\ast}-1\}=\emptyset$ by assumption.\\
Hence, there exists $j^{\ast}\in I$ such that $\operatorname{supp}(f)\cap\{j^{\ast}-L,\dots,j^{\ast}-1\}\neq\emptyset$ and we may apply Theorem \ref{thm:CdHoleL}.
\end{proof}

When compared to the literature, our main results of this section (Theorems \ref{thm:CdHoleL+1} and \ref{thm:CdHoleL}) appear to be unique in the sense that they deliberately accept $V_gg(0,l)$ to vanish for some $0\leq l\leq d-1$.\\

Next, one may wonder if we might be able to reduce the required number of consecutive zeros for $f$ even further. When choosing $r<L+1$ in Example \ref{ex:cexprshift}, we obtain $L$-connected, non-phase-retrievable signals with $r-1$ consecutive zeros. However, since $r\mid L+1$, it follows $r\leq\frac{L+1}{2}$ and thus $r-1<\frac{L}{2}$.
This justifies the conjecture that $\frac{L}{2}$ consecutive zeros might be enough to link phase retrievability to $L$-connectivity. Whether this conjecture holds true, remains an open question at this point.

\section{STFT phase retrieval on $\R^d$}\label{sec:Rd}

Finally, we will take a look at continuous phase retrieval on $G=\R^d$. The main parallel to the cyclic case is the fact that it holds $\widehat{G}\cong G$. Namely, the characters are given by
\begin{equation*}
\chi_{\omega}:\R^d\to\T,\quad t\mapsto e^{2\pi i\left\langle t,\omega\right\rangle}\quad\quad(\omega\in\R^d).
\end{equation*}

Therefore, we obtain the classical Fourier transform
\begin{equation*}
\widehat{h}(\omega)=\int_{\R^d} h(t)e^{-2\pi i\left\langle t,\omega\right\rangle}~\mathrm{d}t\quad (\omega\in\R^d)
\end{equation*}
and STFT
\begin{equation*}
V_gf(x,\omega)=\int_{\R^d} f(t)\overline{g(t-x)}e^{-2\pi i\left\langle t,\omega\right\rangle}~\mathrm{d}t\quad (x,\omega\in\R^d).
\end{equation*}
Even though we won't explicitly need the inversion formula, we note that it holds $c_G=1$ (cf. \cite{Groech}).

As suspected before, proving uniqueness results for the continuous STFT phase retrieval problem requires a lot more effort than in the previously discussed settings, due to the non-discreteness of the reconstruction.
We will therefore spread the proof of our main result over multiple subsections.

\subsection{Admissible translations and connectivity}

In this subsection, we want to identify a suitable window class for phase retrieval and prove some initial results that will help us at reconstructing a given signal.

Since we consider the case $G=\Z$ to be a blueprint for STFT phase retrieval, it is worthwhile to take a look back at the most general results for this setting. With Theorems \ref{thm:Zoneside} and \ref{thm:Zgeneral}, we have shown a complete characterization for phase retrieval with one-sided windows and windows of exponential decay.
The goal of this section is to prove analogous results for the continuous setting. (However, we restrict ourselves to dimension one when considering one-sided windows.)

More precisely, we will analyze phase retrieval for the following window classes.

\begin{definition}\label{def:expdecay}
 Let $p\geq 1$.
\begin{enumerate}[a)]
\item The class of all $L^p$-functions of exponential decay is given by
\begin{equation*}
L^{p,\operatorname{exp}}(\R^d):=\left\{h\in L^p(\R^d)~\middle|~\exists\sigma>0:~\int_{\R^d} |h(x)|^p\cdot e^{\sigma\cdot\|x\|_2}~\mathrm{d}x<\infty\right\}.
\end{equation*}
\item The class of all (one-dimensional) one-sided $L^p$-functions is given by
\begin{equation*}
L^{p,+}(\R):=\left\{h\in L^p(\R)~\middle|~h\cdot\chi_{(-\infty,0)}\equiv 0~\text{a.e.}\right\}.
\end{equation*}
\end{enumerate}
\end{definition}

\begin{remark}
Let $h:\R^d\to\C$ be measurable and suppose that there exist constants $C,\sigma>0$ satisfying
\begin{equation*}
|h(x)|\leq C\cdot e^{-\sigma\cdot\|x\|_2}
\end{equation*}
for almost every $x\in\R^d$. Then clearly, it holds $h\in L^{p,\exp}(\R^d)$ for every $p\geq 1$.
\end{remark}

It is easy to see that these properties carry over to shifts and products in the following way.

\begin{lemma}\label{lm:expdecaygTxg}\mbox{}
\begin{enumerate}[a)]
\item For every $g\in L^{2,\operatorname{exp}}(\R^d)$, it follows $g\cdot \overline{T_xg}\in L^{1,\operatorname{exp}}(\R^d)$.
\item For every $g\in L^{2,+}(\R)$, it follows $g\cdot \overline{T_xg}\in L^{1,+}(\R)$.
\end{enumerate}
\end{lemma}

\begin{proof}
Since the statement in b) is obvious, we only have to show a).\\
Let $\sigma>0$ be such that
\begin{equation*}
\int_{\R^d} |g(x)|^2\cdot e^{\sigma\cdot\|x\|_2}~\mathrm{d}x<\infty.
\end{equation*}
Hölder's inequality then implies
\begin{align*}
&~\int_{\R^d} |g(x)|\cdot|g(y-x)|\cdot e^{\sigma\cdot\|x\|_2}~\mathrm{d}x\\
=&~\int_{\R^d} |g(x)|\cdot e^{\nicefrac{\sigma}{2}\cdot\|x\|_2}\cdot |g(y-x)|\cdot e^{\nicefrac{\sigma}{2}\cdot\|x\|_2}~\mathrm{d}x\\
\leq&~ e^{\nicefrac{\sigma}{2}\cdot\|y\|_2}\cdot\int_{\R^d} |g(x)|\cdot e^{\nicefrac{\sigma}{2}\cdot\|x\|_2}\cdot |g(y-x)|\cdot e^{\nicefrac{\sigma}{2}\cdot\|x-y\|_2}~\mathrm{d}x\\
\leq&~ e^{\nicefrac{\sigma}{2}\cdot\|y\|_2}\cdot\left(\int_{\R^d} |g(x)|^2\cdot e^{\sigma\cdot\|x\|_2}~\mathrm{d}x\right)^{\nicefrac{1}{2}}\cdot\left(\int_{\R^d} |g(y-x)|^2\cdot e^{\sigma\cdot\|y-x\|_2}~\mathrm{d}x\right)^{\nicefrac{1}{2}}<\infty.
\end{align*}
\end{proof}

For further motivation, we will take a brief look at some known results:

When $d=1$, \cite{GJM} provides us (among others) with the following key examples of windows $g\in L^2(\R^d)$ whose ambiguity functions $V_gg$ do not vanish, and which therefore do phase retrieval by Theorem \ref{thm:Vggnonvan}:
\begin{enumerate}[1)]
\item The Gaussian functions defined by $\gamma_a(t):=e^{-a\pi t^2}$ for some $a>0$.
\item The functions $\eta_a$, defined by $\eta_a(t):=e^{-at}\chi_{(0,\infty)}(t)$ for some $a>0$.
\item The functions $c_a$, defined by $c_a(t):=\frac{1}{a+2\pi it}$ for some $a>0$.
\end{enumerate}
First, we observe that the functions $c_a$ are neither one-sided nor of exponential decay, despite doing phase retrieval. This showcases the fact that there are probably plenty of windows doing phase retrieval that lie outside of the class that is studied in this section.
Note however, that phase retrieval with windows $c_a\cdot\chi_{(0,\infty)}$ will be covered by our results below.\\
As for examples 1) and 2), we may draw a parallel to section \ref{sec:Z}, where we observed that phase retrieval seems to be less restrictive (at least proof-wise) for one-sided windows. (However, the results of this section will show that the ``two-sided'' versions of $\eta_a$ are also suitable for phase retrieval.)\\

The main goal of this subsection is to characterize the set $\Omega(g)$ for the relevant windows as well as to establish yet another notion of connectivity related to $\Omega(g)$.\\

While doing so, we will identify certain ``admissible'' translations, in order to ensure that phases can be recovered on every connectivity component of the signal.
For a given $x\in\R^d$, our first goal should be to recover $f\cdot \overline{T_x f}$ almost everywhere. By continuity of $V_ff$ and Lemma \ref{lm:fourinj}, this is possible whenever it holds $V_gg(x,\omega)\neq 0$ for almost every $\omega\in\R^d$.\\
It turns out that we may guarantee this for a large enough (and easily computable) set of translations, whenever the window $g$ is of exponential decay or one-sided (in dimension one). The proofs will again (and more heavily than in section \ref{sec:Z}) involve some results from complex analysis that we will state separately.\\
Note that complex analysis is a  common tool in the discussion of phase retrieval: Not only does it appear in \cite{AlaiBL}, where the authors follow a similar approach to ours, but also (among various other papers) in \cite{MW} for wavelet phase retrieval, and as early as in \cite{AkuI,AkuII} for classical Fourier phase retrieval.\\

We will begin with windows of exponential decay. This decay property has previously been used in the context of phase retrieval (cf. e.g. \cite{ADGT,JKP}), but mostly as a restriction for signals instead of STFT windows.

The first lemma gives a result on analytic continuation of the Fourier transform for functions of exponential decay. For dimension $1$, it can be found in \cite[Chapter 4, Theorem 3.1]{StSh}.
Moreover, its $L^2$-version is a well-known (cf. \cite{Katz} for dimension $1$) variant of the classical Paley-Wiener theorem for bandlimited or compactly supported functions and has already been used in the context of phase retrieval in \cite{ADGT}.\\
Our version requires no additional insights, but we will give the proof for convenience.

\begin{lemma}\label{lm:ancont}
Let $h\in L^{1,\operatorname{exp}}(\R^d)$ and $\sigma>0$ be such that
\begin{equation*}
\int_{\R^d}|h(x)|\cdot e^{\sigma\cdot\|x\|_2}~\mathrm{d}x<\infty.
\end{equation*}
Then, the function
\begin{equation*}
F:\mathcal{S}_{\sigma}\to\C,\quad \xi\mapsto\int_{\R^d} h(z)e^{-2\pi i\left\langle z,\xi\right\rangle}~\mathrm{d}z
\end{equation*}
is well-defined and holomorphic, where
\begin{equation*}
\mathcal{S}_{\sigma}:=\left\{\xi\in\C^d~\middle|~\left\|\operatorname{Im}(\xi)\right\|_2<\frac{\sigma}{2\pi}\right\}.
\end{equation*}
\end{lemma}

\begin{proof}
Note that for every $z\in\R^d$ and $\xi\in\mathcal{S}_{\sigma}$, it holds
\begin{align*}
\left|h(z)e^{-2\pi i\left\langle z,\xi\right\rangle}\right|	&= |h(z)|\cdot e^{\operatorname{Re}\left(-2\pi i\left\langle z,\xi\right\rangle\right)}=|h(z)|\cdot e^{-2\pi\left\langle z,\operatorname{Im}(\xi)\right\rangle}\\
									&\leq |h(z)|\cdot e^{2\pi\cdot\left\|\operatorname{Im}(\xi)\right\|_2\cdot\|z\|_2}\leq |h(z)|\cdot e^{\sigma\cdot\|z\|_2},
\end{align*}
which is integrable by assumption. Hence, $F$ is well-defined and continuous (by the dominated convergence theorem).\\
According to Hartogs's theorem, it suffices to check holomorphicity separately for each coordinate. Thus, consider $1\leq j\leq d$ and $\xi\in\mathcal{S}_{\sigma}$. For $s\in\C$, define $\xi^{(j,s)}\in\C^d$ by
\begin{equation*}
\xi_l^{(j,s)}:=\begin{cases}\xi_l,&\text{if } l\neq j,\\ s,&\text{if } l=j.\end{cases}
\end{equation*}
Clearly, $M:=\left\{s\in\C~\middle|~\xi^{(j,s)}\in\mathcal{S}_{\sigma}\right\}\subseteq\C$ is open and
\begin{equation*}
s\mapsto h(z)e^{-2\pi i\left\langle\xi^{(j,s)},z\right\rangle}
\end{equation*}
is holomorphic on $M$. For every closed curve $\gamma\subseteq M$, it follows
\begin{equation*}
\int_{\gamma} h(z)e^{-2\pi i\left\langle\xi^{(j,s)},z\right\rangle}~\mathrm{d}s=0.
\end{equation*}
An application of Fubini's theorem yields
\begin{equation*}
\int_{\gamma}\int_{\R^d} h(z)e^{-2\pi i\left\langle\xi^{(j,s)},z\right\rangle}~\mathrm{d}z\mathrm{d}s=\int_{\R^d}\int_{\gamma} h(z)e^{-2\pi i\left\langle\xi^{(j,s)},z\right\rangle}~\mathrm{d}s\mathrm{d}z=0,
\end{equation*}
which implies by Morera's theorem that $F$ is holomorphic in the $j$-th coordinate. Altogether, $F$ is holomorphic on $\mathcal{S}_{\sigma}$.
\end{proof}

The usage of this Paley-Wiener type theorem shows once again the connection of our approach to the one taken in \mbox{\cite{AlaiBL}}.

The second lemma takes a look at the zero sets of (multi-dimensional) holomorphic functions.

\begin{lemma}\label{lm:holomzeros}
Let $U\subseteq\C^d$ be open such that $\R^d\subseteq U$ and let $f:U\to\C$ be holomorphic. If $f$ doesn't vanish identically on $\R^d$, the set $\{z\in\R^d~|~f(z)=0\}\subseteq\R^d$ has measure zero.
\end{lemma}

\begin{proof}
We prove the lemma by induction, where the case $d=1$ is clear since a holomorphic function $f:\C\to\C$ that doesn't vanish identically (on $\R$), can only have countably many zeros.\\
Now suppose that the statement is true for a fixed dimension $d$ and consider an open subset $U\subseteq\C^{d+1}$ as well as a holomorphic function $f:U\to\C$ that doesn't vanish identically on $\R^{d+1}$. For every $x\in\R^d$, define
\begin{equation*}
\mathcal{N}_x:=\left\{y\in\R~\middle|~f(x,y)=0\right\}.
\end{equation*}
Whenever $\mathcal{N}_x$ has positive measure, it is in particular uncountable, and since
\begin{equation*}
U_x:=\left\{y\in\C~\middle|~(x,y)\in U\right\}\supseteq\R
\end{equation*}
is open and
\begin{equation*}
f_x:U_x\to\C,\quad z\mapsto f(x,z)
\end{equation*}
is holomorphic, it follows that $\mathcal{N}_x=\R$.\\
Now, let $\xi\in\R^{d+1}$ such that $f(\xi)\neq 0$ and define $y:=\xi_{d+1}$. Then, the function
\begin{equation*}
G_y:\widetilde{U}_y\to\C,\quad z\mapsto f(z,y),
\end{equation*}
where
\begin{equation*}
\widetilde{U}_y:=\left\{z\in\C^d~\middle|~(z,y)\in U\right\}\supseteq\R^d
\end{equation*}
is open, doesn't vanish identically on $\R^d$ and satisfies $G_y(x)=0$ whenever $y\in\mathcal{N}_x$. The induction hypothesis implies that
\begin{equation*}
\left\{x\in\R^d~\middle|~\leb\left(\mathcal{N}_x\right)>0\right\}=\left\{x\in\R^d~\middle|~\mathcal{N}_x=\R\right\}\subseteq\left\{x\in\R^d~\middle|~y\in\mathcal{N}_x\right\}
\end{equation*}
has measure zero. It follows
\begin{equation*}
\leb\left(\left\{z\in\R^{d+1}~\middle|~f(z)=0\right\}\right)=\int_{\R^d}\leb\left(\mathcal{N}_x\right)~\mathrm{d}x=0.
\end{equation*}
\end{proof}

Before we combine these results in order to establish a useful characterization of $\Omega(g)$, we will first turn to the case of a one-sided window in dimension one.

Just as before, our first step is a Paley-Wiener-type theorem. Once again, the $L^2$-version can be found in \cite{Katz} and the proof of the $L^1$-version is largely analogous to that of Lemma \ref{lm:ancont}.

\begin{lemma}\label{lm:ancont2}
Let $h\in L^{1,+}(\R)$ and define
\begin{equation*}
\mathbb{H}:=\left\{z\in\C~\middle|~\operatorname{Im}(z)>0\right\}.
\end{equation*}
Then, the function
\begin{equation*}
F:\overline{\mathbb{H}}\to\C,\quad \xi\mapsto \int_{\R}h(z)e^{-2\pi i\left\langle z,\xi\right\rangle}~\mathrm{d}z
\end{equation*}
is well-defined and continuous, as well as holomorphic on $\mathbb{H}$. Moreover, $F$ is bounded with $\left\|F\right\|_{\infty}\leq\|F\|_1$.
\end{lemma}

In order to understand the zero sets of the functions obtained in Lemma \ref{lm:ancont2}, we will (similarly as in section \ref{sec:Z}) apply the theory of Hardy spaces. This is also the reason why we only consider the case $d=1$.

\begin{lemma}\label{lm:zeroshalfplane}
Let $F:\overline{\mathbb{H}}\to\C$ be continuous and bounded, as well as holomorphic on $\mathbb{H}$. If $F$ doesn't vanish identically on $\mathbb{H}$, it holds $F(y)\neq 0$ for almost every $y\in\R$.
\end{lemma}

\begin{proof}
Let
\begin{equation*}
\widetilde{F}:K_1(0)\to\C,\quad z\mapsto F\left(-i\cdot\frac{1+z}{1-z}\right).
\end{equation*}
It is easy to see that $\widetilde{F}$ is well-defined. Furthermore, it holds
\begin{equation*}
\intpar{0}{2\pi}{\left|\widetilde{F}\left(re^{it}\right)\right|^2}{t}\leq 2\pi\cdot\|F\|_{\infty}^2
\end{equation*}
for every $0<r<1$, and consequently $\widetilde{F}$ is an element of the Hardy space $H^2$ (as defined in \mbox{\cite{Simon3}}). By \mbox{\cite[Theorem 5.2.6]{Simon3}}, it follows $\lim\limits_{r\uparrow 1} \widetilde{F}\left(rz\right)\neq 0$ for almost every $z\in\T$.\\
Finally, for every $y\in\R$, let $z_y:=\frac{iy+1}{iy-1}\in\T\setminus\{1\}$. Since the mapping $y\mapsto z_y$ is injective, we obtain
\begin{equation*}
0\neq\lim\limits_{r\uparrow 1} \widetilde{F}(rz_y)=\lim\limits_{r\uparrow 1} F\left(i\cdot\frac{1+rz_y}{1-rz_y}\right)=F\left(i\cdot\frac{1+z}{1-z}\right)=F(y)
\end{equation*}
for almost every $y\in\R$, where we applied continuity of $F$
\end{proof}

With these technical results at hand, we may now prove a Theorem concerning $\Omega(g)$ for windows $g$ from either of the two classes.

\begin{theorem}\label{thm:RdOmegag}
Let $g\in L^{2,\operatorname{exp}}(\R^d)$ or $d=1$ and $g\in L^{2,+}(\R)$. Then, the following statements hold true.
\begin{enumerate}[a)]
\item For every $x\in\R^d$ satisfying $g\cdot \overline{T_x g}\not\equiv 0$, it holds $V_gg(x,\omega)\neq 0$ for almost every $\omega\in\R^d$.
\item The set
\begin{equation*}
X:=\left\{x\in\R^d~\middle|~g\cdot \overline{T_x g}\not\equiv 0\right\}=\left\{x\in\R^d~\middle|~V_gg(x,\omega)\neq 0 \text{ for almost every } \omega\in\R^d\right\}
\end{equation*}
of all {\normalfont admissible translations} is an open subset of $\R^d$ satisfying $0\in X$.
\end{enumerate}
\end{theorem}

\begin{proof}\mbox{}
\begin{enumerate}[a)]
\item Let $x\in\R^d$ be such that $h_x:=g\cdot\overline{T_xg}\not\equiv 0$. By Lemma \ref{lm:expdecaygTxg}, it holds $h_x\in L^{1,\operatorname{exp}}(\R^d)$ or $h_x\in L^{1,+}(\R)$ respectively.
Combining Lemmas \ref{lm:ancont} and \ref{lm:holomzeros} (for exponential decay) or Lemmas \ref{lm:ancont2} and \ref{lm:zeroshalfplane} (for the one-sided case) yields
\begin{equation*}
\leb\left(\left\{\omega\in\R^d~\middle|~\widehat{h_x}(\omega)=0\right\}\right)=0,
\end{equation*}
i.e. $V_gg(x,\omega)\neq 0$ for almost every $\omega\in\R^d$.
\item Since $|g|^2\not\equiv 0$, it follows $0\in X$. Now, consider an arbitrary $x\in X$. Then there exists $\omega\in\R^d$ such that $V_gg(x,\omega)\neq 0$. Since $V_gg$ is continuous (in particular in the first component), there exists a neighbourhood $U$ of $x$ satisfying $V_gg\left(x^{\prime},\omega\right)\neq 0$ for all $x^{\prime}\in U$.
This obviously implies $g\cdot T_{x^{\prime}}\overline{g}\not\equiv 0$ and therefore $x^{\prime}\in X$ for all $x^{\prime}\in U$.
\end{enumerate}
\end{proof}

The set $X$ and most of its properties (at least for compactly supported $g$) have already been introduced in \cite{Jam}, where the author analyzes phase retrieval for the radar ambiguity function, i.e. the problem of reconstructing $f$ from $\left|V_ff\right|$.
Another parallel to \cite{Jam} manifests itself in the usage of Lebesgue's density theorem below.\\

As mentioned before, every admissible translation $x\in\R^d$ grants us full access to the function $f\cdot \overline{T_x f}$ -- but of course only in a $L^1$-sense. Thus, we have to be somewhat careful when defining a notion of connectivity for this setting.

First, we define $X$-connectivity for measurable sets.

\begin{definition} Consider a measurable subset $\mathcal{A}\subseteq\R^d$.
\begin{enumerate}[a)]
\item Two points $z,z^{\prime}\in\mathcal{A}$ are called $X$-connected and we write $z\sim_{\mathcal{A},X}z^{\prime}$, iff there exist\linebreak $z_0,\dots,z_n\in\mathcal{A}$ satisfying $z_0=z$, $z_n=z^{\prime}$ and $z_{m+1}-z_m\in X$ for every $0\leq m<n$.
\item Clearly, $\sim_{\mathcal{A},X}$ defines an equivalence relation on $\mathcal{A}$. The equivalence classes of $\mathcal{A}$ are called $X$-connectivity components of $\mathcal{A}$, and $\mathcal{A}$ is called $X$-connected, iff there exists at most one connectivity component.
\end{enumerate}
\end{definition}

For convenience, we will -- from now on and through the rest of the section -- consider the norm $\|\cdot\|_{\infty}$ on $\R^d$ instead of $\|\cdot\|_2$.

For any measurable set $M\subseteq\R^d$ and every $x\in\R^d$, we define
\begin{equation*}
d_M(x):=\lim\limits_{\varepsilon\to 0}\frac{\leb(M\cap U_{\varepsilon}(x))}{(2\varepsilon)^d}.
\end{equation*}
According to Lebesgue's density theorem, the set
\begin{equation*}
M^{\ast}:=\{x\in\R^d~|~d_M(x)=1\}
\end{equation*}
is measurable and satisfies $\leb(M\Delta M^{\ast})=0$, where
\begin{equation*}
M\Delta M^{\ast}:=\left(M\setminus M^{\ast}\right)\cup\left(M^{\ast}\setminus M\right)
\end{equation*}
denotes the symmetric difference of $M$ and $M^{\ast}$.

For every measurable function $f:\R^d\to\C$, let
\begin{equation*}
C(f):=\{z\in\R^d~|~f(z)\neq 0\}
\end{equation*}
as well as $C^{\ast}(f):=C(f)^{\ast}$.\\
We want to call $f$ connected, iff $C^{\ast}(f)$ is connected. The following lemma guarantees that this is well-defined even in an $L^2$-sense.

\begin{lemma}\label{lm:welldefconnect}
Let $f_1,f_2:\R^d\to\C$ be measurable such that $f_1\equiv f_2$ almost everywhere. Then, $C^{\ast}(f_1)=C^{\ast}(f_2)$.
\end{lemma}

\begin{proof}
Since $f_1$ and $f_2$ agree almost everywhere, it follows $\leb\left(C(f_1)\cap K\right)=\leb\left(C(f_2)\cap K\right)$ for every measurable subset $K\subseteq\R^d$ of finite measure. In particular, we obtain
\begin{equation*}
d_{C(f_1)}(x)=\lim\limits_{\varepsilon\to 0}\frac{\leb(C(f_1)\cap U_{\varepsilon}(x))}{(2\varepsilon)^d}=\lim\limits_{\varepsilon\to 0}\frac{\leb(C(f_2)\cap U_{\varepsilon}(x))}{(2\varepsilon)^d}=d_{C(f_2)}(x)
\end{equation*}
for every $x\in\R^d$, which yields $C^{\ast}(f_1)=C^{\ast}(f_2)$.
\end{proof}

Thus, we are finally able to make the following definition.

\begin{definition}
A signal $f\in L^2(\R^d)$ is called $X$-connected, iff $C^{\ast}(f)$ is $X$-connected. The $X$-connectivity components of $C^{\ast}(f)$ are called $X$-connectivity components of $f$.
\end{definition}

The use of $C^{\ast}(f)$ instead of $C(f)$ makes sure that connectivity is actually obtained in an $L^2$-sense. As always, it is not hard to prove connectivity to be a necessary condition for phase retrieval.

\begin{theorem}\label{thm:necconRd}
Let $g\in L^{2,\operatorname{exp}}(\R^d)$ or $d=1$ and $g\in L^{2,+}(\R)$. Furthermore, consider two signals $f,\tilde{f}\in L^2\left(\R^d\right)$ satisfying $C^{\ast}\left(\tilde{f}\right)=C^{\ast}(f)$ and let $\left(C_m\right)_{m\in I}$ be the $X$-connectivity components of $f$ (and thus also $\tilde{f}$). Then, the following statements hold true.
\begin{enumerate}[a)]
\item If for every $m\in\N$ there exists $\gamma_m\in\T$ such that $\tilde{f}\vert_{C_m}=\gamma_m f\vert_{C_m}$, it holds $\left|V_gf\right|=\left|V_g\tilde{f}\right|$.
\item If $f$ is not $X$-connected, then $f$ is not phase retrievable with respect to $g$.
\end{enumerate}
\end{theorem}

\begin{proof}\mbox{}
\begin{enumerate}[a)]
\item By Lebesgue's density theorem, it holds $f^{\prime}:=f\cdot\chi_{C^{\ast}(f)}+\chi_{C^{\ast}(f)\setminus C(f)}\equiv f$ almost everywhere. Therefore, Lemma \ref{lm:welldefconnect} implies $C^{\ast}\left(f^{\prime}\right)=C^{\ast}(f)=C\left(f^{\prime}\right)$. Thus, we may w.l.o.g. assume that $C^{\ast}(f)=C(f)$. The same applies for $\tilde{f}$.\\
Now, let $x\in X$ and $z\in\R^d$. Then, either
\begin{equation*}
\tilde{f}(z)\overline{\tilde{f}(z-x)}=f(z)\overline{f(z-x)}=0
\end{equation*}
or $z$ and $z-x$ belong to the same connectivity component $C_m$ of $C(f)=C\left(\tilde{f}\right)$. In this case, it holds
\begin{equation*}
\tilde{f}(z)\overline{\tilde{f}(z-x)}=\gamma_mf(z)\overline{\gamma_m}\overline{f(z-x)}=f(z)\overline{f(z-y)}.
\end{equation*}
Therefore, we obtain $V_ff\vert_{\Omega(g)}=V_{\tilde{f}}\tilde{f}\vert_{\Omega(g)}$, which implies $\left|V_gf\right|=\left|V_g\tilde{f}\right|$.
\item The statement follows immediately from a) as soon as we have shown that each of the connectivity components is measurable with positive measure. In order to prove this, let $m\in I$ and consider $z\in C_m$ as well as $L>0$ such that $U_L(0)\subseteq X$, which exists by Theorem \ref{thm:RdOmegag}.
The fact that $z\in C^{\ast}(f)$ implies
\begin{equation*}
\leb\left(C^{\ast}(f)\cap U_L(z)\right)=\leb\left(C(f)\cap U_L(z)\right)>0
\end{equation*}
and by choice of $L$, it follows $C^{\ast}(f)\cap U_L(z)\subseteq C_m$. Therefore, we obtain
\begin{equation*}
C_m=C^{\ast}(f)\cap\bigcup\limits_{z\in C_m} U_L(z),
\end{equation*}
which is measurable and of positive measure.
\end{enumerate}
\end{proof}

At least after taking the necessary precautions, the proof of Theorem \ref{thm:necconRd} appears to be largely analogous to its counterpart for the discrete setting. However, the opposite implication (i.e. signal recovery) has to be treated differently, since we are no longer able to propagate phases ``from point to point''.
Thus, it requires some effort to prove the following result.

\begin{theorem}\label{thm:Rdexpdec}
Let $g\in L^{2,\exp}(\R^d)$ or $d=1$ and $g\in L^{2,+}(\R)$. Furthermore, consider a signal $f\in L^{2}\left(\R^d\right)$ and let $\left(C_m\right)_{m\in I}$ be the $X$-connectivity components of $f$. Then, the following statements hold true.
\begin{enumerate}[a)]
\item A signal $\tilde{f}\in L^2(\R^d)$ satisfies $\left|V_g f\right|=\left|V_g\tilde{f}\right|$, if and only if $C^{\ast}\left(\tilde{f}\right)=C^{\ast}(f)$ and for every $m\in I$, there exists $\gamma_m\in\T$ such that $\tilde{f}\vert_{C_m}\equiv\gamma_mf\vert_{C_m}$ almost everywhere.
\item $f$ is phase retrievable w.r.t. $g$ if and only if $f$ is $X$-connected.
\end{enumerate}
\end{theorem}

The following two subsections will provide the proof of the implication of part a) that has not already been covered by Theorem \ref{thm:necconRd}. Once this has been shown, part b) follows immediately from a) and Theorem \ref{thm:necconRd}.

Before we turn to the proof, we will briefly emphasize the fact that, as in section \ref{sec:Z}, Lemma \ref{lm:windowshift} allows us to generalize (the one-sided version of) Theorem \ref{thm:Rdexpdec}.

\begin{corollary}\label{cor:Rdoneside}
Let $g\in L^2(\R)$ be such that $g\cdot\chi_I\equiv 0$ (almost everywhere) holds for some unbounded interval $I\subseteq\R$. Furthermore, consider a signal $f\in L^{2}\left(\R\right)$ and let $\left(C_m\right)_{m\in I}$ be the $X$-connectivity components of $f$. Then, the following statements hold true.
\begin{enumerate}[a)]
\item A signal $\tilde{f}\in L^2(\R^d)$ satisfies $\left|V_g f\right|=\left|V_g\tilde{f}\right|$, if and only if $C^{\ast}\left(\tilde{f}\right)=C^{\ast}(f)$ and for every $m\in I$, there exists $\gamma_m\in\T$ such that $\tilde{f}\vert_{C_m}\equiv\gamma_mf\vert_{C_m}$ almost everywhere.
\item $f$ is phase retrievable w.r.t. $g$ if and only if $f$ is $X$-connected.
\end{enumerate}
\end{corollary}

\subsection{Phase initialization}

Fix a window $g\in L^{2,\exp}\left(\R^d\right)$ or assume $d=1$ and fix a window $g\in L^{2,+}(\R)$. Additionally, fix two signals $f,\tilde{f}\in L^2\left(\R^d\right)$ satisfying $\left|V_gf\right|=\left|V_g\tilde{f}\right|$. By Theorem \ref{thm:VgfvsVff}, it follows $V_ff\vert_{\Omega(g)}=V_{\tilde{f}}\tilde{f}\vert_{\Omega(g)}$ as usual.
As in the proof of Theorem \ref{thm:necconRd}, we may assume w.l.o.g. that both $C^{\ast}(f)=C(f)$ and $C^{\ast}\left(\tilde{f}\right)=C\left(\tilde{f}\right)$. Note that, for the sake of clarity, we will only use this assumption within the proofs and not within the formulation of the lemmas along the way.\\

Because of the simple structure of $\Omega(g)$, it is easy to translate the measurement from Fourier domain to time domain.

\begin{lemma}\label{lm:Ax}
For every $x\in X$, there exists a set $\mathcal{N}_x\subseteq\R^d$ of measure zero, such that
\begin{equation*}
\mathcal{A}_x:=\mathcal{N}_x^c\subseteq\left\{z\in\R~\middle|~f(z)\overline{f(z-x)}=\tilde{f}(z)\overline{\tilde{f}(z-x)}\right\}.
\end{equation*}
\end{lemma}

\begin{proof}
Let $x\in X$. By definition, it holds $(x,\omega)\in\Omega(g)$ for almost every $\omega\in\R^d$. Thus,
\begin{equation*}
V_ff(x,\omega)=V_{\tilde{f}}\tilde{f}(x,\omega)\quad\quad\text{for almost every } \omega\in\R^d
\end{equation*}
and by continuity of $V_ff$ and $V_{\tilde{f}}\tilde{f}$ even for every $\omega\in\R^d$.\\
By Lemma \ref{lm:fourinj}, we obtain $f\cdot \overline{T_x f}=\tilde{f}\cdot \overline{T_x\tilde{f}}$ almost everywhere.
\end{proof}

\begin{remark}
Note that $V_ff(x,\cdot)\in L^1\left(\R^d\right)$ doesn't hold true in general. Therefore, we can't apply the Fourier inversion formula in the proof of Lemma \ref{lm:Ax}.
\end{remark}

In order to talk about (common) connectivity components of $f$ and $\tilde{f}$, we will first have to compare $C^{\ast}(f)$ and $C^{\ast}\left(\tilde{f}\right)$.

\begin{lemma}\label{lm:Castagree}
It holds $C^{\ast}(f)=C^{\ast}\left(\tilde{f}\right)$.
\end{lemma}

\begin{proof}
By Theorem \ref{thm:RdOmegag}, it holds $0\in X$. Thus, Lemma \ref{lm:Ax} implies $|f|^2=\left|\tilde{f}\right|^2$ almost everywhere. As a consequence, the functions $f^{\prime}:=\chi_{C(f)}$ and $\tilde{f}^{\prime}:=\chi_{C\left(\tilde{f}\right)}$ agree almost anywhere. Now, Lemma \ref{lm:welldefconnect} yields
\begin{equation*}
C^{\ast}(f)=C^{\ast}\left(f^{\prime}\right)=C^{\ast}\left(\tilde{f}^{\prime}\right)=C^{\ast}\left(\tilde{f}\right).
\end{equation*}
\end{proof}

With Lemma \ref{lm:Castagree} at hand, we know that the connectivity components $\left(C_m\right)_{m\in I}$ of $f$ are also the connectivity components of $\tilde{f}$.\\
Now, we fix one connectivity component $C_m$ in order to show that there exists $\gamma\in\T$ such that $\tilde{f}\vert_{C_m}\equiv \gamma f\vert_{C_m}$ almost everywhere.\\
Recall that we assumed w.l.o.g. that $C(f)=C^{\ast}(f)=C^{\ast}\left(\tilde{f}\right)=C\left(\tilde{f}\right)$, which implies that $f(z)\neq 0\neq\tilde{f}(z)$ holds for every $z\in C_m$.

The goal for the remainder of this subsection will be to fix $\tilde{f}=\gamma f$ on a measurable subset of $\mathcal{A}_0\cap C_m$ with positive measure.\\
First, we identify an area in which all translations are admissible. In order to do so, we fix $L>0$ such that $U_L(0)\subseteq X$, which exists since $X$ is open by Theorem \ref{thm:RdOmegag}. Note that we are still referring to a neighbourhood induced by $\|\cdot\|_{\infty}$.

\begin{lemma}\label{lm:I1I2}
There exist open, disjoint sets $I_1,I_2\subseteq\R^d$ satisfying
\begin{equation*}
\leb\left(I_j\cap\mathcal{A}_0\cap C_m\right)>0\quad\quad (j\in\{1,2\})
\end{equation*}
as well as $\operatorname{diam}\left(I_1\cup I_2\right)<L$.
\end{lemma}

\begin{proof}
As in the proof of Theorem \ref{thm:necconRd}, it holds $\leb\left(C_m\right)>0$. Since $\leb\left(\mathcal{A}_0^c\right)=0$, it follows
\begin{equation*}
\leb\left(\mathcal{A}_0\cap C_m\right)>0.
\end{equation*}
Now, consider the cubes
\begin{equation*}
\mathcal{W}_n:=\prod\limits_{j=1}^d \left(n_jL,(n_j+1)L\right)\quad\quad \left(n\in\Z^d\right).
\end{equation*}
As they form a countable covering of almost all of $\R^d$, we can pick $n\in\Z^d$ such that
\begin{equation*}
\leb\left(\mathcal{W}_n\cap\mathcal{A}_0\cap C_m\right)>0.
\end{equation*}
Using continuity of the mapping
\begin{equation*}
x\mapsto\leb\left(\left((n_1L,x)\times\prod\nolimits_{j=2}^d \left(n_jL,(n_j+1)L\right)\right)\cap\mathcal{A}_0\cap C_m\right),
\end{equation*}
we obtain the existence of $x_0\in (0,L)$ such that both
\begin{equation*}
I_1:=\left((n_1L,x_0)\times\prod\nolimits_{j=2}^d \left(n_jL,(n_j+1)L\right)\right)\cap\mathcal{A}_0\cap C_m
\end{equation*}
and
\begin{equation*}
I_2:=\left((x_0,(n_1+1)L)\times\prod\nolimits_{j=2}^d \left(n_jL,(n_j+1)L\right)\right)\cap\mathcal{A}_0\cap C_m
\end{equation*}
have positive measure.\\
The additional property $\operatorname{diam}\left(I_1\cup I_2\right)<L$ follows from $I_1\cup I_2\subseteq\mathcal{W}_n$.
\end{proof}

With the next lemma, we will ``narrow down'' $I_1$ and $I_2$ even further, to areas in which $\mathcal{A}_0\cap C_m$ shows up with ``high density''.

\begin{lemma}\label{lm:tastzast}\mbox{}
\begin{enumerate}[a)]
\item There exist $t^{\ast}\in I_1$, $z^{\ast}\in I_2$ and $\varepsilon>0$ satisfying $U_{\varepsilon}\left(t^{\ast}\right)\subseteq I_1$ and $U_{\varepsilon}\left(z^{\ast}\right)\subseteq I_2$ as well as
\begin{equation*}
\leb\left(\mathcal{A}_0\cap C_m\cap U_{\nicefrac{\varepsilon}{2}}\left(t^{\ast}\right)\right)\geq\frac{\varepsilon^d}{2}
\end{equation*}
and
\begin{equation*}
\leb\left(\mathcal{A}_0\cap C_m\cap U_{\nicefrac{3\varepsilon}{4}}\left(z^{\ast}\right)\right)\geq\left(1-\frac{1}{2\cdot 3^d}\right)\cdot\left(\frac{3}{2}\varepsilon\right)^d
\end{equation*}
\item Letting
\begin{equation*}
I^{\ast}:=U_{\nicefrac{\varepsilon}{4}}\left(z^{\ast}-t^{\ast}\right)\subseteq U_L(0),
\end{equation*}
it holds
\begin{equation*}
U_{\nicefrac{\varepsilon}{2}}\left(t^{\ast}\right)\subseteq U_{\nicefrac{3\varepsilon}{4}}\left(z^{\ast}\right)-x
\end{equation*}
for every $x\in I^{\ast}$.
\end{enumerate}
\end{lemma}

\begin{proof}\mbox{}
\begin{enumerate}[a)]
\item For $j\in\{1,2\}$, let
\begin{equation*}
J_j:=I_j\cap\mathcal{A}_0\cap C_m.
\end{equation*}
By Lemma \ref{lm:I1I2}, it holds $\leb(J_1),\leb(J_2)>0$, and by Lebesgue's density theorem, it follows that there exist $t^{\ast}\in J_1^{\ast}\cap J_1\subseteq I_1$ and $z^{\ast}\in J_2^{\ast}\cap J_2\subseteq I_2$. Now, we can choose $\tilde{\varepsilon}>0$ such that
\begin{equation*}
\frac{\leb\left(\mathcal{A}_0\cap C_m\cap U_{\nicefrac{\varepsilon}{2}}\left(t^{\ast}\right)\right)}{\varepsilon^d}\geq\frac{\leb\left(J_1\cap U_{\nicefrac{\varepsilon}{2}}\right)}{\varepsilon^d}\geq\frac{1}{2}
\end{equation*}
and
\begin{equation*}
\frac{\leb\left(\mathcal{A}_0\cap C_m\cap U_{\nicefrac{3\varepsilon}{4}}\left(z^{\ast}\right)\right)}{\left(\frac{3}{2}\varepsilon\right)^d}\geq\frac{\leb\left(J_2\cap U_{\nicefrac{3\varepsilon}{4}}\right)}{\left(\frac{3}{2}\varepsilon\right)^d}\geq 1-\frac{1}{2\cdot 3^d}
\end{equation*}
hold for all $0<\varepsilon<\tilde{\varepsilon}$. By openness of $I_1$ and $I_2$, we may then choose $0<\varepsilon<\tilde{\varepsilon}$ satisfying $U_{\varepsilon}\left(t^{\ast}\right)\subseteq I_1$ and $U_{\varepsilon}\left(z^{\ast}\right)\subseteq I_2$.
\item Let $x\in I^{\ast}$, i.e. $\left\|x-\left(z^{\ast}-t^{\ast}\right)\right\|<\frac{\varepsilon}{4}$, which implies $\|x\|<\frac{\varepsilon}{4}+\left\|z^{\ast}-t^{\ast}\right\|$. Next, consider
\begin{equation*}
\tilde{z}:=z^{\ast}+\frac{\varepsilon}{2\left\|z^{\ast}-t^{\ast}\right\|}\cdot\left(z^{\ast}-t^{\ast}\right).
\end{equation*}
Since $\left\|\tilde{z}-z^{\ast}\right\|=\frac{\varepsilon}{2}$, it follows $\tilde{z}\in I_2$ from a) and thus $\left\|\tilde{z}-t^{\ast}\right\|<L$ by Lemma \ref{lm:I1I2}. Because of
\begin{equation*}
\left\|\tilde{z}-t^{\ast}\right\|=\left\|\left(z^{\ast}-t^{\ast}\right)\cdot\left(1+\frac{\varepsilon}{2\left\|z^{\ast}-t^{\ast}\right\|}\right)\right\|=\left\|z^{\ast}-t^{\ast}\right\|+\frac{\varepsilon}{2},
\end{equation*}
this implies
\begin{equation*}
\|x\|<\frac{\varepsilon}{4}+\left\|z^{\ast}-t^{\ast}\right\|<L-\frac{\varepsilon}{4}<L.
\end{equation*}
Thus, $I^{\ast}\subseteq U_L(0)$.\\
Now, let $x\in I^{\ast}$ and $t\in U_{\nicefrac{\varepsilon}{2}}\left(t^{\ast}\right)$. It follows
\begin{equation*}
\left\|(x+t)-z^{\ast}\right\|\leq\left\|x-\left(z^{\ast}-t^{\ast}\right)\right\|+\left\|t-t^{\ast}\right\|<\frac{\varepsilon}{4}+\frac{\varepsilon}{2}=\frac{3\varepsilon}{4}
\end{equation*}
and as $t$ was picked arbitrarily, we obtain $U_{\nicefrac{\varepsilon}{2}}\left(t^{\ast}\right)\subseteq U_{\nicefrac{3\varepsilon}{4}}\left(z^{\ast}\right)-x$.
\end{enumerate}
\end{proof}

With these properties at hand, we are finally able to identify a point $t_0\in\mathcal{A}_0$ such that the set of all admissible translations $x\in X$ satisfying $t_0+x\in\mathcal{A}_0\cap\mathcal{A}_x$, is of positive measure. This will then allow us to initialize the phase on a set of positive measure.
Note that this property is not automatically guaranteed for any \textit{arbitrary} $t_0\in\mathcal{A}_0$, since it might very well hold $t_0+x\notin\mathcal{A}_x$ for every $x\in X$.

\begin{lemma}\label{lm:init0}
For every $t\in U_{\nicefrac{\varepsilon}{2}}\left(t^{\ast}\right)$, let
\begin{equation*}
M_t:=\left\{x\in I^{\ast}~\middle|~t\in\mathcal{A}_x-x\right\}.
\end{equation*}
Then, the following statements hold true.
\begin{enumerate}[a)]
\item For almost every $t\in U_{\nicefrac{\varepsilon}{2}}\left(t^{\ast}\right)$, the set $M_t$ is measurable and satisfies $\leb\left(M_t\right)=\left(\frac{\varepsilon}{2}\right)^d$.
\item There exists $t_0\in \mathcal{A}_0\cap C_m\cap U_{\nicefrac{\varepsilon}{2}}\left(t^{\ast}\right)$ such that
\begin{equation*}
\leb\left(C_m\cap\left(M_{t_0}+t_0\right)\right)>0.
\end{equation*}
\end{enumerate}
\end{lemma}

\begin{proof}
First, note that $\mathcal{A}_x$ is well-defined for every $x\in I^{\ast}\subseteq U_L(0)\subseteq X$ by Lemma \ref{lm:tastzast}b).
\begin{enumerate}[a)]
\item For every $x\in I^{\ast}$, we know that $\leb\left(\left(\mathcal{A}_x-x\right)^c\right)=\leb\left(\mathcal{A}_x^c-x\right)=0$ due to the translation invariance of the Lebesgue measure.\\
This implies $\leb\left(U_{\nicefrac{\varepsilon}{2}}\left(t^{\ast}\right)\cap\left(\mathcal{A}_x-x\right)\right)=\leb\left(U_{\nicefrac{\varepsilon}{2}}\left(t^{\ast}\right)\right)=\varepsilon^d$ and we can compute the double integral
\begin{align*}
\int_{I^{\ast}}\int_{U_{\nicefrac{\varepsilon}{2}}\left(t^{\ast}\right)}\chi_{\mathcal{A}_x-x}\left(t\right)~\mathrm{d}t\mathrm{d}x 	&=\int_{I^{\ast}}\leb\left(U_{\nicefrac{\varepsilon}{2}}\left(t^{\ast}\right)\cap\left(\mathcal{A}_x-x\right)\right)~\mathrm{d}x\\
																					&=\int_{I^{\ast}}\varepsilon^d~\mathrm{d}x=\frac{\varepsilon^{2d}}{2^d},
\end{align*}
using Lemma \ref{lm:tastzast}c). Fubini's theorem implies that
\begin{equation*}
\leb\left(M_t\right)=\int_{I^{\ast}}\chi_{\mathcal{A}_x-x}(t)~\mathrm{d}x
\end{equation*}
exists for almost every $t\in U_{\nicefrac{\varepsilon}{2}}\left(t^{\ast}\right)$. Furthermore, we may change the order of integration to obtain
\begin{equation*}
\int_{U_{\nicefrac{\varepsilon}{2}}\left(t^{\ast}\right)}\leb\left(M_t\right)~\mathrm{d}t =\int_{U_{\nicefrac{\varepsilon}{2}}\left(t^{\ast}\right)}\int_{I^{\ast}} \chi_{\mathcal{A}_x-x}(t)~\mathrm{d}t\mathrm{d}x=\int_{I^{\ast}}\int_{U_{\nicefrac{\varepsilon}{2}}\left(t^{\ast}\right)} \chi_{\mathcal{A}_x-x}(t)~\mathrm{d}t\mathrm{d}x=\frac{\varepsilon^{2d}}{2^d}.
\end{equation*}
Since $\leb\left(U_{\nicefrac{\varepsilon}{2}}\left(t^{\ast}\right)\right)=\varepsilon^d$ and $\leb\left(M_t\right)\leq\leb\left(I^{\ast}\right)=\left(\frac{\varepsilon}{2}\right)^d$ holds for every $t\in U_{\nicefrac{\varepsilon}{2}}\left(t^{\ast}\right)$, this implies $\leb\left(M_t\right)=\left(\frac{\varepsilon}{2}\right)^d$ for almost every $t\in U_{\nicefrac{\varepsilon}{2}}\left(t^{\ast}\right)$.
\item According to Lemma \ref{lm:tastzast}a), it holds
\begin{equation*}
\leb\left(\mathcal{A}_0\cap C_m\cap U_{\nicefrac{\varepsilon}{2}}\left(t^{\ast}\right)\right)>0.
\end{equation*}
Therefore, part a) implies the existence of $t_0\in\mathcal{A}_0\cap C_m\cap U_{\nicefrac{\varepsilon}{2}}\left(t^{\ast}\right)$ satisfying $\leb\left(M_{t_0}\right)=\left(\frac{\varepsilon}{2}\right)^d$.
Since $M_{t_0}\subseteq I^{\ast}$, it follows from Lemma \ref{lm:tastzast}b) that
\begin{equation*}
M_{t_0}+t_0\subseteq U_{\nicefrac{3\varepsilon}{4}}\left(z^{\ast}\right),
\end{equation*}
which implies
\begin{align*}
\leb\left(U_{\nicefrac{3\varepsilon}{4}}\left(z^{\ast}\right)\setminus\left(M_{t_0}+t_0\right)\right)&=\leb\left(U_{\nicefrac{3\varepsilon}{4}}\left(z^{\ast}\right)\right)-\leb\left(M_{t_0}+t_0\right)\\&=\left(\frac{3}{2}\varepsilon\right)^d-\left(\frac{\varepsilon}{2}\right)^d=\left(\frac{3}{2}\varepsilon\right)^d\cdot\left(1-\frac{1}{3^d}\right).
\end{align*}
Lemma \ref{lm:tastzast}a) allows us to compute
\begin{align*}
&~\leb\left(C_m\cap\left(M_{t_0}+t_0\right)\right)\\
=&~\leb\left(C_m\cap\left(M_{t_0}+t_0\right)\cap U_{\nicefrac{3\varepsilon}{4}}\left(z^{\ast}\right)\right)\\
=&~\leb\left(C_m\cap U_{\nicefrac{3\varepsilon}{4}}\left(z^{\ast}\right)\right)-\leb\left(C_m\cap\left(U_{\nicefrac{3\varepsilon}{4}}\left(z^{\ast}\right)\setminus\left(M_{t_0}+t_0\right)\right)\right)\\
\geq&~\leb\left(\mathcal{A}_0\cap C_m\cap U_{\nicefrac{3\varepsilon}{4}}\left(z^{\ast}\right)\right)-\leb\left(U_{\nicefrac{3\varepsilon}{4}}\left(z^{\ast}\right)\setminus\left(M_{t_0}+t_0\right)\right)\\
\geq&~\left(1-\frac{1}{2\cdot 3^d}\right)\cdot\left(\frac{3}{2}\varepsilon\right)^d-\left(1-\frac{1}{3^d}\right)\cdot\left(\frac{3}{2}\varepsilon\right)^d\\
=&~\frac{1}{2\cdot 3^d}\cdot\left(\frac{3}{2}\varepsilon\right)^d>0.
\end{align*}
\end{enumerate}
\end{proof}

We will now propagate the phase from $t_0$ onto $M_{t_0}+t_0$ to obtain the following theorem.

\begin{theorem}\label{thm:phaseiniRd}
There exist $\gamma\in\T$ and a measurable subset $\mathcal{B}_0\subseteq\R^d$ satisfying $\leb\left(\mathcal{B}_0\cap C_m\right)>0$ as well as
\begin{equation*}
\tilde{f}(z)=\gamma\cdot f(z)\quad\text{for all } z\in\mathcal{B}_0.
\end{equation*}
\end{theorem}

\begin{proof}
As indicated before, we will choose $\mathcal{B}_0:=M_{t_0}+t_0$ with $t_0$ from Lemma \ref{lm:init0}b). By Lemma \ref{lm:init0}b), it holds $\leb\left(\mathcal{B}_0\cap C_m\right)>0$.\\
Since $t_0\in\mathcal{A}_0\cap C_m$, i.e. $\left|f(t_0)\right|=\left|\tilde{f}(t_0)\right|\neq 0$, we may define
\begin{equation*}
\gamma:=\frac{\tilde{f}(t_0)}{f(t_0)}\in\T.
\end{equation*}
Now, for every $z\in\mathcal{B}_0$, it holds $x:=z-t_0\in M_{t_0}$, i.e. $x\in I^{\ast}\subseteq X$ and $z=t_0+x\in\mathcal{A}_x$. This implies
\begin{equation*}
\tilde{f}(z)=f(z)\cdot\frac{\overline{f(t_0)}}{\overline{\tilde{f}(t_0)}}=\gamma f(z).
\end{equation*}
\end{proof}

Note that from now on, we will have no more need of $\mathcal{A}_0$, since the condition $\tilde{f}(z)=\gamma f(z)$ already imlies $|f(z)|=\left|\tilde{f}(z)\right|$.

\subsection{Phase propagation}

In this subsection, we will start things off from Theorem \ref{thm:phaseiniRd} and propagate the phase $\gamma$ onto almost all of $C_m$. We will still proceed iteratively, along the sets $\left(\mathcal{D}_k\right)_{k\in\N}$, defined by
\begin{equation*}
\mathcal{D}_0:=\left(\mathcal{B}_0\cap C_m\right)^{\ast}
\end{equation*}
as well as
\begin{equation*}
\mathcal{D}_{k+1}:=\left(\mathcal{D}_k+X\right)\cap C_m\quad\quad(k\in\N_0).
\end{equation*}
The most important properties are summarized in the following lemma.

\begin{lemma}\label{lm:propDk}
The following statements hold true.
\begin{enumerate}[a)]
\item It holds $\emptyset\neq\mathcal{D}_0\subseteq C_m$.
\item For every $k\in\N_0$, it holds $\mathcal{D}_k\subseteq\mathcal{D}_k^{\ast}$.
\item It holds $C_m=\bigcup_{k\in\N_0} \mathcal{D}_k$.
\end{enumerate}
\end{lemma}

\begin{proof}\mbox{}
\begin{enumerate}[a)]
\item $\mathcal{D}_0\neq\emptyset$ follows from the fact that $\leb\left(\mathcal{D}_0\right)=\leb\left(\mathcal{B}_0\cap C_m\right)>0$ holds by Theorem \ref{thm:phaseiniRd}.\\
Moreover, it holds
\begin{equation*}
\mathcal{D}_0=\left(\mathcal{B}_0\cap C_m\right)^{\ast}\subseteq C_m^{\ast}\subseteq C^{\ast}(f).
\end{equation*}
For every $z\in\mathcal{D}_0$, the fact that $z\in C_m^{\ast}$ implies that $C_m\cap U_L(z)\neq\emptyset$. Hence, there exists $y\in C_m$ such that $x:=z-y\in U_L(0)\subseteq X$. Together with $z\in C^{\ast}(f)$, this yields $z\in C_m$.
\item We will prove this statement by induction. For $k=0$, note that every subset $M\subseteq\R^d$ satisfies
\begin{equation*}
d_{M^{\ast}}(z)=\lim\limits_{\varepsilon\to 0} \frac{\leb\left(M^{\ast}\cap U_{\varepsilon}(z)\right)}{(2\varepsilon)^d}=\lim\limits_{\varepsilon\to 0} \frac{\leb\left(M\cap U_{\varepsilon}(z)\right)}{(2\varepsilon)^d}=d_M(z)
\end{equation*}
for every $z\in\R^d$ by Lebesgue's density theorem. In particular, $\mathcal{D}_0^{\ast}=\mathcal{D}_0$.\\
Now, assume that $\mathcal{D}_k\subseteq\mathcal{D}_k^{\ast}$ holds true for some $k\in\N_0$ and pick $z\in\mathcal{D}_{k+1}$. Moreover, let $0<\alpha<1$.\\
First, $z\in\mathcal{D}_{k+1}$ implies the existence of $x\in X$ such that
\begin{equation*}
z-x\in\mathcal{D}_k\subseteq\mathcal{D}_k^{\ast}
\end{equation*}
by the induction hypothesis. Therefore, there exists $\varepsilon_1>0$ such that
\begin{equation*}
\frac{\leb\left(\mathcal{D}_k\cap U_{\varepsilon}(z-x)\right)}{(2\varepsilon)^d}\geq 1-\frac{\alpha}{2}
\end{equation*}
holds for every $0<\varepsilon<\varepsilon_1$. Consequently, we can estimate
\begin{equation*}
\frac{\leb\left(\left(\mathcal{D}_k+X\right)\cap U_{\varepsilon}(z)\right)}{(2\varepsilon)^d}\geq\frac{\leb\left(\left(\mathcal{D}_k+x\right)\cap U_{\varepsilon}(z)\right)}{(2\varepsilon)^d}=\frac{\leb\left(\mathcal{D}_k\cap U_{\varepsilon}(z-x)\right)}{(2\varepsilon)^d}\geq 1-\frac{\alpha}{2}.
\end{equation*}
On the other hand, it holds $z\in C_m\subseteq C^{\ast}(f)$ as well as $C_m\cap U_{\varepsilon}(z)=C^{\ast}(f)\cap U_{\varepsilon}(z)$ for every $0<\varepsilon<L$, since every point of $C^{\ast}(f)\cap U_{\varepsilon}(z)$ belongs to $C_m$ because of $U_L(0)\subseteq X$. Hence, we may choose $0<\varepsilon_2<L$ such that
\begin{equation*}
\frac{\leb\left(C_m\cap U_{\varepsilon}(z)\right)}{(2\varepsilon)^d}=\frac{\leb\left(C^{\ast}(f)\cap U_{\varepsilon}(z)\right)}{(2\varepsilon)^d}=\frac{\leb\left(C(f)\cap U_{\varepsilon}(z)\right)}{(2\varepsilon)^d}\geq 1-\frac{\alpha}{2}
\end{equation*}
holds for every $0<\varepsilon<\varepsilon_2$. Together, it follows
\begin{equation*}
\frac{\leb\left(\mathcal{D}_{k+1}\cap U_{\varepsilon}(z)\right)}{(2\varepsilon)^d}\geq 1-\alpha
\end{equation*}
for every $0<\varepsilon<\min\left\{\varepsilon_1,\varepsilon_2\right\}$. Since $0<\alpha<1$ was picked arbitrarily, we obtain $z\in\mathcal{D}_{k+1}^{\ast}$.
\item The inclusion ``$\supseteq$'' is clear, since $\mathcal{D}_k\subseteq C_m$ holds by definition when $k\in\N$ and by part a) when $k=0$.\\
Fix $z\in\mathcal{D}_0$ which exists by part a). Additionally, it follows $z\in C_m$.\\
Now, let $z^{\prime}\in C_m$ be arbitrary. Then, there exist $z_1,\dots, z_n\in C_m$ satisfying
\begin{equation*}
x_j:=z_j-z_{j-1}\in X\quad\quad\text{for every } 1\leq j\leq n+1,
\end{equation*}
where we let $z_0:=z$ and $z_{n+1}:=z^{\prime}$. Inductively, this implies $z^{\prime}\in\mathcal{D}_{n+1}\subseteq\bigcup_{k\in\N_0} \mathcal{D}_k$.
\end{enumerate}
\end{proof}

Recall that Theorem \ref{thm:phaseiniRd} states the existence of $\gamma\in\T$ satisfying $\tilde{f}(z)=\gamma\cdot f(z)$ for all $z\in\mathcal{B}_0\cap C_m\subseteq\mathcal{B}_0$.
Since $\leb\left(\mathcal{D}_0\right)=\leb\left(\mathcal{B}_0\cap C_m\right)$ holds by Lebesgue's density theorem, we obtain that $\tilde{f}(z)=\gamma\cdot f(z)$ holds for almost every $z\in\mathcal{D}_0$.\\
As mentioned before, our main goal will now be to propagate the phase in each step from (most of) $\mathcal{D}_k$ to (most of) $\mathcal{D}_{k+1}$. We will first prove a local version of this iteration.

\begin{lemma}\label{lm:locpropRd}
Let $k\in\N_0$ and assume that there exists a subset $\mathcal{E}_k\subseteq\mathcal{D}_k$ satisfying\linebreak $\leb\left(\mathcal{D}_k\setminus\mathcal{E}_k\right)=0$ as well as
\begin{equation*}
\tilde{f}(t)=\gamma\cdot f(t)\quad\quad\text{for every } t\in\mathcal{E}_k.
\end{equation*}
Then, for every $z\in\mathcal{D}_{k+1}$, there exists an open neighbourhood $U_z$ of $z$ such that $\tilde{f}\left(z^{\prime}\right)=\gamma\cdot f\left(z^{\prime}\right)$ holds for almost every $z^{\prime}\in U_z$.
\end{lemma}

\begin{proof}
Let $z\in\mathcal{D}_{k+1}$. By definition, there exists $y\in\mathcal{D}_k$ such that $x:=z-y\in X$. Furthermore, since $X$ is open by Theorem \ref{thm:RdOmegag}, we can choose $\varepsilon>0$ satisfying $U_{\varepsilon}(x)\subseteq X$.\\
Now, pick $0<\alpha<1$. We will show that
\begin{equation*}
\leb\left(\left\{z^{\prime}\in U_{\varepsilon}(z)~\middle|~\tilde{f}\left(z^{\prime}\right)\neq\gamma f\left(z^{\prime}\right)\right\}\right)\leq (1-\alpha)\cdot (2\varepsilon)^d.
\end{equation*}
By Lemma \ref{lm:propDk}, it holds $y\in\mathcal{D}_{k}^{\ast}$. Hence, we may choose $\delta>0$ such that
\begin{equation*}
\leb\left(\mathcal{D}_k\cap U_{\nicefrac{\delta}{2}}(y)\right)\geq \alpha\cdot\delta^d
\end{equation*}
and assume w.l.o.g. that $n:=\frac{\varepsilon}{\delta}-\frac{1}{2}\in\N$.\\
For all $j\in\Z^d$ satisfying $\|j\|=\|j\|_{\infty}\leq n$, this implies
\begin{equation*}
\|\delta j\|\leq\varepsilon-\frac{\delta}{2}<\varepsilon
\end{equation*}
and therefore $x+\delta j\in U_{\varepsilon}(x)\subseteq X$. Additionally, we obtain
\begin{equation*}
\bigcup\limits_{\|j\|\leq n} U_{\nicefrac{\delta}{2}}(y)+(x+\delta j)\subseteq U_{\varepsilon}(z)=:U_z,
\end{equation*}
where the union on the left-hand side is disjoint. As a consequence,
\begin{equation*}
M:=\bigcup\limits_{\|j\|\leq n}\left(\left(\mathcal{E}_k\cap U_{\nicefrac{\delta}{2}}(y)\right)+(x+\delta j)\right)\cap\mathcal{A}_{x+\delta j}
\end{equation*}
is also a subset of $U_z$ and since this union is still disjoint, we may compute
\begin{align*}
\leb(M)	&=\sum\limits_{\|j\|\leq n} \leb\left(\left(\left(\mathcal{E}_k\cap U_{\nicefrac{\delta}{2}}(y)\right)+(x+\delta j)\right)\cap\mathcal{A}_{x+\delta j}\right)\\
		&=\sum\limits_{\|j\|\leq n} \leb\left(\left(\mathcal{E}_k\cap U_{\nicefrac{\delta}{2}}(y)\right)+(x+\delta j)\right)\\
		&=\sum\limits_{\|j\|\leq n} \leb\left(\mathcal{E}_k\cap U_{\nicefrac{\delta}{2}}(y)\right)\\
		&=\sum\limits_{\|j\|\leq n} \leb\left(\mathcal{D}_k\cap U_{\nicefrac{\delta}{2}}(y)\right)\\
		&\geq\sum\limits_{\|j\|\leq n} \alpha\cdot\delta^d=\alpha\cdot(2n+1)^d\delta^d=\alpha\cdot(2\varepsilon)^d.
\end{align*}
Finally, pick $z^{\prime}\in M$. By definition, there exists $j\in\Z^d$ satisfying $\|j\|\leq n$ and $t\in\mathcal{E}_k$ such that $z^{\prime}\in\mathcal{A}_{x+\delta j}$ as well as $z^{\prime}=t+(x+\delta j)$, i.e.
\begin{equation*}
f\left(z^{\prime}\right)\overline{f(t)}=\tilde{f}\left(z^{\prime}\right)\overline{\tilde{f}(t)}.
\end{equation*}
The fact that $\mathcal{E}_k\subseteq\mathcal{D}_k\subseteq C_m\subseteq C(f)$ implies
\begin{equation*}
\tilde{f}(t)=\gamma\cdot f(t)\neq 0.
\end{equation*}
Thus, we obtain $\tilde{f}\left(z^{\prime}\right)=\gamma\cdot f\left(z^{\prime}\right)$ and are able to conclude that
\begin{equation*}
\leb\left(\left\{z^{\prime}\in U_{\varepsilon}(z)~\middle|~\tilde{f}\left(z^{\prime}\right)\neq\gamma\cdot f\left(z^{\prime}\right)\right\}\right)\leq (1-\alpha)\cdot (2\varepsilon)^d.
\end{equation*}
As $0<\alpha<1$ was picked arbitrarily, it follows $\tilde{f}\left(z^{\prime}\right)=\gamma\cdot f\left(z^{\prime}\right)$ for almost every $z^{\prime}\in U_z$.
\end{proof}

In order to obtain global propagation from (most of) $\mathcal{D}_k$ onto (most of) $\mathcal{D}_k$, we have to ensure that the sets
\begin{equation*}
\left\{z^{\prime}\in U_z~\middle|~\tilde{f}\left(z^{\prime}\right)\neq\gamma\cdot f\left(z^{\prime}\right)\right\},
\end{equation*}
who are individually of measure zero for every single $z\in\mathcal{D}_k$, don't add up towards
\begin{equation*}
\left\{z^{\prime}\in\mathcal{D}_{k+1}~\middle|~\tilde{f}\left(z^{\prime}\right)\neq\gamma\cdot f\left(z^{\prime}\right)\right\}
\end{equation*}
having positive measure. In order to do so, we introduce some kind of compactness.

\begin{lemma}\label{lm:globpropRd}
Let $k\in\N_0$ and assume that there exists a subset $\mathcal{E}_k\subseteq\mathcal{D}_k$ satisfying\linebreak $\leb\left(\mathcal{D}_k\setminus\mathcal{E}_k\right)=0$ as well as
\begin{equation*}
\tilde{f}(t)=\gamma\cdot f(t)\quad\quad\text{for every } t\in\mathcal{E}_k.
\end{equation*}
Then, there exists a subset $\mathcal{E}_{k+1}\subseteq\mathcal{D}_{k+1}$ satisfying $\leb\left(\mathcal{D}_{k+1}\setminus\mathcal{E}_{k+1}\right)=0$ as well as $\tilde{f}(t)=\gamma\cdot f(t)$ for every $t\in\mathcal{E}_{k+1}$.
\end{lemma}

\begin{proof}
By Lemma \ref{lm:locpropRd}, for every $z\in\mathcal{D}_{k+1}$, there exists an open neighbourhood $U_z$ of $z$ as well as a subset $\mathcal{E}_z\subseteq U_z$ satisfying $\leb\left(U_z\setminus\mathcal{E}_z\right)=0$ and
\begin{equation*}
\tilde{f}\left(z^{\prime}\right)=\gamma\cdot f\left(z^{\prime}\right)\quad\quad\text{for every } z^{\prime}\in\mathcal{E}_z.
\end{equation*}
Since $\mathcal{D}_{k+1}$ is measurable (by induction), there exists a sequence $\left(A_n\right)_{n\in\N}$ of closed subsets of $\mathcal{D}_{k+1}$, such that
\begin{equation*}
\leb\left(\mathcal{D}_{k+1}\setminus A_n\right)<\frac{1}{n}
\end{equation*}
holds for every $n\in\N$. We may assume w.l.o.g. that $A_n\subseteq A_{n+1}$ holds for every $n\in\N$. (Otherwise, we consider $\tilde{A}_n:=\bigcup_{j\leq n} A_j$ instead.)\\
Now, let $K_n:=A_n\cap [-n,n]^d$. Clearly, $K_n$ is compact for every $n\in\N$ and it holds\linebreak $\bigcup_{n\in\N} K_n=\bigcup_{n\in\N} A_n$. Moreover, we obtain
\begin{equation*}
\leb\left(\mathcal{D}_{k+1}\setminus\bigcup\limits_{n\in\N} K_n\right)=\leb\left(\mathcal{D}_{k+1}\setminus\bigcup\limits_{n\in\N} A_n\right)\leq\leb\left(\mathcal{D}_{k+1}\setminus A_j\right)<\frac{1}{j}
\end{equation*}
for every $j\in\N$. This implies
\begin{equation*}
\leb\left(\mathcal{D}_{k+1}\setminus\bigcup\limits_{n\in\N} K_n\right)=0.
\end{equation*}
Let $n\in\N$. Together with $K_n\subseteq\bigcup_{z\in K_n} U_z$, the compactness of $K_n$ implies the existence of a finite set $J_n\subseteq K_n$ satisfying
\begin{equation*}
K_n\subseteq\bigcup\limits_{z\in J_n} U_z.
\end{equation*}
Letting
\begin{equation*}
\mathcal{E}^{(n)}:=\bigcup\limits_{z\in J_n}\mathcal{E}_z,
\end{equation*}
we obtain $\tilde{f}\left(z^{\prime}\right)=\gamma\cdot f\left(z^{\prime}\right)$ for every $z^{\prime}\in\mathcal{E}^{(n)}$ as well as
\begin{equation*}
\leb\left(K_n\setminus\mathcal{E}^{(n)}\right)\leq\leb\left(\bigcup\limits_{z\in J_n} U_z\setminus\mathcal{E}^{(n)}\right)\leq\sum\limits_{z\in J_n}\leb\left(U_z\setminus\mathcal{E}^{(n)}\right)\leq\sum\limits_{z\in J_n}\leb\left(U_z\setminus\mathcal{E}_z\right)=0.
\end{equation*}
Finally, define
\begin{equation*}
\mathcal{E}_{k+1}:=\left(\bigcup\limits_{n\in\N}\mathcal{E}^{(n)}\right)\cap\mathcal{D}_{k+1}.
\end{equation*}
Clearly, it holds $\tilde{f}\left(z^{\prime}\right)=\gamma\cdot f\left(z^{\prime}\right)$ for every $z^{\prime}\in\mathcal{E}_{k+1}$. Furthermore, it follows
\begin{equation*}
\leb\left(\mathcal{D}_{k+1}\setminus\mathcal{E}_{k+1}\right)=\leb\left(\mathcal{D}_{k+1}\setminus\bigcup\limits_{n\in\N} K_n\right)+\leb\left(\bigcup\limits_{n\in\N} K_n\setminus\mathcal{E}_{k+1}\right)\leq\sum\limits_{n\in\N} \leb\left(K_n\setminus\mathcal{E}^{(n)}\right)=0,
\end{equation*}
which proves the claim.
\end{proof}

Recall that $\tilde{f}(z)=\gamma\cdot f(z)$ holds for every $z\in\mathcal{B}_0\cap C_m$ by Theorem \ref{thm:phaseiniRd}.\\
Since $\leb\left(\mathcal{D}_0\setminus\left(\mathcal{B}_0\cap C_m\right)\right)=0$ holds by Lebesgue's density theorem, this implies that there is a measurable subset $\mathcal{E}_0\subseteq\mathcal{D}_0$ satisfying $\leb\left(\mathcal{D}_0\setminus\mathcal{E}_0\right)=0$ as well as $\tilde{f}(z)=\gamma\cdot f(z)$ for every $z\in\mathcal{E}_0$.
For every $k\in\N$, Lemma \ref{lm:globpropRd} inductively implies the existence of a measurable set $\mathcal{E}_k\subseteq\mathcal{D}_k$ satisfying $\leb\left(\mathcal{D}_k\setminus\mathcal{E}_k\right)=0$ as well as $\tilde{f}(z)=\gamma\cdot f(z)$ for every $z\in\mathcal{E}_k$.\\
Finally, define
\begin{equation*}
\mathcal{E}:=\bigcup\limits_{k\in\N_0}\mathcal{E}_k.
\end{equation*}
Clearly, $\tilde{f}(z)=\gamma\cdot f(z)$ holds for every $z\in\mathcal{E}$. Additionally, Lemma \ref{lm:propDk}c) yields
\begin{equation*}
\leb\left(C_m\setminus\mathcal{E}\right)=\leb\left(\bigcup\limits_{k\in\N_0}\mathcal{D}_k\setminus\mathcal{E}\right)\leq\sum\limits_{k\in\N_0}\leb\left(\mathcal{D}_k\setminus\mathcal{E}\right)\leq\sum\limits_{k\in\N_0}\leb\left(\mathcal{D}_k\setminus\mathcal{E}_k\right)=0,
\end{equation*}
which concludes the proof of Theorem \ref{thm:Rdexpdec}.

\subsection{Examples}

We conclude this section by applying Theorem \ref{thm:Rdexpdec} to a few interesting windows. First however, we make a general remark towards understanding the result.

\begin{remark}
Even though Theorem \ref{thm:Rdexpdec} is easy to understand from an abstract point of view, it can sometimes be hard to decide whether a given signal is $X$-connected. Consider e.g. a closed, symmetric set $A\subseteq\R^d$ of positive measure, satisfying $0\notin A$. For $X:=\R^d\setminus A$, we may ask the following question:
Does there exist $f\in L^2(\R^d)$ which is not $X$-connected? If the answer is always ``yes'', this would imply that, for the considered window classes, phase retrieval is equivalent to $V_gg(x,\omega)$ being non-zero for almost all $x,\omega\in\R^d$. Otherwise, any set $X$ for which the answer is ``no'' could give us a clue towards constructing a window that does phase retrieval despite satisfying $V_gg\vert_M\equiv 0$ on a subset $M$ of $\R^{2d}$ of positive measure.
However, we are currently unable to answer this question.
\end{remark}

We begin by revisiting an example of a window that is already known to do phase retrieval.

\begin{example}
Let $\gamma_a\in L^2(\R)$ be defined by $\gamma_a(t):=e^{-a\pi t^2}$ for some $a>0$. Then $g_a$ does phase retrieval.
\end{example}

As mentioned before, this result is already well-known, but can also be obtained by Theorem \ref{thm:Rdexpdec}. The theorem therefore covers one of the most prominent windows for phase retrieval (or generally in time-frequency analysis).
The same applies (partially) to the next example.

\begin{example}
Let $\eta_a:\R\to\R$ be defined by $\eta_a(t):=e^{-at}\chi_{(0,\infty)}(t)$ for some $a>0$. Then $\eta_a$ does phase retrieval. Moreover, the two-sided version $\widetilde{\eta}_a:=\eta_a+\mathcal{I}\eta_a$, i.e. $\widetilde{\eta}_a(t)=e^{-a|t|}$ for almost every $t\in\R$, does phase retrieval as well.
\end{example}

The fact that $\eta_a$ does phase retrieval is also known (cf. again the beginning of this section), but can be obtained through the one-sided version of Theorem \ref{thm:Rdexpdec} as well.
On the other hand, the statement about $\widetilde{\eta}_a$ appears to be new. It becomes clear that $\eta_a$ doing phase retrieval doesn't ``depend'' on it being one-sided.
More generally, Theorem \ref{thm:Rdexpdec} shows immediately that \textit{every} $g\in L^{2,\exp}(\R^d)$ that vanishes almost nowhere, does phase retrieval. Even though the exponential decay gives some restrictions, we still obtain a large variety of windows doing phase retrieval.\\

Conversely, we may also apply our results to the one-sided versions of (known) two-sided phase retrieval windows.

\begin{example}
Let $c_a:\R\to\R$ be defined by $c_a(t):=\frac{1}{a+2\pi it}$ for some $a>0$. Then $c_a$ does phase retrieval. Moreover, the one-sided version $\tilde{c}_a:=c_a\cdot\chi_{(0,\infty)}$ does phase retrieval as well.
\end{example}

As mentioned before, phase retrieval for $c_a$ is known from \cite{GJM} and can \textit{not} be obtained using the results of this section. However, phase retrieval for $\tilde{c}_a$ seems to be another new result.
More generally, Theorem \ref{thm:Rdexpdec} yields phase retrieval for \textit{every} $g\in L^{2,+}(\R)$ that vanishes almost nowhere on $(0,\infty)$.\\

Finally, we consider some compactly supported windows. These are obviously inadequate for global phase retrieval, but can still be used for phase retrieval on special signal classes.

\begin{example}\label{ex:gab}
Let $g_{a,b}:=\chi_{[a,b]}$ for some $a,b\in\R$ satisfying $a<b$. Then, a signal $f\in L^2(\R^d)$ is phase retrievable with respect to $g_{a,b}$ if and only if $f$ has no hole of length $b-a$, i.e. if the implication
\begin{equation*}
f\vert_{(t,t+b-a)}\equiv 0 ~\text{a.e.}\quad\quad\Rightarrow\quad\quad f\vert_{(-\infty,t)}\equiv 0~\text{a.e.}\quad\text{or}\quad f\vert_{(t+b-a,\infty)}\equiv 0~\text{a.e.}
\end{equation*}
holds true for every $t\in\R$.
\end{example}

This follows from the fact that obviously, it holds $X=(-(b-a),b-a)$ in this case and from the fact that a signal $f$ is $X$-connected if and only if it has no holes of length $b-a$. As in the discrete setting, we obtain the intuition that $g$ has to be ``long enough'' to bridge over possible holes of $f$.\\
It is possible to ``thin out'' the support of $g$ even further.

\begin{example}
Let $g_{a_1,a_2,b_1,b_2}:=\chi_{[a_1,a_2]}+\chi_{[b_1,b_2]}$ for some $a_1,a_2,b_1,b_2\in\R$ satisfying\linebreak $a_1<a_2<b_1<b_2$ as well as $b_1-a_2<\max\{a_2-a_1,b_2-b_1\}$.
Then, a signal $f\in L^2(\R^d)$ is phase retrievable with respect to $g_{a_1,a_2,b_1,b_2}$ if and only if $f$ has no hole of length $b_2-a_1$ in the sense of Example \ref{ex:gab}.
\end{example}

In order to see this, one needs to accept $X=(-(b_2-a_1),b_2-a_1)$, which follows from the upper bound on $b_1-a_2$.\\
A variation of this example also shows that things can get out of control very easily: When considering $g_{a_1,a_2,b_1,b_2}$ for the case $b_1-a_2>\max\{a_2-a_1,b_2-b_1\}$, we obtain
\begin{equation*}
X=(-(b_2-a_1),-(b_1-a_2))\cup(-m,m)\cup(b_1-a_2,b_2-a_1),
\end{equation*}
where $m:=\max\{a_2-a_1,b_2-b_1\}$. Classifying all $X$-connected signals appears to be a difficult task in this case. (Note that we can say at least that ``no hole of length $m$'' is a sufficient condition and ``no hole of length $b_2-a_1$'' is a necessary condition for $X$-connectivity.)

\section{Concluding remarks}

The goal of this paper was to provide useful criteria for STFT phase retrieval on $\Z$, $\Z_d$ and $\R^d$. Coarsely summarized, we proved a full characterization of phase retrieval for the following combinations of windows and signals:

\begin{enumerate}[1)]
\item $g\in\ell^2(\Z)$ one-sided, $f\in\ell^2(\Z)$ arbitrary
\item $g\in\ell^2(\Z)$ of exponential decay, $f\in\ell^2(\Z)$ arbitrary
\item $g\in\C^d$ of length $L+1\leq\frac{d}{2}$ drawn randomly, $f\in\C^d$ arbitrary (characterization holds with probability 1)
\item $g\in\C^d$ of length $L+1\leq\frac{d}{2}$ fixed, $f\in\C^d$ with at least $L$ consecutive zeros
\item $g\in L^2(\R)$ one-sided, $f\in L^2(\R)$ arbitrary
\item $g\in L^2(\R^d)$ of exponential decay, $f\in L^2(\R^d)$ arbitrary
\end{enumerate}

Accordingly, we usually proved characterizations for phase retrieval that are valid for the whole signal class $L^2(G)$, with scenario 4) being the only exception.
This corresponds to identifying both the possibilities and the limitations of phase retrieval for a single given window (i.e. a given measurement method), which appears to be highly useful in applications.\\

Large similarities between the different settings are evident (especially when comparing $G=\Z$ and $G=\R^d$) with the common theme within all the different scenarios being the fact that the characterizations only depend on the support of window and signal.
This leads to conditions that are usually easy to understand and to check (at least in the discrete settings) and showcase the fact that disconnectedness (in one way or another) is (up to only a few pathological counterexamples in the cyclic case) essentially the only obstacle for STFT phase retrieval.

\section*{Acknowledgements}
The author would like to thank Hartmut Führ for fruitful discussions on the topic as well as careful reading of the manuscript.

\bibliography{PRSTFT.bib}
\bibliographystyle{plain}
\vspace{0.5cm}

\end{document}